\documentclass[a4paper,12pt]{article}
\usepackage{fullpage}
\usepackage{amsthm}
\usepackage{amsfonts}
\usepackage{amsmath}
\usepackage{amscd}
\usepackage{amssymb}
\usepackage[svgnames]{pstricks}
\usepackage{pstricks,pst-node,pst-tree}
\usepackage{stmaryrd}
\usepackage{graphics}
\theoremstyle{plain}
\newtheorem{thm}{Theorem}[section]
\newtheorem{lem}[thm]{Lemma}
\newtheorem{prop}[thm]{Proposition}

\newtheorem{prob}[thm]{Problem}

\theoremstyle{definition}

\theoremstyle{remark}
\newtheorem*{rem}{Remark}

\newtheorem{case}{Case}
\newcommand{\N}{\mathbb{N}}
\newcommand{\Z}{\mathbb{Z}}

\newcommand{\Zn}{\Z/n\Z}
\newcommand{\m}{\mathfrak{m}}

\newcommand{\twoheadlongrightarrow}{\relbar\joinrel\twoheadrightarrow}

\DeclareMathOperator*{\rad}{rad}
\renewcommand{\leq}{\leqslant}
\renewcommand{\geq}{\geqslant}

\newcommand{\circul}{\mathbf{Circ}}
\newcommand{\univ}{\mathrm{US}}
\newcommand{\ST}{\mathrm{ST}}
\newcommand{\PT}{\mathrm{PT}}
\newcommand{\IAP}{\mathrm{IAP}}
\newcommand{\AP}{\mathrm{AP}}
\newcommand{\DAT}{\mathrm{DAT}}
\newcommand{\IDAO}{\mathrm{IDAO}}
\begin{document}
\title{A universal sequence of integers generating balanced Steinhaus figures modulo an odd number}
\author{\textsc{Jonathan Chappelon}\textsuperscript{a,b,c}\\[2ex]
\textsuperscript{a}Univ Lille Nord de France, F-59000 Lille, France\\
\textsuperscript{b}ULCO, LMPA J.~Liouville, B.P. 699, F-62228 Calais, France\\
\textsuperscript{c}CNRS, FR 2956, France\\[2ex]
jonathan.chappelon@lmpa.univ-littoral.fr}
\date{April 15, 2010}
\maketitle
\begin{abstract}
In this paper, we partially solve an open problem, due to J.C.~Molluzzo in 1976, on the existence of balanced Steinhaus triangles modulo a positive integer $n$, that are Steinhaus triangles containing all the elements of $\Zn$ with the same multiplicity. For every odd number $n$, we build an orbit in $\Zn$, by the linear cellular automaton generating the Pascal triangle modulo $n$, which contains infinitely many balanced Steinhaus triangles. This orbit, in $\Zn$, is obtained from an integer sequence called the universal sequence. We show that there exist balanced Steinhaus triangles for at least $2/3$ of the admissible sizes, in the case where $n$ is an odd prime power. Other balanced Steinhaus figures, such as Steinhaus trapezoids, generalized Pascal triangles, Pascal trapezoids or lozenges, also appear in the orbit of the universal sequence modulo $n$ odd. We prove the existence of balanced generalized Pascal triangles for at least $2/3$ of the admissible sizes, in the case where $n$ is an odd prime power, and the existence of balanced lozenges for all admissible sizes, in the case where $n$ is a square-free odd number.\\[2ex]
\end{abstract}

\section{Introduction}

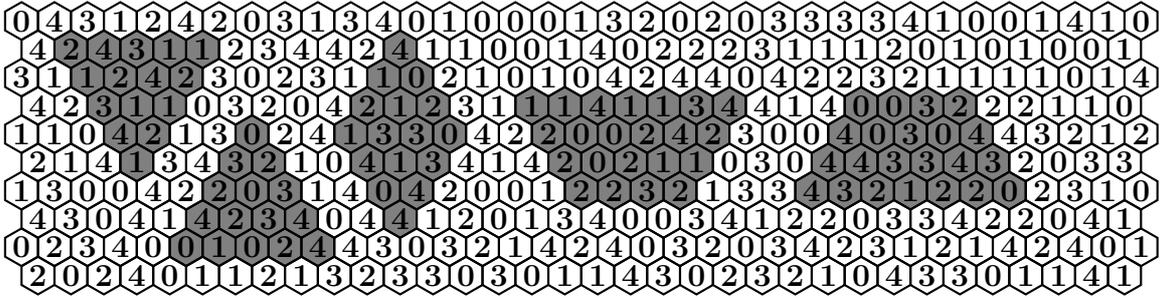
\begin{figure}[!h]
\begin{center}
\begin{pspicture}(15.1554446,3.875)
\multips(0.649519053,3.375)(0.433012702,0){5}{\pspolygon*[linecolor=Gray](0,0)(0.216506351,0.125)(0.433012702,0)(0.433012702,-0.25)(0.216506351,-0.375)(0,-0.25)}
\multips(0.866025404,3)(0.433012702,0){4}{\pspolygon*[linecolor=Gray](0,0)(0.216506351,0.125)(0.433012702,0)(0.433012702,-0.25)(0.216506351,-0.375)(0,-0.25)}
\multips(1.08253175,2.625)(0.433012702,0){3}{\pspolygon*[linecolor=Gray](0,0)(0.216506351,0.125)(0.433012702,0)(0.433012702,-0.25)(0.216506351,-0.375)(0,-0.25)}
\multips(1.29903811,2.25)(0.433012702,0){2}{\pspolygon*[linecolor=Gray](0,0)(0.216506351,0.125)(0.433012702,0)(0.433012702,-0.25)(0.216506351,-0.375)(0,-0.25)}
\multips(1.51554446,1.875)(0.433012702,0){1}{\pspolygon*[linecolor=Gray](0,0)(0.216506351,0.125)(0.433012702,0)(0.433012702,-0.25)(0.216506351,-0.375)(0,-0.25)}

\multips(2.16506351,0.75)(0.433012702,0){5}{\pspolygon*[linecolor=Gray](0,0)(0.216506351,0.125)(0.433012702,0)(0.433012702,-0.25)(0.216506351,-0.375)(0,-0.25)}
\multips(2.38156986,1.125)(0.433012702,0){4}{\pspolygon*[linecolor=Gray](0,0)(0.216506351,0.125)(0.433012702,0)(0.433012702,-0.25)(0.216506351,-0.375)(0,-0.25)}
\multips(2.59807621,1.5)(0.433012702,0){3}{\pspolygon*[linecolor=Gray](0,0)(0.216506351,0.125)(0.433012702,0)(0.433012702,-0.25)(0.216506351,-0.375)(0,-0.25)}
\multips(2.81458256,1.875)(0.433012702,0){2}{\pspolygon*[linecolor=Gray](0,0)(0.216506351,0.125)(0.433012702,0)(0.433012702,-0.25)(0.216506351,-0.375)(0,-0.25)}
\multips(3.03108891,2.25)(0.433012702,0){1}{\pspolygon*[linecolor=Gray](0,0)(0.216506351,0.125)(0.433012702,0)(0.433012702,-0.25)(0.216506351,-0.375)(0,-0.25)}

\multips(0.433012702,-0.75)(0,0){1}{
\multips(4.54663337,1.875)(0.433012702,0){1}{\pspolygon*[linecolor=Gray](0,0)(0.216506351,0.125)(0.433012702,0)(0.433012702,-0.25)(0.216506351,-0.375)(0,-0.25)}
\multips(4.33012702,2.25)(0.433012702,0){2}{\pspolygon*[linecolor=Gray](0,0)(0.216506351,0.125)(0.433012702,0)(0.433012702,-0.25)(0.216506351,-0.375)(0,-0.25)}
\multips(4.11362067,2.625)(0.433012702,0){3}{\pspolygon*[linecolor=Gray](0,0)(0.216506351,0.125)(0.433012702,0)(0.433012702,-0.25)(0.216506351,-0.375)(0,-0.25)}
\multips(3.89711432,3)(0.433012702,0){4}{\pspolygon*[linecolor=Gray](0,0)(0.216506351,0.125)(0.433012702,0)(0.433012702,-0.25)(0.216506351,-0.375)(0,-0.25)}
\multips(4.11362067,3.375)(0.433012702,0){3}{\pspolygon*[linecolor=Gray](0,0)(0.216506351,0.125)(0.433012702,0)(0.433012702,-0.25)(0.216506351,-0.375)(0,-0.25)}
\multips(4.33012702,3.75)(0.433012702,0){2}{\pspolygon*[linecolor=Gray](0,0)(0.216506351,0.125)(0.433012702,0)(0.433012702,-0.25)(0.216506351,-0.375)(0,-0.25)}
\multips(4.54663337,4.125)(0.433012702,0){1}{\pspolygon*[linecolor=Gray](0,0)(0.216506351,0.125)(0.433012702,0)(0.433012702,-0.25)(0.216506351,-0.375)(0,-0.25)}}

\multips(6.71169688,2.625)(0.433012702,0){7}{\pspolygon*[linecolor=Gray](0,0)(0.216506351,0.125)(0.433012702,0)(0.433012702,-0.25)(0.216506351,-0.375)(0,-0.25)}
\multips(6.92820323,2.25)(0.433012702,0){6}{\pspolygon*[linecolor=Gray](0,0)(0.216506351,0.125)(0.433012702,0)(0.433012702,-0.25)(0.216506351,-0.375)(0,-0.25)}
\multips(7.14470958,1.875)(0.433012702,0){5}{\pspolygon*[linecolor=Gray](0,0)(0.216506351,0.125)(0.433012702,0)(0.433012702,-0.25)(0.216506351,-0.375)(0,-0.25)}
\multips(7.36121593,1.5)(0.433012702,0){4}{\pspolygon*[linecolor=Gray](0,0)(0.216506351,0.125)(0.433012702,0)(0.433012702,-0.25)(0.216506351,-0.375)(0,-0.25)}

\multips(11.0418239,2.625)(0.433012702,0){4}{\pspolygon*[linecolor=Gray](0,0)(0.216506351,0.125)(0.433012702,0)(0.433012702,-0.25)(0.216506351,-0.375)(0,-0.25)}
\multips(10.8253175,2.25)(0.433012702,0){5}{\pspolygon*[linecolor=Gray](0,0)(0.216506351,0.125)(0.433012702,0)(0.433012702,-0.25)(0.216506351,-0.375)(0,-0.25)}
\multips(10.6088112,1.875)(0.433012702,0){6}{\pspolygon*[linecolor=Gray](0,0)(0.216506351,0.125)(0.433012702,0)(0.433012702,-0.25)(0.216506351,-0.375)(0,-0.25)}
\multips(10.3923048,1.5)(0.433012702,0){7}{\pspolygon*[linecolor=Gray](0,0)(0.216506351,0.125)(0.433012702,0)(0.433012702,-0.25)(0.216506351,-0.375)(0,-0.25)}

\multips(0,3.375)(0,-0.75){5}{\multips(0,0)(0.433012702,0){35}{\pspolygon(0,0.125)(0,0.375)(0.216506351,0.5)(0.433012702,0.375)(0.433012702,0.125)(0.216506351,0)}}
\multips(0.216506351,3)(0,-0.75){5}{\multips(0,0)(0.433012702,0){34}{\pspolygon(0,0.125)(0,0.375)(0.216506351,0.5)(0.433012702,0.375)(0.433012702,0.125)(0.216506351,0)}}

\rput(0.216506,3.625){$\mathbf{0}$}
\rput(0.649519,3.625){$\mathbf{4}$}
\rput(1.08253,3.625){$\mathbf{3}$}
\rput(1.51554,3.625){$\mathbf{1}$}
\rput(1.94856,3.625){$\mathbf{2}$}
\rput(2.38157,3.625){$\mathbf{4}$}
\rput(2.81458,3.625){$\mathbf{2}$}
\rput(3.2476,3.625){$\mathbf{0}$}
\rput(3.68061,3.625){$\mathbf{3}$}
\rput(4.11362,3.625){$\mathbf{1}$}
\rput(4.54663,3.625){$\mathbf{3}$}
\rput(4.97965,3.625){$\mathbf{4}$}
\rput(5.41266,3.625){$\mathbf{0}$}
\rput(5.84567,3.625){$\mathbf{1}$}
\rput(6.27868,3.625){$\mathbf{0}$}
\rput(6.7117,3.625){$\mathbf{0}$}
\rput(7.14471,3.625){$\mathbf{0}$}
\rput(7.57772,3.625){$\mathbf{1}$}
\rput(8.01073,3.625){$\mathbf{3}$}
\rput(8.44375,3.625){$\mathbf{2}$}
\rput(8.87676,3.625){$\mathbf{0}$}
\rput(9.30977,3.625){$\mathbf{2}$}
\rput(9.74279,3.625){$\mathbf{0}$}
\rput(10.1758,3.625){$\mathbf{3}$}
\rput(10.6088,3.625){$\mathbf{3}$}
\rput(11.0418,3.625){$\mathbf{3}$}
\rput(11.4748,3.625){$\mathbf{3}$}
\rput(11.9078,3.625){$\mathbf{4}$}
\rput(12.3409,3.625){$\mathbf{1}$}
\rput(12.7739,3.625){$\mathbf{0}$}
\rput(13.2069,3.625){$\mathbf{0}$}
\rput(13.6399,3.625){$\mathbf{1}$}
\rput(14.0729,3.625){$\mathbf{4}$}
\rput(14.5059,3.625){$\mathbf{1}$}
\rput(14.9389,3.625){$\mathbf{0}$}

\rput(0.433013,3.25){$\mathbf{4}$}
\rput(0.866025,3.25){$\mathbf{2}$}
\rput(1.29904,3.25){$\mathbf{4}$}
\rput(1.73205,3.25){$\mathbf{3}$}
\rput(2.16506,3.25){$\mathbf{1}$}
\rput(2.59808,3.25){$\mathbf{1}$}
\rput(3.03109,3.25){$\mathbf{2}$}
\rput(3.4641,3.25){$\mathbf{3}$}
\rput(3.89711,3.25){$\mathbf{4}$}
\rput(4.33013,3.25){$\mathbf{4}$}
\rput(4.76314,3.25){$\mathbf{2}$}
\rput(5.19615,3.25){$\mathbf{4}$}
\rput(5.62917,3.25){$\mathbf{1}$}
\rput(6.06218,3.25){$\mathbf{1}$}
\rput(6.49519,3.25){$\mathbf{0}$}
\rput(6.9282,3.25){$\mathbf{0}$}
\rput(7.36122,3.25){$\mathbf{1}$}
\rput(7.79423,3.25){$\mathbf{4}$}
\rput(8.22724,3.25){$\mathbf{0}$}
\rput(8.66025,3.25){$\mathbf{2}$}
\rput(9.09327,3.25){$\mathbf{2}$}
\rput(9.52628,3.25){$\mathbf{2}$}
\rput(9.95929,3.25){$\mathbf{3}$}
\rput(10.3923,3.25){$\mathbf{1}$}
\rput(10.8253,3.25){$\mathbf{1}$}
\rput(11.2583,3.25){$\mathbf{1}$}
\rput(11.6913,3.25){$\mathbf{2}$}
\rput(12.1244,3.25){$\mathbf{0}$}
\rput(12.5574,3.25){$\mathbf{1}$}
\rput(12.9904,3.25){$\mathbf{0}$}
\rput(13.4234,3.25){$\mathbf{1}$}
\rput(13.8564,3.25){$\mathbf{0}$}
\rput(14.2894,3.25){$\mathbf{0}$}
\rput(14.7224,3.25){$\mathbf{1}$}

\rput(0.216506,2.875){$\mathbf{3}$}
\rput(0.649519,2.875){$\mathbf{1}$}
\rput(1.08253,2.875){$\mathbf{1}$}
\rput(1.51554,2.875){$\mathbf{2}$}
\rput(1.94856,2.875){$\mathbf{4}$}
\rput(2.38157,2.875){$\mathbf{2}$}
\rput(2.81458,2.875){$\mathbf{3}$}
\rput(3.2476,2.875){$\mathbf{0}$}
\rput(3.68061,2.875){$\mathbf{2}$}
\rput(4.11362,2.875){$\mathbf{3}$}
\rput(4.54663,2.875){$\mathbf{1}$}
\rput(4.97965,2.875){$\mathbf{1}$}
\rput(5.41266,2.875){$\mathbf{0}$}
\rput(5.84567,2.875){$\mathbf{2}$}
\rput(6.27868,2.875){$\mathbf{1}$}
\rput(6.7117,2.875){$\mathbf{0}$}
\rput(7.14471,2.875){$\mathbf{1}$}
\rput(7.57772,2.875){$\mathbf{0}$}
\rput(8.01073,2.875){$\mathbf{4}$}
\rput(8.44375,2.875){$\mathbf{2}$}
\rput(8.87676,2.875){$\mathbf{4}$}
\rput(9.30977,2.875){$\mathbf{4}$}
\rput(9.74279,2.875){$\mathbf{0}$}
\rput(10.1758,2.875){$\mathbf{4}$}
\rput(10.6088,2.875){$\mathbf{2}$}
\rput(11.0418,2.875){$\mathbf{2}$}
\rput(11.4748,2.875){$\mathbf{3}$}
\rput(11.9078,2.875){$\mathbf{2}$}
\rput(12.3409,2.875){$\mathbf{1}$}
\rput(12.7739,2.875){$\mathbf{1}$}
\rput(13.2069,2.875){$\mathbf{1}$}
\rput(13.6399,2.875){$\mathbf{1}$}
\rput(14.0729,2.875){$\mathbf{0}$}
\rput(14.5059,2.875){$\mathbf{1}$}
\rput(14.9389,2.875){$\mathbf{4}$}

\rput(0.433013,2.5){$\mathbf{4}$}
\rput(0.866025,2.5){$\mathbf{2}$}
\rput(1.29904,2.5){$\mathbf{3}$}
\rput(1.73205,2.5){$\mathbf{1}$}
\rput(2.16506,2.5){$\mathbf{1}$}
\rput(2.59808,2.5){$\mathbf{0}$}
\rput(3.03109,2.5){$\mathbf{3}$}
\rput(3.4641,2.5){$\mathbf{2}$}
\rput(3.89711,2.5){$\mathbf{0}$}
\rput(4.33013,2.5){$\mathbf{4}$}
\rput(4.76314,2.5){$\mathbf{2}$}
\rput(5.19615,2.5){$\mathbf{1}$}
\rput(5.62917,2.5){$\mathbf{2}$}
\rput(6.06218,2.5){$\mathbf{3}$}
\rput(6.49519,2.5){$\mathbf{1}$}
\rput(6.9282,2.5){$\mathbf{1}$}
\rput(7.36122,2.5){$\mathbf{1}$}
\rput(7.79423,2.5){$\mathbf{4}$}
\rput(8.22724,2.5){$\mathbf{1}$}
\rput(8.66025,2.5){$\mathbf{1}$}
\rput(9.09327,2.5){$\mathbf{3}$}
\rput(9.52628,2.5){$\mathbf{4}$}
\rput(9.95929,2.5){$\mathbf{4}$}
\rput(10.3923,2.5){$\mathbf{1}$}
\rput(10.8253,2.5){$\mathbf{4}$}
\rput(11.2583,2.5){$\mathbf{0}$}
\rput(11.6913,2.5){$\mathbf{0}$}
\rput(12.1244,2.5){$\mathbf{3}$}
\rput(12.5574,2.5){$\mathbf{2}$}
\rput(12.9904,2.5){$\mathbf{2}$}
\rput(13.4234,2.5){$\mathbf{2}$}
\rput(13.8564,2.5){$\mathbf{1}$}
\rput(14.2894,2.5){$\mathbf{1}$}
\rput(14.7224,2.5){$\mathbf{0}$}

\rput(0.216506,2.125){$\mathbf{1}$}
\rput(0.649519,2.125){$\mathbf{1}$}
\rput(1.08253,2.125){$\mathbf{0}$}
\rput(1.51554,2.125){$\mathbf{4}$}
\rput(1.94856,2.125){$\mathbf{2}$}
\rput(2.38157,2.125){$\mathbf{1}$}
\rput(2.81458,2.125){$\mathbf{3}$}
\rput(3.2476,2.125){$\mathbf{0}$}
\rput(3.68061,2.125){$\mathbf{2}$}
\rput(4.11362,2.125){$\mathbf{4}$}
\rput(4.54663,2.125){$\mathbf{1}$}
\rput(4.97965,2.125){$\mathbf{3}$}
\rput(5.41266,2.125){$\mathbf{3}$}
\rput(5.84567,2.125){$\mathbf{0}$}
\rput(6.27868,2.125){$\mathbf{4}$}
\rput(6.7117,2.125){$\mathbf{2}$}
\rput(7.14471,2.125){$\mathbf{2}$}
\rput(7.57772,2.125){$\mathbf{0}$}
\rput(8.01073,2.125){$\mathbf{0}$}
\rput(8.44375,2.125){$\mathbf{2}$}
\rput(8.87676,2.125){$\mathbf{4}$}
\rput(9.30977,2.125){$\mathbf{2}$}
\rput(9.74279,2.125){$\mathbf{3}$}
\rput(10.1758,2.125){$\mathbf{0}$}
\rput(10.6088,2.125){$\mathbf{0}$}
\rput(11.0418,2.125){$\mathbf{4}$}
\rput(11.4748,2.125){$\mathbf{0}$}
\rput(11.9078,2.125){$\mathbf{3}$}
\rput(12.3409,2.125){$\mathbf{0}$}
\rput(12.7739,2.125){$\mathbf{4}$}
\rput(13.2069,2.125){$\mathbf{4}$}
\rput(13.6399,2.125){$\mathbf{3}$}
\rput(14.0729,2.125){$\mathbf{2}$}
\rput(14.5059,2.125){$\mathbf{1}$}
\rput(14.9389,2.125){$\mathbf{2}$}

\rput(0.433013,1.75){$\mathbf{2}$}
\rput(0.866025,1.75){$\mathbf{1}$}
\rput(1.29904,1.75){$\mathbf{4}$}
\rput(1.73205,1.75){$\mathbf{1}$}
\rput(2.16506,1.75){$\mathbf{3}$}
\rput(2.59808,1.75){$\mathbf{4}$}
\rput(3.03109,1.75){$\mathbf{3}$}
\rput(3.4641,1.75){$\mathbf{2}$}
\rput(3.89711,1.75){$\mathbf{1}$}
\rput(4.33013,1.75){$\mathbf{0}$}
\rput(4.76314,1.75){$\mathbf{4}$}
\rput(5.19615,1.75){$\mathbf{1}$}
\rput(5.62917,1.75){$\mathbf{3}$}
\rput(6.06218,1.75){$\mathbf{4}$}
\rput(6.49519,1.75){$\mathbf{1}$}
\rput(6.9282,1.75){$\mathbf{4}$}
\rput(7.36122,1.75){$\mathbf{2}$}
\rput(7.79423,1.75){$\mathbf{0}$}
\rput(8.22724,1.75){$\mathbf{2}$}
\rput(8.66025,1.75){$\mathbf{1}$}
\rput(9.09327,1.75){$\mathbf{1}$}
\rput(9.52628,1.75){$\mathbf{0}$}
\rput(9.95929,1.75){$\mathbf{3}$}
\rput(10.3923,1.75){$\mathbf{0}$}
\rput(10.8253,1.75){$\mathbf{4}$}
\rput(11.2583,1.75){$\mathbf{4}$}
\rput(11.6913,1.75){$\mathbf{3}$}
\rput(12.1244,1.75){$\mathbf{3}$}
\rput(12.5574,1.75){$\mathbf{4}$}
\rput(12.9904,1.75){$\mathbf{3}$}
\rput(13.4234,1.75){$\mathbf{2}$}
\rput(13.8564,1.75){$\mathbf{0}$}
\rput(14.2894,1.75){$\mathbf{3}$}
\rput(14.7224,1.75){$\mathbf{3}$}

\rput(0.216506,1.375){$\mathbf{1}$}
\rput(0.649519,1.375){$\mathbf{3}$}
\rput(1.08253,1.375){$\mathbf{0}$}
\rput(1.51554,1.375){$\mathbf{0}$}
\rput(1.94856,1.375){$\mathbf{4}$}
\rput(2.38157,1.375){$\mathbf{2}$}
\rput(2.81458,1.375){$\mathbf{2}$}
\rput(3.2476,1.375){$\mathbf{0}$}
\rput(3.68061,1.375){$\mathbf{3}$}
\rput(4.11362,1.375){$\mathbf{1}$}
\rput(4.54663,1.375){$\mathbf{4}$}
\rput(4.97965,1.375){$\mathbf{0}$}
\rput(5.41266,1.375){$\mathbf{4}$}
\rput(5.84567,1.375){$\mathbf{2}$}
\rput(6.27868,1.375){$\mathbf{0}$}
\rput(6.7117,1.375){$\mathbf{0}$}
\rput(7.14471,1.375){$\mathbf{1}$}
\rput(7.57772,1.375){$\mathbf{2}$}
\rput(8.01073,1.375){$\mathbf{2}$}
\rput(8.44375,1.375){$\mathbf{3}$}
\rput(8.87676,1.375){$\mathbf{2}$}
\rput(9.30977,1.375){$\mathbf{1}$}
\rput(9.74279,1.375){$\mathbf{3}$}
\rput(10.1758,1.375){$\mathbf{3}$}
\rput(10.6088,1.375){$\mathbf{4}$}
\rput(11.0418,1.375){$\mathbf{3}$}
\rput(11.4748,1.375){$\mathbf{2}$}
\rput(11.9078,1.375){$\mathbf{1}$}
\rput(12.3409,1.375){$\mathbf{2}$}
\rput(12.7739,1.375){$\mathbf{2}$}
\rput(13.2069,1.375){$\mathbf{0}$}
\rput(13.6399,1.375){$\mathbf{2}$}
\rput(14.0729,1.375){$\mathbf{3}$}
\rput(14.5059,1.375){$\mathbf{1}$}
\rput(14.9389,1.375){$\mathbf{0}$}

\rput(0.433013,1){$\mathbf{4}$}
\rput(0.866025,1){$\mathbf{3}$}
\rput(1.29904,1){$\mathbf{0}$}
\rput(1.73205,1){$\mathbf{4}$}
\rput(2.16506,1){$\mathbf{1}$}
\rput(2.59808,1){$\mathbf{4}$}
\rput(3.03109,1){$\mathbf{2}$}
\rput(3.4641,1){$\mathbf{3}$}
\rput(3.89711,1){$\mathbf{4}$}
\rput(4.33013,1){$\mathbf{0}$}
\rput(4.76314,1){$\mathbf{4}$}
\rput(5.19615,1){$\mathbf{4}$}
\rput(5.62917,1){$\mathbf{1}$}
\rput(6.06218,1){$\mathbf{2}$}
\rput(6.49519,1){$\mathbf{0}$}
\rput(6.9282,1){$\mathbf{1}$}
\rput(7.36122,1){$\mathbf{3}$}
\rput(7.79423,1){$\mathbf{4}$}
\rput(8.22724,1){$\mathbf{0}$}
\rput(8.66025,1){$\mathbf{0}$}
\rput(9.09327,1){$\mathbf{3}$}
\rput(9.52628,1){$\mathbf{4}$}
\rput(9.95929,1){$\mathbf{1}$}
\rput(10.3923,1){$\mathbf{2}$}
\rput(10.8253,1){$\mathbf{2}$}
\rput(11.2583,1){$\mathbf{0}$}
\rput(11.6913,1){$\mathbf{3}$}
\rput(12.1244,1){$\mathbf{3}$}
\rput(12.5574,1){$\mathbf{4}$}
\rput(12.9904,1){$\mathbf{2}$}
\rput(13.4234,1){$\mathbf{2}$}
\rput(13.8564,1){$\mathbf{0}$}
\rput(14.2894,1){$\mathbf{4}$}
\rput(14.7224,1){$\mathbf{1}$}

\rput(0.216506,0.625){$\mathbf{0}$}
\rput(0.649519,0.625){$\mathbf{2}$}
\rput(1.08253,0.625){$\mathbf{3}$}
\rput(1.51554,0.625){$\mathbf{4}$}
\rput(1.94856,0.625){$\mathbf{0}$}
\rput(2.38157,0.625){$\mathbf{0}$}
\rput(2.81458,0.625){$\mathbf{1}$}
\rput(3.2476,0.625){$\mathbf{0}$}
\rput(3.68061,0.625){$\mathbf{2}$}
\rput(4.11362,0.625){$\mathbf{4}$}
\rput(4.54663,0.625){$\mathbf{4}$}
\rput(4.97965,0.625){$\mathbf{3}$}
\rput(5.41266,0.625){$\mathbf{0}$}
\rput(5.84567,0.625){$\mathbf{3}$}
\rput(6.27868,0.625){$\mathbf{2}$}
\rput(6.7117,0.625){$\mathbf{1}$}
\rput(7.14471,0.625){$\mathbf{4}$}
\rput(7.57772,0.625){$\mathbf{2}$}
\rput(8.01073,0.625){$\mathbf{4}$}
\rput(8.44375,0.625){$\mathbf{0}$}
\rput(8.87676,0.625){$\mathbf{3}$}
\rput(9.30977,0.625){$\mathbf{2}$}
\rput(9.74279,0.625){$\mathbf{0}$}
\rput(10.1758,0.625){$\mathbf{3}$}
\rput(10.6088,0.625){$\mathbf{4}$}
\rput(11.0418,0.625){$\mathbf{2}$}
\rput(11.4748,0.625){$\mathbf{3}$}
\rput(11.9078,0.625){$\mathbf{1}$}
\rput(12.3409,0.625){$\mathbf{2}$}
\rput(12.7739,0.625){$\mathbf{1}$}
\rput(13.2069,0.625){$\mathbf{4}$}
\rput(13.6399,0.625){$\mathbf{2}$}
\rput(14.0729,0.625){$\mathbf{4}$}
\rput(14.5059,0.625){$\mathbf{0}$}
\rput(14.9389,0.625){$\mathbf{1}$}

\rput(0.433013,0.25){$\mathbf{2}$}
\rput(0.866025,0.25){$\mathbf{0}$}
\rput(1.29904,0.25){$\mathbf{2}$}
\rput(1.73205,0.25){$\mathbf{4}$}
\rput(2.16506,0.25){$\mathbf{0}$}
\rput(2.59808,0.25){$\mathbf{1}$}
\rput(3.03109,0.25){$\mathbf{1}$}
\rput(3.4641,0.25){$\mathbf{2}$}
\rput(3.89711,0.25){$\mathbf{1}$}
\rput(4.33013,0.25){$\mathbf{3}$}
\rput(4.76314,0.25){$\mathbf{2}$}
\rput(5.19615,0.25){$\mathbf{3}$}
\rput(5.62917,0.25){$\mathbf{3}$}
\rput(6.06218,0.25){$\mathbf{0}$}
\rput(6.49519,0.25){$\mathbf{3}$}
\rput(6.9282,0.25){$\mathbf{0}$}
\rput(7.36122,0.25){$\mathbf{1}$}
\rput(7.79423,0.25){$\mathbf{1}$}
\rput(8.22724,0.25){$\mathbf{4}$}
\rput(8.66025,0.25){$\mathbf{3}$}
\rput(9.09327,0.25){$\mathbf{0}$}
\rput(9.52628,0.25){$\mathbf{2}$}
\rput(9.95929,0.25){$\mathbf{3}$}
\rput(10.3923,0.25){$\mathbf{2}$}
\rput(10.8253,0.25){$\mathbf{1}$}
\rput(11.2583,0.25){$\mathbf{0}$}
\rput(11.6913,0.25){$\mathbf{4}$}
\rput(12.1244,0.25){$\mathbf{3}$}
\rput(12.5574,0.25){$\mathbf{3}$}
\rput(12.9904,0.25){$\mathbf{0}$}
\rput(13.4234,0.25){$\mathbf{1}$}
\rput(13.8564,0.25){$\mathbf{1}$}
\rput(14.2894,0.25){$\mathbf{4}$}
\rput(14.7224,0.25){$\mathbf{1}$}
\end{pspicture}

\caption{\label{fig1}\footnotesize Examples of Steinhaus figures in $\Z/5\Z$: the Steinhaus triangle $\nabla(2,4,3,1,1)$, the Pascal triangle $\Delta(4,2,1,3,0,2,4,1,3)$, the lozenge $\diamondsuit(4,4,2,4,1,1,0)$, the Steinhaus trapezoid $\ST((1,1,4,1,1,3,4),4)$ and the Pascal trapezoid $\PT((2,0,3,3,3,3,4,1,0,0,1,4,1),4)$.}
\end{center}
\end{figure}

Let $n$ be a positive integer and let $\Zn$ denote the finite cyclic group of order $n$. Let $S=(a_j)_{j\in\Z}$ be a doubly infinite sequence of elements in $\Zn$. The \textit{derived sequence} $\partial S$ of $S$ is the sequence obtained by pairwise adding consecutive terms of $S$, that is $\partial S=(a_j+a_{j+1})_{j\in\Z}$. This operation of derivation can be repeated and then, the $i$th derived sequence $\partial^{i}S$ is recursively defined by $\partial^{0} S=S$ and $\partial^{i}S=\partial\partial^{i-1}S$ for all integers $i\geq1$. The sequence of all the iterated derived sequences of $S$ is called the \textit{orbit} $\mathcal{O}_S=(\partial^{i}S)_{i\in\N}$ of $S$. For all $i\in\N$ and all $j\in\Z$, we denote by $a_{i,j}$ the $j$th term of $\partial^iS$. Since $a_{i+1,j}=a_{i,j}+a_{i,j+1}$ by the linear local rule of this cellular automaton, the orbit of $S$ can be seen as the $(\N\times\Z)$-indexed sequence of elements in $\Zn$ defined by
$$
\mathcal{O}_S = \left(a_{i,j}=\sum_{k=0}^{i}\binom{i}{k}a_{j+k}\ \middle|\ i\in\N,j\in\Z\right),
$$
where $\binom{i}{k}$ is the binomial coefficient $\binom{i}{k}=\frac{i!}{(i-k)!k!}$. For every $i\in\N$, the $i$th \textit{row} of $\mathcal{O}_S$ is the sequence $R_i=\partial^iS=(a_{i,j})_{j\in\Z}$ and, for every $j\in\Z$, the $j$th \textit{diagonal} and the $j$th \textit{anti-diagonal} of $\mathcal{O}_S$ are the sequences $D_j=(a_{i,j})_{i\in\N}$ and  $AD_j=(a_{i,j-i})_{i\in\N}$ respectively. Orbits of integer sequences and the canonical projection map $\pi_n:\Z\twoheadlongrightarrow\Zn$ are also considered in this paper. Elementary figures appear in this linear cellular automaton. Examples of them in $\Z/5\Z$ are depicted in Figure~\ref{fig1}.

In this paper, Steinhaus figures are viewed as finite \textit{multisets} in $\Zn$, that are sets in $\Zn$ for which repeated elements are allowed. A finite multiset $M$ in $\Zn$ corresponds to a function $\m_M:\Zn\longrightarrow\N$, the \textit{multiplicity function} associated with $M$, which assigns its multiplicity in $M$ to each element of $\Zn$. The \textit{cardinality} of $M$, denoted by $|M|$, is the number of elements of $M$ counted with multiplicity, that is the non-negative integer $|M|=\sum_{x\in\Zn}\m_M(x)$.

Now, let $S_m=(a_0,\ldots,a_{m-1})$ be a finite sequence of length $m\geq1$ in $\Zn$. The \textit{Steinhaus triangle} $\nabla S_m$ associated with $S_m$ is the collection of all the iterated derived sequences of $S_m$, that is the finite orbit $\nabla S_m=\mathcal{O}_{S_m}=\{S_m,\partial S_m,\ldots,\partial^{m-1}S_m\}$. Namely, it is the multiset in $\Zn$ defined by
$$
\nabla S_m = \left\{ \sum_{k=0}^{i}\binom{i}{k}a_{j+k}\ \middle|\ 0\leq i\leq m-1,0\leq j\leq m-1-i\right\}.
$$
We shall say that the triangle $\nabla S_m$ is of \textit{order} $m$. A Steinhaus triangle of order $m$ has cardinality $\binom{m+1}{2}$. These triangles have been named in honor of H.~Steinhaus, who proposed this construction, for the binary case $\Z/2\Z$, in his book on elementary mathematical problems \cite{Steinhaus1963}.
The \textit{Steinhaus trapezoid} $\ST(S_m,h)$ of order $m$ and of height $h$, with $1\leq h\leq m$, is the collection of the first $h$ derived sequences of $S_m$, that is,
$$
\ST(S_m,h)=\bigcup_{i=0}^{h-1}\partial^{i}S_m=\nabla S_m\setminus\nabla\partial^{h}S_m.
$$
A Steinhaus trapezoid of order $m$ and of height $h$ has cardinality $h(2m-h+1)/2$. Now, let $S_{2m-1}=(a_0,\ldots,a_{2m-2})$ be a finite sequence of length $2m-1\geq1$ in $\Zn$. The \textit{generalized Pascal triangle} (or \textit{Pascal triangle} for short) $\Delta S_{2m-1}$ associated with $S_{2m-1}$ is the triangle of height $m$, built from the top to the base, appearing in the center of the Steinhaus triangle $\nabla S_{2m-1}$. Namely, it is the multiset in $\Zn$ defined by
$$
\Delta S_{2m-1} = \left\{ \sum_{k=0}^{i}\binom{i}{k}a_{m-1-j-k}\ \middle|\ 0\leq j\leq i\leq m-1\right\}.
$$
Obviously, the generalized Pascal triangle associated with the sequence of length $2m-1$ with a $1$ in the middle and $0$ elsewhere corresponds to the first $m$ rows of the standard Pascal triangle modulo $n$. A Pascal triangle of order $2m-1$ has cardinality $\binom{m+1}{2}$. The \textit{Pascal trapezoid} $\PT(S_{2m-1},h)$ of order $2m-1$ and of height $h$ is the collection of the last $h$ rows of the Pascal triangle $\Delta S_{2m-1}$, that is,
$$
\PT(S_{2m-1},h)=\Delta S_{2m-1}\setminus\Delta(a_j)_{h\leq j\leq 2m-h-2}.
$$
A Pascal trapezoid of order $2m-1$ and of height $h$ has cardinality $h(2m-h+1)/2$. Finally, the \textit{lozenge} $\lozenge S_{2m-1}$ associated with the sequence $S_{2m-1}$ is the multiset union of the Pascal triangle $\Delta S_{2m-1}$ and of the Steinhaus triangle $\nabla\partial^{m}S_{2m-1}$. The lozenge $\lozenge S_{2m-1}$ is then the multiset in $\Zn$ defined by
$$
\lozenge S_{2m-1} = \Delta S_{2m-1}\bigcup\nabla\partial^{m}S_{2m-1} = \left\{ \sum_{k=0}^{i+j}\binom{i+j}{k}a_{m-1-j-k}\ \middle|\ 0\leq i,j\leq m-1\right\}.
$$
A lozenge of order $2m-1$ has cardinality $m^2$.

In 1963 \cite{Steinhaus1963}, H.~Steinhaus posed the elementary problem of determining if there exists, for every $m\geq1$ such that $(m+1)m/2$ is even, a binary Steinhaus triangle of order $m$ containing as many $0$'s as $1$'s. This problem was solved, for the first time, by H.~Harborth in 1972 \cite{Harborth1972}. For every $m\equiv0$ or $3\pmod{4}$, he explicitly built at least four such binary Steinhaus triangles of order $m$. Other solutions of the Steinhaus problem appear in the literature \cite{Eliahou2007,Eliahou2004,Eliahou2005}. A generalization of this problem in any finite cyclic group was posed by J.C.~Molluzzo in 1976 \cite{Molluzzo1978}.

A finite multiset $M$ in $\Zn$ is said to be \textit{balanced} if each element of $\Zn$ appears in $M$ with the same multiplicity. Thus, the multiset $M$ is balanced if and only if $\m_M$ is the constant function on $\Zn$ equal to $|M|/n$.

\begin{prob}[Molluzzo,1976]
Let $n$ be a positive integer. For every $m\geq1$ such that the binomial coefficient $\binom{m+1}{2}$ is divisible by $n$, does there exist a balanced Steinhaus triangle of order $m$ in $\Zn$?
\end{prob}

In this paper, for every odd number $n$, we explicitly build balanced Steinhaus triangles of order $m$ in $\Zn$ for every $m\equiv0\pmod{n}$ or $m\equiv-1\pmod{3n}$. This answers in the affirmative Problem~2 of \cite{Chappelon2008}. In \cite{Chappelon2008}, the author completely and positively solved this Molluzzo problem in $\Z/3^k\Z$ for all $k\geq1$. Moreover, for $n$ odd, he showed that there exist at least $\varphi(n)n$ balanced Steinhaus triangles of order $m$ in $\Zn$ for every $m\equiv0$ or $-1\pmod{\varphi(\rad(n))n}$, where $\varphi$ is the \textit{Euler totient function} and $\rad(n)$ is the \textit{radical} of $n$, that is the product of the distinct prime factors of $n$. As observed in \cite{Chappelon2008b}, this problem of Molluzzo does not always admit a positive solution. Indeed, it can be verified, by exhaustive search, that there is no balanced Steinhaus triangle of order $m=5$ in $\Z/15\Z$ or of order $m=6$ in $\Z/21\Z$. Here, we are also interested in the generalization of the Molluzzo problem on each kind of Steinhaus figure defined above, not only on Steinhaus triangles.

\begin{prob}\label{prob2}
Let $n$ be a positive integer. For each kind of Steinhaus figure, do there exist balanced Steinhaus figures in $\Zn$ for all admissible sizes, i.e., for all Steinhaus figures whose cardinality is divisible by $n$? In other words,
\begin{itemize}
\item
for every $m\geq1$ such that $\binom{m+1}{2}$ is divisible by $n$, does there exist a balanced Steinhaus triangle of order $m$?
\item
for every $m\geq1$ and every $h\leq m$ such that $h(2m-h+1)/2$ is divisible by $n$, does there exist a balanced Steinhaus trapezoid of order $m$ and of height $h$?
\item
for every $m\geq1$ such that $\binom{m+1}{2}$ is divisible by $n$, does there exist a balanced Pascal triangle of order $2m-1$?
\item
for every $m\geq1$ and every $h\leq m$ such that $h(2m-h+1)/2$ is divisible by $n$, does there exist a balanced Pascal trapezoid of order $2m-1$ and of height $h$?
\item
for every $m\geq1$ such that $m^2$ is divisible by $n$, does there exist a balanced lozenge of order $2m-1$?
\end{itemize}
\end{prob}

For all positive integers $n$ and $k$ and for all $k$-tuples of elements $A=(a_0,\ldots,a_{k-1})$ and $D=(d_0,\ldots,d_{k-1})$ in $\Zn$, or in $\Z$, the \textit{$k$-interlaced arithmetic progression} $\IAP(A,D)$ is the sequence with first terms $(a_0,\ldots,a_{k-1})$ and with common differences $(d_0,\ldots,d_{k-1})$, that is the doubly infinite sequence $\IAP(A,D)=(a_{j})_{j\in\Z}$ defined by $a_{j_0+jk} = a_{j_0}+jd_{j_0}$, for all $j\in\Z$ and for every $j_0\in\{0,1,\ldots,k-1\}$. For $k=1$, we denote by $\AP(a_0,d_0)$ the arithmetic progression with first element $a_0$ and with common difference $d_0$. 

Let $S=\IAP((0,-1,1),(1,-2,1))$. We shall show that this sequence has the remarkable property that the orbit of its projection $\pi_n(S)$ contains infinitely many balanced Steinhaus figures for every odd $n$. For this reason, we shall call this sequence the \textit{universal sequence} and denote it by $\univ$. The first few terms of $\univ$, where $\mathbf{0}$ is the term of index $0$, are given below:
$$
\univ = \left( \ldots\ldots, -3, -3, 5, -2, -2, 3, -1, -1, 1, 0, \mathbf{0}, -1, 1, 1, -3, 2, 2, -5, 3, 3, -7, \ldots\ldots \right).
$$
The following theorem is the main goal of this article.

\begin{thm}\label{thm0}
Let $n\in\N$ be odd. Then, the orbit of the projection $\pi_n(\univ)$ of the universal sequence in $\Zn$ contains:
\begin{itemize}
\item
balanced Steinhaus triangles of order $m$ for every $m\equiv0\pmod{n}$ or $m\equiv-1\pmod{3n}$. This partially solves the Molluzzo problem for $2/3$ of the admissible orders $m$, in the case where $n$ is an odd prime power.
\item
balanced Steinhaus trapezoids of order $m$ and of height $h$ for every $m\equiv0\pmod{n}$ or $m\equiv-1\pmod{3n}$ and for every $h\equiv m\pmod{n}$ or $h\equiv m+1\pmod{3n}$.
\item
balanced Pascal triangles of order $2m-1$ for every $m\equiv-1\pmod{n}$ or $m\equiv0\pmod{3n}$. This also gives a partial solution of Problem~\ref{prob2} for $2/3$ of the admissible orders $2m-1$, in the case where $n$ is an odd prime power.
\item
balanced Pascal trapezoids of order $2m-1$ and of height $h$ for every $m\equiv-1\pmod{n}$ or $m\equiv0\pmod{3n}$ and for every $h\equiv m+1\pmod{n}$ or $h\equiv m\pmod{3n}$.
\item
balanced lozenges of order $2m-1$ for every $m\equiv0\pmod{n}$. This completely solves Problem~\ref{prob2}, for lozenges, in the case where $n$ is a square-free odd number.
\end{itemize}
\end{thm}

It would be highly desirable to have a similar result for $n$ even, but this is widely open. Here are a few results on Steinhaus figures in the binary case $\Z/2\Z$. The five smallest and the three greatest possible numbers of $1$'s in a binary Steinhaus triangle of fixed size was determined by G.J.~Chang \cite{Chang1983}. H.~Harborth and G.~Hurlbert \cite{Harborth2005} proved that every positive integer is realizable as the number of $1$'s in a generalized binary Pascal triangle, that is, for every natural $k$, there exists a binary sequence $S$ of length $2m_k-1$ such that $\Delta S$ contains exactly $k$ elements equal to $1$. They also determined the minimum value for $m_k$. The maximum number of $1$'s in binary Steinhaus figures (like Steinhaus triangles, generalized Pascal triangles, parallelograms or trapezoids) was studied by M.~Bartsch in her Dissertation \cite{Bartsch1985}. Symmetries in binary Steinhaus triangles and in binary generalized Pascal triangles were explored in \cite{Barb'e2000,Brunat2009}.

This paper is organized as follows. In Section~2, we study doubly arithmetic triangles ($\DAT$ for short) in $\Zn$. These are triangles where all the rows are arithmetic progressions with the same common difference and where all the diagonals are also arithmetic progressions with the same common difference. We show that these triangles constitute a source of balanced multisets in $\Zn$, for $n$ odd, while they are never balanced in $\Zn$, for $n$ even. Moreover, we prove that the orbit associated with the sequence of zeros is the only doubly arithmetic orbit in $\Zn$. In Section~3, interlaced doubly arithmetic orbits, i.e., orbits that are an interlacing of doubly arithmetic structures, are considered. We determine all the interlaced doubly arithmetic orbits in $\Z$ and, in Section~4, we show that the projection of these particular orbits in $\Zn$, for $n$ odd, contains infinitely many balanced Steinhaus figures. This result is refined in Section~5, by considering antisymmetric sequences. In Section~6, a particular case of this antisymmetric refinement leads to the universal sequence $\univ$ and we prove Theorem~\ref{thm0}. Finally, in Section~7, we analyse the results on the generalized Molluzzo problem that we have obtained in this paper and we pose new open problems on the existence of balanced Steinhaus figures in additive cellular automata of dimension $1$ and in the cellular automaton of dimension $2$ where the standard Pascal tetrahedron appears.

\section{DAT: a source of balanced multisets}

For all positive integers $n$ and $m$ and for all elements $a$, $d_1$ and $d_2$ in $\Zn$, the \textit{doubly arithmetic triangle} $\DAT(a,d_1,d_2,m)$ is the triangle of order $m$ in $\Zn$, with first element $a$ and where each diagonal and each row are arithmetic progressions with respective common differences $d_1$ and $d_2$, that is the multiset in $\Zn$ defined by
$$
\DAT(a,d_1,d_2,m) = \left\{ a+id_1+jd_2\ \middle|\ 0\leq i\leq m-1\ ,\ 0\leq j\leq m-1-i \right\}.
$$
In this section, we show that doubly arithmetic triangles constitute a source of balanced multisets in $\Zn$, for $n$ odd. Obviously, we can see that the anti-diagonals of a $\DAT$ are arithmetic progressions with common difference $d_1-d_2$. We begin by determining a necessary condition, on the common differences $d_1$ and $d_2$, to obtain a balanced $\DAT$ in $\Zn$.

\begin{prop}\label{prop1}
Let $n$ be a positive integer and let $a,d_1,d_2\in\Zn$. If the doubly arithmetic triangle $\DAT(a,d_1,d_2,m)$ of order $m\geq1$ is balanced, then its common differences $d_1$, $d_2$ and $d_1-d_2$ are invertible in $\Zn$.
\end{prop}

\begin{proof}
For $n=1$ or $m=1$, it is clear. Suppose now that $n>1$ and $m>1$. Let $\DAT(a,d_1,d_2,m)$ be a doubly arithmetic triangle in $\Zn$ where at least one of the common differences $d_1$, $d_2$ and $d_1-d_2$ is not invertible. Without loss of generality, suppose that it is $d_2$. If not, we can consider the rotations $\DAT(a,d_2,d_1,m)$ or $\DAT(a+(m-1)d_2,-d_2,d_1-d_2,m)$ of $\DAT(a,d_1,d_2,m)$. Let $\delta_1$ and $\delta_2$ be two integers whose respective residue classes modulo $n$ are $d_1$ and $d_2$. We distinguish different cases according to the value of the greatest common divisor of $\delta_1$, $\delta_2$ and $n$.
\begin{case}
If $q=\gcd(\delta_1,\delta_2,n)\neq1$, then we consider the projection map $\pi_q: \Zn\twoheadlongrightarrow\Z/q\Z$. All elements of the triangle $\pi_q(\DAT(a,d_1,d_2,m))=\DAT(\pi_q(a),0,0,m)$ are equal to $\pi_q(a)$. Therefore, the triangle $\DAT(a,d_1,d_2,m)$ is not balanced in $\Z/n\Z$ since its projection in $\Z/q\Z$ is not.
\end{case}
\begin{case}
If $\gcd(\delta_1,\delta_2,n)=1$, then we set $q=\gcd(\delta_2,n)\neq1$ and we consider the projection $\nabla=\pi_q(\DAT(a,d_1,d_2,m))=\DAT(\pi_q(a),\pi_q(d_1),0,m)$ in $\Z/q\Z$, where $\pi_q(d_1)$ is invertible in $\Z/q\Z$. Since the $(kq+l)$th row of $\nabla$ is the constant sequence, of length $m-kq-l+1$, equal to $\pi_q(a)+l\pi_q(d_1)$, for all $l\in\{0,1,\ldots,q-1\}$ and for all $k\in\N$ such that $kq+l\leq m-1$, it follows that we have
$$
\m_{\nabla}(\pi_q(a)) > \m_{\nabla}(\pi_q(a)+\pi_q(d_1)) > \m_{\nabla}(\pi_q(a)+2\pi_q(d_1)) \geq \ldots\ldots \geq \m_{\nabla}(\pi_q(a)+(q-1)\pi_q(d_1)).
$$
Therefore $\nabla$ is not balanced in $\Z/q\Z$ and thus $\DAT(a,d_1,d_2,m)$ is not in $\Zn$.
\end{case}
\end{proof}

\begin{rem}
For $n$ even, there is no balanced $\DAT$ in $\Zn$ since at least one element of $\{d_1,d_2,d_1-d_2\}$ is not invertible in $\Zn$, by the parity of $n$.
\end{rem}

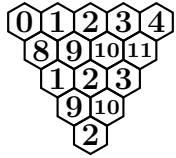
\begin{figure}[!h]
\begin{center}
\begin{pspicture}(2.16506351,2)
\multips(0,1.875)(0.433012702,0){5}{\pspolygon(0,0)(0.216506351,0.125)(0.433012702,0)(0.433012702,-0.25)(0.216506351,-0.375)(0,-0.25)}
\multips(0.216506351,1.5)(0.433012702,0){4}{\pspolygon(0,0)(0.216506351,0.125)(0.433012702,0)(0.433012702,-0.25)(0.216506351,-0.375)(0,-0.25)}
\multips(0.433012702,1.125)(0.433012702,0){3}{\pspolygon(0,0)(0.216506351,0.125)(0.433012702,0)(0.433012702,-0.25)(0.216506351,-0.375)(0,-0.25)}
\multips(0.649519053,0.75)(0.433012702,0){2}{\pspolygon(0,0)(0.216506351,0.125)(0.433012702,0)(0.433012702,-0.25)(0.216506351,-0.375)(0,-0.25)}
\multips(0.866025404,0.375)(0.433012702,0){1}{\pspolygon(0,0)(0.216506351,0.125)(0.433012702,0)(0.433012702,-0.25)(0.216506351,-0.375)(0,-0.25)}
\rput(0.216506351,1.75){$\mathbf{0}$}
\rput(0.649519053,1.75){$\mathbf{1}$}
\rput(1.08253175,1.75){$\mathbf{2}$}
\rput(1.51554446,1.75){$\mathbf{3}$}
\rput(1.94855716,1.75){$\mathbf{4}$}
\rput(0.433012702,1.375){$\mathbf{8}$}
\rput(0.866025404,1.375){$\mathbf{9}$}
\rput(1.29903811,1.375){\scriptsize $\mathbf{10}$}
\rput(1.73205081,1.375){\scriptsize $\mathbf{11}$}
\rput(0.649519053,1){$\mathbf{1}$}
\rput(1.08253175,1){$\mathbf{2}$}
\rput(1.51554446,1){$\mathbf{3}$}
\rput(0.866025404,0.625){$\mathbf{9}$}
\rput(1.29903811,0.625){\scriptsize $\mathbf{10}$}
\rput(1.08253175,0.25){$\mathbf{2}$}
\end{pspicture}
\end{center}
\caption{\label{fig2}The doubly arithmetic triangle $\DAT(0,8,1,5)$ in $\Z/15\Z$.}
\end{figure}

\begin{rem}
Another necessary condition for a $\DAT$ of order $m$ to be balanced in $\Zn$ is that its cardinality, that is the binomial coefficient $\binom{m+1}{2}$, must be divisible by $n$. But these two necessary conditions are not sufficient: as depicted in Figure~\ref{fig2}, the triangle $\DAT(0,8,1,5)$ is not balanced in $\Z/15\Z$, although its cardinality $\binom{6}{2}=15$ is divisible by $n=15$ and its common differences $8$, $1$ and $7$ are invertible in $\Z/15\Z$.
\end{rem}

The following theorem is the main result of this section.

\begin{thm}\label{thm1}
Let $n\in\N$ be odd and let $d_1,d_2\in\Zn$ be invertible such that $d_1-d_2$ is also invertible. Then, the doubly arithmetic triangle $\DAT(a,d_1,d_2,m)$ is balanced in $\Zn$ for all $m\equiv0$ or $-1\pmod{n}$.
\end{thm}

\begin{proof}
Let $m$ be a multiple of $n$. We denote by $R_i$ the $i$th row of $\DAT(a,d_1,d_2,m)$, that is $R_i=(a+id_1+jd_2)_{0\leq j\leq m-1-i}$. We prove that, for $0\leq\lambda\leq m/n-2$, the consecutive $n$ rows $\{R_{\lambda n},R_{\lambda n+1},\ldots,R_{(\lambda+1)n-1}\}$ are balanced. Consider the permutation $\sigma$ of the set $\{0,1,\ldots,n-1\}$ defined by
$$
\sigma(i)\equiv i(d_1-d_2){d_1}^{-1}\pmod n
$$
for all $i\in\{0,1,\ldots,n-1\}$. Denote by $k_i$ the cardinality of the orbit of $i$ under $\sigma$. Let $\nabla(i,j)=a+id_1+jd_2$ denote the $j$th term in the $i$th row of $\DAT(a,d_1,d_2,m)$. Now, we show that, for every $i\in\{0,1,\ldots,n-1\}$, the concatenation $\cup_{l=0}^{k_i-1}R_{\lambda n+\sigma^{l}(i)}$ is balanced in $\Zn$. Since
$$
\begin{array}{l}
\displaystyle\nabla\left( \lambda n+\sigma^{l}(i) , m-1-\lambda n-\sigma^{l}(i) \right) + d_2 = a + (\lambda n+\sigma^{l}(i))d_1 + (m-\lambda n-\sigma^{l}(i))d_2\\[2ex]
\displaystyle = a + \sigma^{l}(i)(d_1-d_2) = a + \sigma^{l+1}(i)d_1 = \nabla\left(\lambda n+\sigma^{l+1}(i),0\right),
\end{array}
$$
for all $l\in\{0,1,\ldots,k_i-1\}$, it follows that the concatenation $\cup_{l=0}^{k_i-1}R_{\lambda n+\sigma^{l}(i)}$ is an arithmetic progression with invertible common difference $d_2$ and of length a multiple of $n$. Therefore, its multiplicity function is constant on $\Zn$. Finally, since $\left\{0,1,\ldots,n-1\right\}$ is a disjoint union of orbits under $\sigma$, the multiplicity function of $\cup_{i=0}^{n-1}R_{\lambda_n+i}$ is constant on $\Zn$ and thus the triangle $\DAT(a,d_1,d_2,m)$ is balanced in $\Zn$.
\par For $m\equiv-1\pmod{n}$, the doubly arithmetic triangle $\DAT(a,d_1,d_2,m)$ is obtained from the balanced triangle $\DAT(a,d_1,d_2,m+1)$ by rejecting its right side. Since it is an arithmetic progression with invertible common difference $d_1-d_2$ and of length $m+1\equiv0\pmod{n}$, it follows that this right side contains all the elements of $\Zn$ with the same multiplicity. This completes the proof.
\end{proof}

\begin{rem}
For $n$ odd and for every $d\in\Zn$ invertible, the doubly arithmetic triangles $\DAT(a,d,-d,m)$, $\DAT(a,d,2d,m)$ and $\DAT(a,2d,d,m)$ are balanced in $\Zn$, for all $m\equiv0$ or $-1\pmod{n}$.
\end{rem}

Let $n$ be a positive integer and let $d_1$ and $d_2$ be two elements of $\Zn$. The orbit $\mathcal{O}_S$, associated with a doubly infinite sequence $S$ in $\Zn$, is said to be \textit{$(d_1,d_2)$-doubly arithmetic} if each subtriangle appearing in it is a $\DAT$ with common differences $(d_1,d_2)$, that is if $\mathcal{O}_S$ is an orbit where all the diagonals are arithmetic progressions with the same common difference $d_1$ and where all the rows are arithmetic progressions with the same common difference $d_2$.

Now, we prove that, for every positive integer $n$, there does not exist a doubly arithmetic orbit in $\Zn$, except the trivial orbit generated by the sequence of zeros in $\Zn$.

\begin{prop}
Let $n$ be a positive integer. The orbit associated with the sequence of zeros is the only doubly arithmetic orbit in $\Zn$.
\end{prop}

\begin{proof}
It is clear that if $\mathcal{O}_S$ is $(d_1,d_2)$-doubly arithmetic, then $S$ is an arithmetic progression with common difference $d_2$. We set $S=\AP(a,d_2)$. It is known \cite{Chappelon2008}, and easy to retrieve, that the derived sequence $\partial S$ of $S$ is an arithmetic progression with common difference $2d_2$. Moreover, it is also $d_2$, by the doubly arithmetic structure of the orbit $\mathcal{O}_S$ and thus, the common difference $d_2$ vanishes. By the local rule in $\mathcal{O}_S$, we obtain that $a+d_1=2a$ and $a+2d_1=4a$. Therefore, we have $a=d_1=0$ and $S$ is the sequence of zeros. This completes the proof.
\end{proof}

Even if there does not exist a non-trivial doubly arithmetic orbit, the results of this section will be useful in next sections, where orbits with an interlaced doubly arithmetic structure are studied.

\section{Interlaced doubly arithmetic orbits of integers}

For all positive integers $n$, $k_1$ and $k_2$ and for every doubly infinite sequence $S$ in $\Zn$, or in $\Z$, the orbit $\mathcal{O}_S=\left(a_{i,j}\middle|a_{i+1,j}=a_{i,j}+a_{i,j+1},i\in\N,j\in\Z\right)$ is said to be \textit{$(k_1,k_2)$-interlaced doubly arithmetic} if, for every $i_0\in\{0,1,\ldots,k_1-1\}$ and every $j_0\in\{0,1,\ldots,k_2-1\}$, the subsequence $\left(a_{i_0+ik_1,j_0+jk_2}\middle|i\in\N,j\in\Z\right)$ is doubly arithmetic, i.e., if we have
$$
a_{i_0+ik_1,j_0+jk_2} = a_{i_0,j_0} + i(a_{i_0+k_1,j_0}-a_{i_0,j_0}) + j(a_{i_0,j_0+k_2}-a_{i_0,j_0}),
$$
for all $i\in\N$ and all $j\in\Z$.

Determining all interlaced doubly arithmetic orbits ($\IDAO$ for short) in $\Zn$ seems to be very difficult. Nevertheless, $\IDAO$ in $\Z$ are determined in this section and their projection in $\Zn$ will be considered in subsequent sections. First, it is clear that the sequence $S$ associated with a $(k_1,k_2)$-interlaced doubly arithmetic orbit $\mathcal{O}_S$ is a $k_2$-interlaced arithmetic progression. We begin by showing that the interlaced arithmetic structure of a sequence is preserved under the derivation process.

\begin{prop}\label{prop2}
Let $n$ be a positive integer. Let $(a_0,\ldots,a_{k-1})$ and $(d_0,\ldots,d_{k-1})$ be two $k$-tuples of elements in $\Zn$, or in $\Z$. Then, we have
$$
\begin{array}{l}
\displaystyle\partial \IAP\left( (a_0,\ldots,a_{k-1}) , (d_0,\ldots,d_{k-1}) \right)\\
\displaystyle \quad = \IAP\left( (a_0+a_1,\ldots,a_{k-2}+a_{k-1},a_{k-1}+a_0+d_0) , (d_0+d_1,\ldots,d_{k-2}+d_{k-1},d_{k-1}+d_{0}) \right).
\end{array}
$$
\end{prop}

\begin{proof}
Consider $S=\IAP\left((a_0,\ldots,a_{k-1}),(d_0,\ldots,d_{k-1})\right)={\left(x_j\right)}_{j\in\N}$ and $\partial S={\left(y_j\right)}_{j\in\N}$. Then, for all $l\in\Z$, we have
$$
y_{j_0+lk} = x_{j_0+lk}+x_{j_0+lk+1} = (a_{j_0}+ld_{j_0})+(a_{j_0+1}+ld_{j_0+1}) = (a_{j_0}+a_{j_0+1})+l(d_{j_0}+d_{j_0+1}),
$$
for all $j_0\in\{0,1,\ldots,k-2\}$, and
$$
y_{(k-1)+lk} = x_{(k-1)+lk}+x_{(l+1)k} = (a_{k-1}+ld_{k-1})+(a_0+(l+1)d_0) = (a_{k-1}+a_0+d_0)+l(d_{k-1}+d_0),
$$
for $j_0=k-1$. This completes the proof.
\end{proof}

We can now explicitly determine all the iterated derived sequences of an interlaced arithmetic progression.

\begin{prop}\label{prop5}
Let $n$ be a positive integer. Let $A$ and $D$ be two $k$-tuples of elements in $\Zn$, or in $\Z$. Then, for every integer $i\geq0$, we have
$$
\partial^{i}\IAP\left(A,D\right)
=
\IAP\left(
A\mathbf{C_i}+D\mathbf{T_i}\ ,\ D\mathbf{C_i}
\right),
$$
where $\mathbf{C_i}$ is the circulant matrix of size $k$ defined by
$$
\mathbf{C_i} = \circul\left( \sum_{l\geq0}\binom{i}{lk} , \sum_{l\geq0}\binom{i}{lk-1} , \ldots , \sum_{l\geq0}\binom{i}{lk+1} \right),
$$
and where $\mathbf{T_i}$ is the Toeplitz matrix of size $k$ where the $(r,s)$-entry of $\mathbf{T_i}$ is, for $1\leq r,s\leq k$, 
$$
\left(\mathbf{T_i}\right)_{r,s} = \sum_{l\geq0}l\binom{i}{r-s+lk}.
$$
\end{prop}

\begin{proof}
By iteration on $i$. Trivial for $i=0$. For $i=1$, Proposition~\ref{prop2} leads to
$$
\mathbf{C_1} = \circul(1,0,\ldots,0,1) = \left(\begin{array}{ccccc}
 1 & 0 & \cdots & 0 & 1 \\
 1 & 1 & & & 0 \\
 0 & \ddots & \ddots & & \vdots \\
 \vdots & & \ddots & \ddots & 0 \\
 0 & \cdots & 0 & 1 & 1
\end{array}\right)\ \text{and}\ \mathbf{T_1} = \left(\begin{array}{ccccc}
 0 & \cdots & \cdots & 0 & 1 \\
 \vdots & & & & 0 \\
 \vdots & & 0 & & \vdots \\
 \vdots & & & & \vdots \\
 0 & \cdots & \cdots & \cdots & 0 
\end{array}\right).
$$
We proceed by induction. Suppose that the result is true for some $i\geq1$. First, the $(i+1)$th derived sequence of $S=\IAP\left(A,D\right)$ is equal to
$$
\partial^{i+1}S
\begin{array}[t]{l}
\displaystyle = \partial\partial^{i}S =
\partial \IAP\left(
A\mathbf{C_i}+D\mathbf{T_i}\ ,\ D\mathbf{C_i}
\right)\\[2ex]
\displaystyle =
\IAP\left(
A\mathbf{C_i}\mathbf{C_1}+D\left(\mathbf{T_i}\mathbf{C_1}+\mathbf{C_i}\mathbf{T_1}\right)
\ ,\ 
D\mathbf{C_i}\mathbf{C_1}
\right).
\end{array}
$$
Since the product of two circulant matrices is also a circulant matrix, it follows that
$$
\mathbf{C_i}\mathbf{C_1}
\begin{array}[t]{l}
\displaystyle = \circul\left(\sum_{l\geq0}\binom{i}{lk},\sum_{l\geq0}\binom{i}{lk-1},\ldots,\sum_{l\geq0}\binom{i}{lk+1}\right)\circul\left(1,0,\ldots,0,1\right)\\[3ex]
\displaystyle = \circul\begin{array}[t]{l}
\displaystyle\left(\sum_{l\geq0}\binom{i}{lk}+\sum_{l\geq0}\binom{i}{lk-1},\sum_{l\geq0}\binom{i}{lk-1}+\sum_{l\geq0}\binom{i}{lk-2},\ldots\ldots\right. \\[3ex]
\displaystyle\left. \ldots\ldots,\sum_{l\geq0}\binom{i}{lk} +\sum_{l\geq0}\binom{i}{lk+1}\right)\end{array}\\[3ex]
\displaystyle = \circul\left(\sum_{l\geq0}\binom{i+1}{lk},\sum_{l\geq0}\binom{i+1}{lk-1},\ldots,\sum_{l\geq0}\binom{i+1}{lk+1}\right) = \mathbf{C_{i+1}}.
\end{array}
$$
Moreover, let $\mathbf{T_i}\mathbf{C_1}+\mathbf{C_i}\mathbf{T_1}=\left(\beta_{r,s}\right)$ for $1\leq r,s\leq k$. Note that
$$
\beta_{r,s} = \left(\mathbf{T_i}\mathbf{C_1}\right)_{r,s} + \left(\mathbf{C_i}\mathbf{T_1}\right)_{r,s} = \sum_{u=1}^{k}{\left(\mathbf{T_i}\right)}_{r,u}{\left(\mathbf{C_1}\right)}_{u,s} + \sum_{v=1}^{k}{\left(\mathbf{C_i}\right)}_{r,v}{\left(\mathbf{T_1}\right)}_{v,s}.
$$
Hence for $s<k$,
$$
\beta_{r,s} \begin{array}[t]{l}
\displaystyle = \left(\mathbf{T_i}\right)_{r,s} + \left(\mathbf{T_i}\right)_{r,s+1} = \sum_{l\geq0}l\binom{i}{r-s+lk}+\sum_{l\geq0}l\binom{i}{r-s+lk-1}\\[3ex]
\displaystyle = \sum_{l\geq0}l\binom{i+1}{r-s+lk} = \left(\mathbf{T_{i+1}}\right)_{r,s}.
\end{array}
$$
For $s=k$,
$$
\beta_{r,s} = \left(\mathbf{T_i}\right)_{r,1} + \left(\mathbf{T_i}\right)_{r,k} + \left(\mathbf{C_i}\right)_{r,1}
\begin{array}[t]{l}
\displaystyle = \sum_{l\geq0}l\binom{i}{r-1+lk} + \sum_{l\geq0}l\binom{i}{r-k+lk} + \sum_{l\geq0}\binom{i}{r-1+lk}\\[3ex]
\displaystyle = \sum_{l\geq0}(l+1)\binom{i}{r-1+lk} + \sum_{l\geq0}l\binom{i}{r-k+lk}\\[3ex]
\displaystyle = \sum_{l\geq0}l\binom{i}{r-1+(l-1)k} + \sum_{l\geq0}l\binom{i}{r+(l-1)k}\\[3ex]
\displaystyle = \sum_{l\geq0}l\binom{i+1}{r-k+lk} = \left(\mathbf{T_{i+1}}\right)_{r,k}.
\end{array}
$$
This completes the proof.
\end{proof}

The main result of this section is the complete characterization of $\IDAO$ in $\Z$.

\begin{thm}\label{thm2}
Every interlaced doubly arithmetic orbit $\mathcal{O}_S$ in $\Z$ is generated by an interlaced arithmetic progression of the form $S=\IAP((a_0,a_1,a_2),(d,-2d-3\Sigma,d+3\Sigma))$, where $a_0$, $a_1$, $a_2$ and $d$ are integers, and $\Sigma:=a_0+a_1+a_2$.
\end{thm}

We begin by showing that the interlaced arithmetic progressions listed in Theorem~\ref{thm2} will generate interlaced doubly arithmetic orbits of integers.

\begin{prop}\label{prop4}
Let $a_0,a_1,a_2,d\in\Z$ and let $\Sigma=a_0+a_1+a_2$. Then, the orbit $\mathcal{O}_S$ associated with $S=\IAP((a_0,a_1,a_2)(d,-2d-3\Sigma,d+3\Sigma))$ is $(6,3)$-interlaced doubly arithmetic.
\end{prop}

\begin{proof}
Let $\mathcal{O}_S=\left(a_{i,j}\middle|a_{i+1,j}=a_{i,j}+a_{i,j+1},i\in\N,j\in\Z\right)$ be the orbit associated with $S$ and let $S_{i_0,j_0}$ be the subsequence $S_{i_0,j_0}=\left(a_{i_0+6i,j_0+3j}\middle|i\in\N,j\in\Z\right)$, for all $i_0\in\{0,1,2,3,4,5\}$ and all $j_0\in\{0,1,2\}$. We can prove, by induction on $i$, that, for all $j\in\Z$, we have
$$
\begin{array}{lll}
a_{6i,3j} & = & a_0-2i(d+3\Sigma)+jd,\\
a_{6i,3j+1} & = & a_{1}-2id-j(2d+3\Sigma),\\
a_{6i,3j+2} & = & a_{2}+2i(2d+3\Sigma)+j(d+3\Sigma),\\
a_{6i+1,3j} & = & (a_0+a_1)-2i(2d+3\Sigma)-j(d+3\Sigma),\\
a_{6i+1,3j+1} & = & (a_1+a_2)+2i(d+3\Sigma)-jd,\\
a_{6i+1,3j+2} & = & (a_0+a_2+d)+2id+j(2d+3\Sigma),\\
a_{6i+2,3j} & = & (a_1+\Sigma)-2id-j(2d+3\Sigma),\\
a_{6i+2,3j+1} & = & (a_2+\Sigma+d)+2i(2d+3\Sigma)+j(d+3\Sigma),\\
a_{6i+2,3j+2} & = & (a_0-2\Sigma)-2i(d+3\Sigma)+jd,\\
a_{6i+3,3j} & = & (a_1+a_2+2\Sigma+d)+2i(d+3\Sigma)-jd,\\
a_{6i+3,3j+1} & = & (a_0+a_2-\Sigma+d)+2id+j(2d+3\Sigma),\\
a_{6i+3,3j+2} & = & (a_0+a_1-4\Sigma-2d)-2i(2d+3\Sigma)-j(d+3\Sigma),\\
a_{6i+4,3j} & = & (a_2+2\Sigma+2d)+2i(2d+3\Sigma)+j(d+3\Sigma),\\
a_{6i+4,3j+1} & = & (a_0-4\Sigma-d)-2i(d+3\Sigma)+jd,\\
a_{6i+4,3j+2} & = & (a_1-\Sigma-2d)-2id-j(2d+3\Sigma),\\
a_{6i+5,3j} & = & (a_0+a_2-2\Sigma+d)+2id+j(2d+3\Sigma),\\
a_{6i+5,3j+1} & = & (a_0+a_1-5\Sigma-3d)-2i(2d+3\Sigma)-j(d+3\Sigma),\\
a_{6i+5,3j+2} & = & (a_1+a_2+4\Sigma+d)+2i(d+3\Sigma)-jd.
\end{array}
$$
Thus, these $18$ subsequences $S_{i_0,j_0}$ are doubly arithmetic. This completes the proof.
\end{proof}

Now, we show that there is no other sequence generating $\IDAO$ in $\Z$. Since any $(k_1,k_2)$-$\IDAO$ is also a $(k_1k_2,k_1k_2)$-$\IDAO$, we suppose that we have $k_1=k_2=k$ in the sequel. The problem of determining all $(k,k)$-$\IDAO$ can then be converted into a system of linear equations.

\begin{prop}\label{prop3}
Let $n$ be a positive integer. Let $A$ and $D$ be two $k$-tuples of elements in $\Zn$, or in $\Z$, and let $S=\IAP(A,D)$ be a $k$-interlaced arithmetic progression. Then, the orbit $\mathcal{O}_S$ is $(k,k)$-interlaced doubly arithmetic if and only if $A$ and $D$ satisfy
$$
\left(\begin{array}{cc}
\mathbf{W_k}^2 & \mathbf{W_k}{\mathbf{T_k}}^{\!\mathrm{T}}\\
\mathbf{0_k} & \mathbf{W_k}
\end{array}\right)
\left(\begin{array}{c}
{A}^{\mathrm{T}}\\
{D}^{\mathrm{T}}
\end{array}\right)
=
0,
$$
where $\mathbf{W_k}=\mathbf{C_k}-\mathbf{I_k}=\circul\left(\binom{k}{0},\binom{k}{1},\ldots,\binom{k}{k-1}\right)$, that is the Wendt matrix of size $k$.
\end{prop}

The proof of this proposition is based on the following two lemmas.

\begin{lem}\label{lem3}
Let $n$ be a positive integer. Let $S$ be a $k$-interlaced arithmetic progression in $\Zn$, or in $\Z$. Then, the orbit $\mathcal{O}_S=\left(a_{i,j}\middle|a_{i+1,j}=a_{i,j}+a_{i,j+1},i\in\N,j\in\Z\right)$ is $(k,k)$-interlaced doubly arithmetic if and only if we have $(1)$: for every $i_0\in\{0,1,\ldots,k-1\}$ and for every $i\in\N$, the row $R_{ik+i_0}$ is of the same common differences as $R_{i_0}$, and $(2)$: for every $j_0\in\{0,1,\ldots,k-1\}$, the sequence $(a_{ik,j_0})_{i\in\N}$ is an arithmetic progression.
\end{lem}

\begin{proof}
If the orbit $\mathcal{O}_S$ is $(k,k)$-interlaced doubly arithmetic, then it is clear that the assertions $(1)$ and $(2)$ are verified. Suppose now that $(1)$ and $(2)$ hold. We begin by showing $(3)$: for every $j_0\in\{0,1,\ldots,k-1\}$ and for every $j\in\Z$, the sequence $(a_{ik,j_0+jk})_{i\in\N}$ is an arithmetic progression. Indeed, for every $i\in\N$, we have
$$
a_{ik,j_0+jk}
\begin{array}[t]{l}
\displaystyle\stackrel{(1)}{=} a_{ik,j_0} + j(a_{0,j_0+k}-a_{0,j_0})\\[2ex]
\displaystyle\stackrel{(2)}{=} a_{0,j_0} + i(a_{k,j_0}-a_{0,j_0}) + j(a_{0,j_0+k}-a_{0,j_0})\\[2ex]
\displaystyle = a_{0,j_0+jk} + i(a_{k,j_0}-a_{0,j_0}).
\end{array}
$$
Moreover, since $a_{i_0+ik,j_0}=\sum_{l=0}^{i_0}\binom{i_0}{l}a_{ik,j_0+l}$ by the local rule of the automaton, it follows that we have
$$
a_{i_0+ik,j_0+jk}
\begin{array}[t]{l}
\displaystyle\stackrel{(1)}{=} a_{i_0+ik,j_0} + j(a_{i_0,j_0+k}-a_{i_0,j_0}) = \sum_{l=0}^{i_0}\binom{i_0}{l}a_{ik,j_0+l} + j(a_{i_0,j_0+k}-a_{i_0,j_0})\\[2ex]
\displaystyle\stackrel{(3)}{=} \sum_{l=0}^{i_0}\binom{i_0}{l}(a_{0,j_0+l}+i(a_{k,j_0+l}-a_{0,j_0+l})) + j(a_{i_0,j_0+k}-a_{i_0,j_0})\\[2ex]
\displaystyle = \sum_{l=0}^{i_0}\binom{i_0}{l}a_{0,j_0+l}+i\left(\sum_{l=0}^{i_0}\binom{i_0}{l}a_{k,j_0+l}-\sum_{l=0}^{i_0}\binom{i_0}{l}a_{0,j_0+l}\right) + j(a_{i_0,j_0+k}-a_{i_0,j_0})\\[2ex]
\displaystyle = a_{i_0,j_0} + i(a_{i_0+k,j_0}-a_{i_0,j_0}) + j(a_{i_0,j_0+k}-a_{i_0,j_0}),
\end{array}
$$
for all $i_0$, $j_0\in\{0,1,\ldots,k-1\}$ and for all $i$, $j\in\Z$.
\end{proof}

\begin{lem}\label{lem2}
For all $i$, $j\in\N$ such that $0\leq j\leq i$, we have
\vspace{-0.5cm}
\begin{enumerate}
\item
$\mathbf{C_i}=\mathbf{C_1}^i$ and so $\mathbf{C_{ik}}=\mathbf{C_k}^{i}$,
\item
$\mathbf{T_i}=\mathbf{T_j}\mathbf{C_{i-j}}+\mathbf{C_j}\mathbf{T_{i-j}}$.
\end{enumerate}
\end{lem}

\begin{proof}
The first assertion follows from the recursive definition of $\mathbf{C_i}$. For the second assertion, we proceed by induction on $i$. The result is trivial for $i=0$ and for $i=1$. Suppose it is true until $i$ and prove it for $i+1$. It is clear for $j=0$ and for $j=i+1$. Let $j$ be an integer such that $1\leq j\leq i$. By the induction hypothesis and the recursive definition of $\mathbf{T_{i+1}}$ found in the proof of Proposition~\ref{prop5}, we obtain $\mathbf{T_{i+1}}=\mathbf{T_i}\mathbf{C_1}+\mathbf{C_i}\mathbf{T_1}=(\mathbf{T_j}\mathbf{C_{i-j}}+\mathbf{C_j}\mathbf{T_{i-j}})\mathbf{C_1}+\mathbf{C_i}\mathbf{T_1}=\mathbf{T_j}\mathbf{C_{i-j+1}}+\mathbf{C_j}(\mathbf{T_{i-j}}\mathbf{C_1}+\mathbf{C_{i-j}}\mathbf{T_1})=\mathbf{T_j}\mathbf{C_{i-j+1}}+\mathbf{C_j}\mathbf{T_{i-j+1}}$.
\end{proof}

We are now ready to prove Proposition~\ref{prop3}.

\begin{proof}[Proof of Proposition~\ref{prop3}]
Let $S=\IAP(A,D)$ be a $k$-interlaced arithmetic progression of elements in $\Zn$, or in $\Z$. We know that the orbit $\mathcal{O}_S$ is $(k,k)$-interlaced doubly arithmetic if and only if the assertions $(1)$ and $(2)$ are satisfied by Lemma~\ref{lem3}. We consider the equations $(1')$ and $(2')$:
$$
\begin{array}{l}
(1'):\ D\mathbf{W_k} = 0,\\[1.5ex]
(2'):\ A\mathbf{W_k}^2 + D\mathbf{T_k}\mathbf{W_k} = 0.
\end{array}
$$
First, by Proposition~\ref{prop5}, the assertions $(1)$ and $(1')$ are equivalent:
$$
(1)\ \begin{array}[t]{cl}
\displaystyle\stackrel{Prop.\ref{prop5}}{\Longleftrightarrow} & D\mathbf{C_{i_0+ik}}=D\mathbf{C_{i_0}},\ \text{for\ all}\ i\in\N\ \text{and}\ i_0\in\{0,1,\ldots,k-1\},\\[2ex]
\displaystyle\Longleftrightarrow & D(\mathbf{C_{ik}}-\mathbf{I_k})\mathbf{C_{i_0}} = 0,\ \text{for\ all}\ i\in\N\ \text{and}\ i_0\in\{0,1,\ldots,k-1\},\\[2ex]
\displaystyle\stackrel{Lem.\ref{lem2}}{\Longleftrightarrow} & D(\mathbf{C_k}-\mathbf{I_k})\displaystyle\sum_{l=0}^{i-1}\mathbf{C_k}^l\mathbf{C_{i_0}} = 0,\ \text{for\ all}\ i\in\N\ \text{and}\ i_0\in\{0,1,\ldots,k-1\},\\[3ex]
\displaystyle\Longleftrightarrow & D\mathbf{W_k}=0\quad(1').
\end{array}
$$
Proposition~\ref{prop5} also permits to put assertion $(2)$ in equation as follows:
$$
(2)\ \begin{array}[t]{cl}
\Longleftrightarrow & (a_{ik,0},a_{ik,1},\ldots,a_{ik,k-1})_{i\in\N}\ \text{is\ arithmetic},\\[2ex]
\stackrel{Prop.\ref{prop5}}{\Longleftrightarrow} & (A\mathbf{C_{ik}}+D\mathbf{T_{ik}})_{i\in\N}\ \text{is\ arithmetic},\\[2ex]
\Longleftrightarrow & A(\mathbf{C_{(i+2)k}}-2\mathbf{C_{(i+1)k}}+\mathbf{C_{ik}})+D(\mathbf{T_{(i+2)k}}-2\mathbf{T_{(i+1)k}}+\mathbf{T_{ik}})=0,\ \text{for\ all}\ i\in\N.
\end{array}
$$
Moreover, Lemma~\ref{lem2} leads to
$$
\mathbf{C_{(i+2)k}}-2\mathbf{C_{(i+1)k}}+\mathbf{C_{ik}} = ({\mathbf{C_k}}^2-2\mathbf{C_k}+\mathbf{I_k})\mathbf{C_{ik}} = \mathbf{W_k}^2\mathbf{C_{ik}},
$$
and
$$
\mathbf{T_{(i+2)k}}-2\mathbf{T_{(i+1)k}}+\mathbf{T_{ik}} = (\mathbf{T_{2k}}\mathbf{C_{ik}}+\mathbf{C_{2k}}\mathbf{T_{ik}})-2(\mathbf{T_k}\mathbf{C_{ik}}+\mathbf{C_k}\mathbf{T_{ik}})+\mathbf{T_{ik}}.
$$
Finally, since $D\mathbf{C_k}=D$ by assertion $(1')$, it follows that
$$
D\left(\mathbf{T_{(i+2)k}}-2\mathbf{T_{(i+1)k}}+\mathbf{T_{ik}}\right) = D(\mathbf{T_{2k}}-2\mathbf{T_k})\mathbf{C_{ik}} = D(\mathbf{T_{k}}\mathbf{C_{k}}+\mathbf{C_k}\mathbf{T_k}-2\mathbf{T_k})\mathbf{C_{ik}} = D\mathbf{T_k}\mathbf{W_k}\mathbf{C_{ik}}.
$$
Hence $A\mathbf{W_k}^2\mathbf{C_{ik}}+D\mathbf{T_k}\mathbf{W_k}\mathbf{C_{ik}}=0$, for all $i\in\N$ and so we have $(2')$. This completes the proof.
\end{proof}

In \cite{Wendt1894}, E.~Wendt investigated the resultant of $X^k-1$ and $(X+1)^k-1$, which corresponds to the determinant of $\mathbf{W_k}$. E.~Lehmer was the first to prove that the determinant of $\mathbf{W_k}$ vanishes if and only if $k$ is divisible by $6$ \cite{Lehmer1935}. It is also easy to deduce from her proof that the Wendt matrix $\mathbf{W_k}$ is of rank $k$ if $k$ is not divisible by $6$ and of rank $k-2$ otherwise.

\begin{prop}\label{prop14}
$$
\mathrm{rank}(\mathbf{W_k}) = \left\{\begin{array}{ll}
k & \text{if}\ k\not\equiv0\pmod{6},\\
k-2 & \text{if}\ k\equiv0\pmod{6}.
\end{array}\right.
$$
\end{prop}

We are now able to prove the main theorem of this section.

\begin{proof}[Proof of Theorem~\ref{thm2}]
If $k$ is not divisible by $6$, then the Wendt matrix $\mathbf{W_k}$ is of rank $k$ by Proposition~\ref{prop14}. This implies that $A=D=(0,\ldots,0)$ and thus $S$ is the sequence of zeros. Otherwise, if $k$ is divisible by $6$, then Proposition~\ref{prop4} implies that the vector space of $(k,k)$-interlaced doubly arithmetic orbits is of dimension greater than or equal to $4$. Moreover, since $\mathrm{rank}(\mathbf{W_k}^2)=\mathrm{rank}(\mathbf{W_k})=k-2$ by Proposition~\ref{prop14}, it follows that the matrix
$$
\left(\begin{array}{cc}
\mathbf{W_k}^2 & \mathbf{W_k}{\mathbf{T_k}}^{\!\mathrm{T}}\\
\mathbf{0_k} & \mathbf{W_k}
\end{array}\right)
$$
is of rank greater than or equal to $2k-4$. Therefore, there is no other $(k,k)$-$\IDAO$ than those listed in Theorem~\ref{thm2}. This completes the proof.
\end{proof}

\section{Balanced Steinhaus figures modulo an odd number}

In this section, we show that, for $n$ odd, the projection in $\Zn$ of an $\IDAO$ in $\Z$, obtained in the previous section, contains infinitely many balanced Steinhaus figures.

\begin{thm}\label{thm3}
Let $n\in\N$ be odd and let $a_0,a_1,a_2,d\in\Zn$. Define $\Sigma:=a_0+a_1+a_2$. If $d$, $d+3\Sigma$, and $2d+3\Sigma$ are invertible, then, the following Steinhaus figures, contained in the orbit of $S=\IAP((a_0,a_1,a_2),(d,-2d-3\Sigma,d+3\Sigma))$, are balanced:
\begin{itemize}
\item
every Steinhaus triangle of order $m$ in $\mathcal{O}_S$, for every $m\equiv0$ or $-1\pmod{6n}$,
\item
every Steinhaus trapezoid of order $m$ and of height $h$ in $\mathcal{O}_S$, for every $m\equiv0$ or $-1\pmod{6n}$ and for every $h\equiv m$ or $m+1\pmod{6n}$,
\item
every Pascal triangle of order $2m-1$ in $\mathcal{O}_S$, for every $m\equiv0$ or $-1\pmod{6n}$,
\item
every Pascal trapezoid of order $2m-1$ and of height $h$ in $\mathcal{O}_S$, for every $m\equiv0$ or $-1\pmod{6n}$ and for every $h\equiv m$ or $m+1\pmod{6n}$,
\item
every lozenge of order $2m-1$ in $\mathcal{O}_S$, for every $m\equiv0\pmod{6n}$.
\end{itemize}
\end{thm}

\begin{proof}
Let $\mathcal{O}_S=\left(a_{i,j}\middle|a_{i+1,j}=a_{i,j}+a_{i,j+1},i\in\N,j\in\Z\right)$ be the orbit associated with $S$. Consider the subsequences $S_{i_0,j_0}=\left(a_{i_0+6i,j_0+6j}\middle|i\in\N,j\in\Z\right)$, for $i_0$ and $j_0$ in $\{0,1,2,3,4,5\}$. Each of these $36$ subsequences is doubly arithmetic since the orbit $\mathcal{O}_S$ is $(6,3)$-interlaced doubly arithmetic by Proposition~\ref{prop4}. The following table gives their common differences $d_1$, $d_2$, $d_1-d_2$.
\begin{figure}[!h]
\begin{center}
\begin{tabular}{|c|c|c|c|}
\hline
$S_{i_0,j_0}$ & $d_1$ & $d_2$ & $d_1-d_2$\\
\hline\hline
$S_{1,2}$ , $S_{1,5}$ , $S_{3,1}$ , $S_{3,4}$ , $S_{5,0}$ , $S_{5,3}$ & $2d$ & $2(2d+3\Sigma)$ & $-2(d+3\Sigma)$\\
\hline
$S_{0,1}$ , $S_{0,4}$ , $S_{2,0}$ , $S_{2,3}$ , $S_{4,2}$ , $S_{4,5}$ & $-2d$ & $-2(2d+3\Sigma)$ & $2(d+3\Sigma)$\\
\hline
$S_{1,1}$ , $S_{1,4}$ , $S_{3,0}$ , $S_{3,3}$ , $S_{5,2}$ , $S_{5,5}$ & $2(d+3\Sigma)$ & $-2d$ & $2(2d+3\Sigma)$\\
\hline
$S_{0,0}$ , $S_{0,3}$ , $S_{2,2}$ , $S_{2,5}$ , $S_{4,1}$ , $S_{4,4}$ & $-2(d+3\Sigma)$ & $2d$ & $-2(2d+3\Sigma)$\\
\hline
$S_{0,2}$ , $S_{0,5}$ , $S_{2,1}$ , $S_{2,4}$ , $S_{4,0}$ , $S_{4,3}$ & $2(2d+3\Sigma)$ & $2(d+3\Sigma)$ & $2d$\\
\hline
$S_{1,0}$ , $S_{1,3}$ , $S_{3,2}$ , $S_{3,5}$ , $S_{5,1}$ , $S_{5,4}$ & $-2(2d+3\Sigma)$ & $-2(d+3\Sigma)$ & $-2d$\\
\hline
\end{tabular}
\end{center}
\end{figure}
Thus, each subsequence $S_{i_0,j_0}$ is doubly arithmetic, with invertible common differences $d_1$, $d_2$ and $d_1-d_2$. Let $\lambda\geq1$ and let $\nabla$ be a Steinhaus triangle of order $m=6\lambda n$ or $m=6\lambda n-1$, that appears in $\mathcal{O}_S$. Since $\nabla\cap S_{i_0,j_0}$, for $i_0$ and $j_0$ in $\{0,1,2,3,4,5\}$, is a doubly arithmetic triangle of order $\lambda n$ or $\lambda n-1$ and with invertible common differences $d_1$, $d_2$ and $d_1-d_2$, it follows from Theorem~\ref{thm1} that the $36$ subtriangles are balanced. Therefore their union, the Steinhaus triangle $\nabla$, is also balanced in $\Zn$. Similarly, every Pascal triangle of order $2m-1$ in $\mathcal{O}_S$ is balanced, for all $m\equiv0$ or $-1\pmod{6n}$, since it can be decomposed into $36$ subtriangles, which are balanced doubly arithmetic triangles by Theorem~\ref{thm1} again. For trapezoids, a Steinhaus trapezoid (resp. Pascal trapezoid) of order $m$ (resp. $2m-1$) and of height $h$ in $\mathcal{O}_S$ can be seen as the multiset difference between a Steinhaus triangle of order $m$ and a Steinhaus triangle of order $m-h$ (resp. between a Pascal triangle of order $2m-1$ and a Pascal triangle of order $2(m-h)-1$). Therefore, these trapezoids are balanced, for all $m\equiv0$ or $-1\pmod{6n}$ and for all $h\equiv m$ or $m+1\pmod{6n}$. Finally, a lozenge of order $2m-1$ in $\mathcal{O}_S$ is balanced, for all $m\equiv0\pmod{6n}$, since it is the multiset union of a Pascal triangle of order $2m-1$ and of a Steinhaus triangle of order $m-1$, which are both balanced in $\Zn$.
\end{proof}

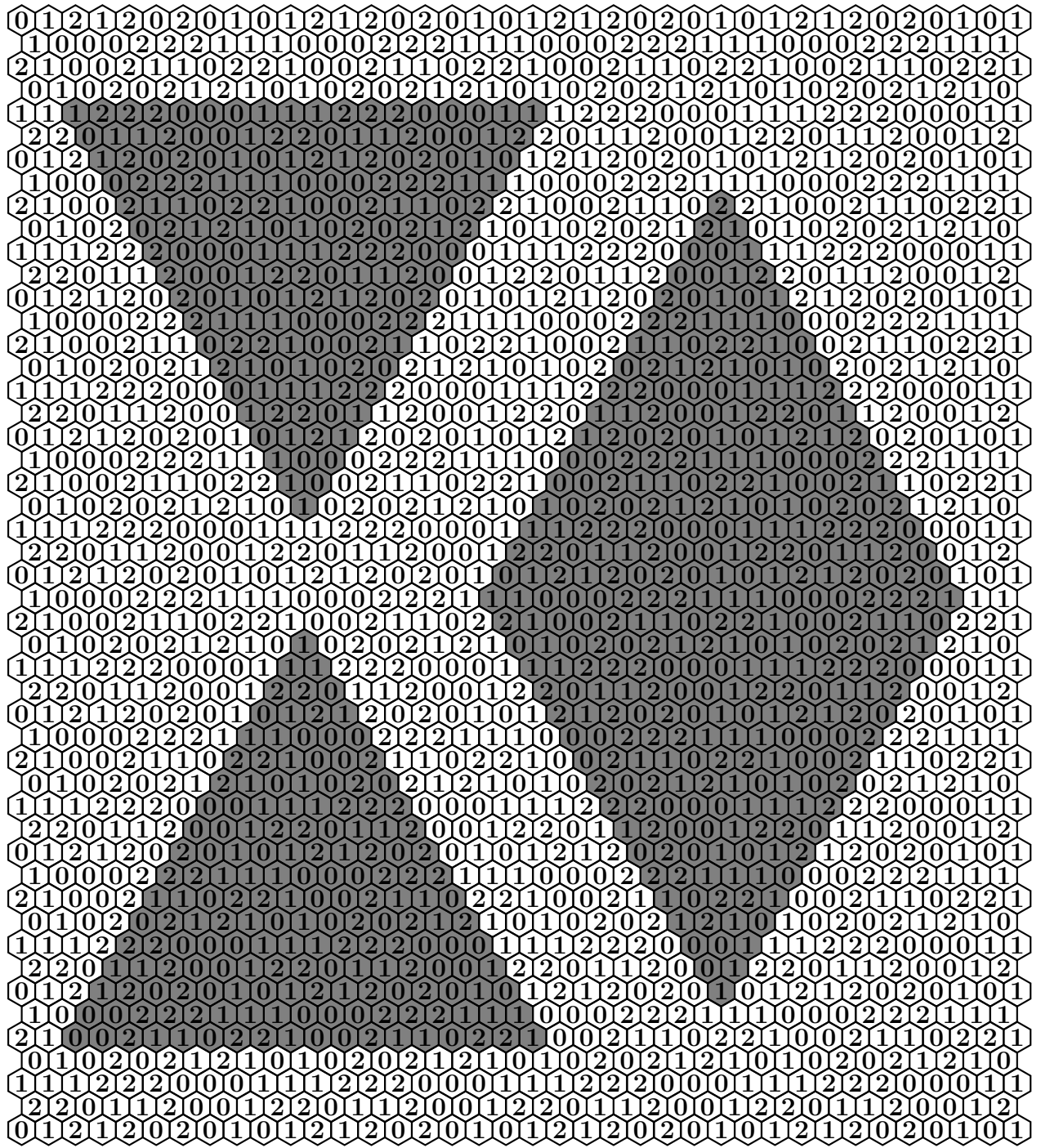
\begin{figure}[!p]
\begin{center}
\begin{pspicture}(16.4544827,18.5)

\multips(0.649519053,10.125)(0,0){1}{
\multips(3.89711432,0)(0.433012702,0){1}{\pspolygon*[linecolor=Gray](0,0.125)(0,0.375)(0.216506351,0.5)(0.433012702,0.375)(0.433012702,0.125)(0.216506351,0)}
\multips(3.68060797,0.375)(0.433012702,0){2}{\pspolygon*[linecolor=Gray](0,0.125)(0,0.375)(0.216506351,0.5)(0.433012702,0.375)(0.433012702,0.125)(0.216506351,0)}
\multips(3.46410162,0.75)(0.433012702,0){3}{\pspolygon*[linecolor=Gray](0,0.125)(0,0.375)(0.216506351,0.5)(0.433012702,0.375)(0.433012702,0.125)(0.216506351,0)}
\multips(3.24759526,1.125)(0.433012702,0){4}{\pspolygon*[linecolor=Gray](0,0.125)(0,0.375)(0.216506351,0.5)(0.433012702,0.375)(0.433012702,0.125)(0.216506351,0)}
\multips(3.03108891,1.5)(0.433012702,0){5}{\pspolygon*[linecolor=Gray](0,0.125)(0,0.375)(0.216506351,0.5)(0.433012702,0.375)(0.433012702,0.125)(0.216506351,0)}
\multips(2.81458256,1.875)(0.433012702,0){6}{\pspolygon*[linecolor=Gray](0,0.125)(0,0.375)(0.216506351,0.5)(0.433012702,0.375)(0.433012702,0.125)(0.216506351,0)}
\multips(2.59807621,2.25)(0.433012702,0){7}{\pspolygon*[linecolor=Gray](0,0.125)(0,0.375)(0.216506351,0.5)(0.433012702,0.375)(0.433012702,0.125)(0.216506351,0)}
\multips(2.38156986,2.625)(0.433012702,0){8}{\pspolygon*[linecolor=Gray](0,0.125)(0,0.375)(0.216506351,0.5)(0.433012702,0.375)(0.433012702,0.125)(0.216506351,0)}
\multips(2.16506351,3)(0.433012702,0){9}{\pspolygon*[linecolor=Gray](0,0.125)(0,0.375)(0.216506351,0.5)(0.433012702,0.375)(0.433012702,0.125)(0.216506351,0)}
\multips(1.94855716,3.375)(0.433012702,0){10}{\pspolygon*[linecolor=Gray](0,0.125)(0,0.375)(0.216506351,0.5)(0.433012702,0.375)(0.433012702,0.125)(0.216506351,0)}
\multips(1.73205081,3.75)(0.433012702,0){11}{\pspolygon*[linecolor=Gray](0,0.125)(0,0.375)(0.216506351,0.5)(0.433012702,0.375)(0.433012702,0.125)(0.216506351,0)}
\multips(1.51554446,4.125)(0.433012702,0){12}{\pspolygon*[linecolor=Gray](0,0.125)(0,0.375)(0.216506351,0.5)(0.433012702,0.375)(0.433012702,0.125)(0.216506351,0)}
\multips(1.29903811,4.5)(0.433012702,0){13}{\pspolygon*[linecolor=Gray](0,0.125)(0,0.375)(0.216506351,0.5)(0.433012702,0.375)(0.433012702,0.125)(0.216506351,0)}
\multips(1.08253175,4.875)(0.433012702,0){14}{\pspolygon*[linecolor=Gray](0,0.125)(0,0.375)(0.216506351,0.5)(0.433012702,0.375)(0.433012702,0.125)(0.216506351,0)}
\multips(0.866025404,5.25)(0.433012702,0){15}{\pspolygon*[linecolor=Gray](0,0.125)(0,0.375)(0.216506351,0.5)(0.433012702,0.375)(0.433012702,0.125)(0.216506351,0)}
\multips(0.649519053,5.625)(0.433012702,0){16}{\pspolygon*[linecolor=Gray](0,0.125)(0,0.375)(0.216506351,0.5)(0.433012702,0.375)(0.433012702,0.125)(0.216506351,0)}
\multips(0.433012702,6)(0.433012702,0){17}{\pspolygon*[linecolor=Gray](0,0.125)(0,0.375)(0.216506351,0.5)(0.433012702,0.375)(0.433012702,0.125)(0.216506351,0)}
\multips(0.216506351,6.375)(0.433012702,0){18}{\pspolygon*[linecolor=Gray](0,0.125)(0,0.375)(0.216506351,0.5)(0.433012702,0.375)(0.433012702,0.125)(0.216506351,0)}}

\multips(7.57772228,8.625)(0,0){1}{
\multips(0,0)(0.433012702,0){18}{\pspolygon*[linecolor=Gray](0,0.125)(0,0.375)(0.216506351,0.5)(0.433012702,0.375)(0.433012702,0.125)(0.216506351,0)}
\multips(0.216506351,0.375)(0.433012702,0){17}{\pspolygon*[linecolor=Gray](0,0.125)(0,0.375)(0.216506351,0.5)(0.433012702,0.375)(0.433012702,0.125)(0.216506351,0)}
\multips(0.433012702,0.75)(0.433012702,0){16}{\pspolygon*[linecolor=Gray](0,0.125)(0,0.375)(0.216506351,0.5)(0.433012702,0.375)(0.433012702,0.125)(0.216506351,0)}
\multips(0.649519053,1.125)(0.433012702,0){15}{\pspolygon*[linecolor=Gray](0,0.125)(0,0.375)(0.216506351,0.5)(0.433012702,0.375)(0.433012702,0.125)(0.216506351,0)}
\multips(0.866025404,1.5)(0.433012702,0){14}{\pspolygon*[linecolor=Gray](0,0.125)(0,0.375)(0.216506351,0.5)(0.433012702,0.375)(0.433012702,0.125)(0.216506351,0)}
\multips(1.08253175,1.875)(0.433012702,0){13}{\pspolygon*[linecolor=Gray](0,0.125)(0,0.375)(0.216506351,0.5)(0.433012702,0.375)(0.433012702,0.125)(0.216506351,0)}
\multips(1.29903811,2.25)(0.433012702,0){12}{\pspolygon*[linecolor=Gray](0,0.125)(0,0.375)(0.216506351,0.5)(0.433012702,0.375)(0.433012702,0.125)(0.216506351,0)}
\multips(1.51554446,2.625)(0.433012702,0){11}{\pspolygon*[linecolor=Gray](0,0.125)(0,0.375)(0.216506351,0.5)(0.433012702,0.375)(0.433012702,0.125)(0.216506351,0)}
\multips(1.73205081,3)(0.433012702,0){10}{\pspolygon*[linecolor=Gray](0,0.125)(0,0.375)(0.216506351,0.5)(0.433012702,0.375)(0.433012702,0.125)(0.216506351,0)}
\multips(1.94855716,3.375)(0.433012702,0){9}{\pspolygon*[linecolor=Gray](0,0.125)(0,0.375)(0.216506351,0.5)(0.433012702,0.375)(0.433012702,0.125)(0.216506351,0)}
\multips(2.16506351,3.75)(0.433012702,0){8}{\pspolygon*[linecolor=Gray](0,0.125)(0,0.375)(0.216506351,0.5)(0.433012702,0.375)(0.433012702,0.125)(0.216506351,0)}
\multips(2.38156986,4.125)(0.433012702,0){7}{\pspolygon*[linecolor=Gray](0,0.125)(0,0.375)(0.216506351,0.5)(0.433012702,0.375)(0.433012702,0.125)(0.216506351,0)}
\multips(2.59807621,4.5)(0.433012702,0){6}{\pspolygon*[linecolor=Gray](0,0.125)(0,0.375)(0.216506351,0.5)(0.433012702,0.375)(0.433012702,0.125)(0.216506351,0)}
\multips(2.81458256,4.875)(0.433012702,0){5}{\pspolygon*[linecolor=Gray](0,0.125)(0,0.375)(0.216506351,0.5)(0.433012702,0.375)(0.433012702,0.125)(0.216506351,0)}
\multips(3.03108891,5.25)(0.433012702,0){4}{\pspolygon*[linecolor=Gray](0,0.125)(0,0.375)(0.216506351,0.5)(0.433012702,0.375)(0.433012702,0.125)(0.216506351,0)}
\multips(3.24759526,5.625)(0.433012702,0){3}{\pspolygon*[linecolor=Gray](0,0.125)(0,0.375)(0.216506351,0.5)(0.433012702,0.375)(0.433012702,0.125)(0.216506351,0)}
\multips(3.46410162,6)(0.433012702,0){2}{\pspolygon*[linecolor=Gray](0,0.125)(0,0.375)(0.216506351,0.5)(0.433012702,0.375)(0.433012702,0.125)(0.216506351,0)}
\multips(3.68060797,6.375)(0.433012702,0){1}{\pspolygon*[linecolor=Gray](0,0.125)(0,0.375)(0.216506351,0.5)(0.433012702,0.375)(0.433012702,0.125)(0.216506351,0)}}

\multips(7.36121593,2.25)(0,0){1}{
\multips(3.89711432,0)(0.433012702,0){1}{\pspolygon*[linecolor=Gray](0,0.125)(0,0.375)(0.216506351,0.5)(0.433012702,0.375)(0.433012702,0.125)(0.216506351,0)}
\multips(3.68060797,0.375)(0.433012702,0){2}{\pspolygon*[linecolor=Gray](0,0.125)(0,0.375)(0.216506351,0.5)(0.433012702,0.375)(0.433012702,0.125)(0.216506351,0)}
\multips(3.46410162,0.75)(0.433012702,0){3}{\pspolygon*[linecolor=Gray](0,0.125)(0,0.375)(0.216506351,0.5)(0.433012702,0.375)(0.433012702,0.125)(0.216506351,0)}
\multips(3.24759526,1.125)(0.433012702,0){4}{\pspolygon*[linecolor=Gray](0,0.125)(0,0.375)(0.216506351,0.5)(0.433012702,0.375)(0.433012702,0.125)(0.216506351,0)}
\multips(3.03108891,1.5)(0.433012702,0){5}{\pspolygon*[linecolor=Gray](0,0.125)(0,0.375)(0.216506351,0.5)(0.433012702,0.375)(0.433012702,0.125)(0.216506351,0)}
\multips(2.81458256,1.875)(0.433012702,0){6}{\pspolygon*[linecolor=Gray](0,0.125)(0,0.375)(0.216506351,0.5)(0.433012702,0.375)(0.433012702,0.125)(0.216506351,0)}
\multips(2.59807621,2.25)(0.433012702,0){7}{\pspolygon*[linecolor=Gray](0,0.125)(0,0.375)(0.216506351,0.5)(0.433012702,0.375)(0.433012702,0.125)(0.216506351,0)}
\multips(2.38156986,2.625)(0.433012702,0){8}{\pspolygon*[linecolor=Gray](0,0.125)(0,0.375)(0.216506351,0.5)(0.433012702,0.375)(0.433012702,0.125)(0.216506351,0)}
\multips(2.16506351,3)(0.433012702,0){9}{\pspolygon*[linecolor=Gray](0,0.125)(0,0.375)(0.216506351,0.5)(0.433012702,0.375)(0.433012702,0.125)(0.216506351,0)}
\multips(1.94855716,3.375)(0.433012702,0){10}{\pspolygon*[linecolor=Gray](0,0.125)(0,0.375)(0.216506351,0.5)(0.433012702,0.375)(0.433012702,0.125)(0.216506351,0)}
\multips(1.73205081,3.75)(0.433012702,0){11}{\pspolygon*[linecolor=Gray](0,0.125)(0,0.375)(0.216506351,0.5)(0.433012702,0.375)(0.433012702,0.125)(0.216506351,0)}
\multips(1.51554446,4.125)(0.433012702,0){12}{\pspolygon*[linecolor=Gray](0,0.125)(0,0.375)(0.216506351,0.5)(0.433012702,0.375)(0.433012702,0.125)(0.216506351,0)}
\multips(1.29903811,4.5)(0.433012702,0){13}{\pspolygon*[linecolor=Gray](0,0.125)(0,0.375)(0.216506351,0.5)(0.433012702,0.375)(0.433012702,0.125)(0.216506351,0)}
\multips(1.08253175,4.875)(0.433012702,0){14}{\pspolygon*[linecolor=Gray](0,0.125)(0,0.375)(0.216506351,0.5)(0.433012702,0.375)(0.433012702,0.125)(0.216506351,0)}
\multips(0.866025404,5.25)(0.433012702,0){15}{\pspolygon*[linecolor=Gray](0,0.125)(0,0.375)(0.216506351,0.5)(0.433012702,0.375)(0.433012702,0.125)(0.216506351,0)}
\multips(0.649519053,5.625)(0.433012702,0){16}{\pspolygon*[linecolor=Gray](0,0.125)(0,0.375)(0.216506351,0.5)(0.433012702,0.375)(0.433012702,0.125)(0.216506351,0)}
\multips(0.433012702,6)(0.433012702,0){17}{\pspolygon*[linecolor=Gray](0,0.125)(0,0.375)(0.216506351,0.5)(0.433012702,0.375)(0.433012702,0.125)(0.216506351,0)}
\multips(0.216506351,6.375)(0.433012702,0){18}{\pspolygon*[linecolor=Gray](0,0.125)(0,0.375)(0.216506351,0.5)(0.433012702,0.375)(0.433012702,0.125)(0.216506351,0)}}

\multips(0.866025404,1.5)(0,0){1}{
\multips(0,0)(0.433012702,0){18}{\pspolygon*[linecolor=Gray](0,0.125)(0,0.375)(0.216506351,0.5)(0.433012702,0.375)(0.433012702,0.125)(0.216506351,0)}
\multips(0.216506351,0.375)(0.433012702,0){17}{\pspolygon*[linecolor=Gray](0,0.125)(0,0.375)(0.216506351,0.5)(0.433012702,0.375)(0.433012702,0.125)(0.216506351,0)}
\multips(0.433012702,0.75)(0.433012702,0){16}{\pspolygon*[linecolor=Gray](0,0.125)(0,0.375)(0.216506351,0.5)(0.433012702,0.375)(0.433012702,0.125)(0.216506351,0)}
\multips(0.649519053,1.125)(0.433012702,0){15}{\pspolygon*[linecolor=Gray](0,0.125)(0,0.375)(0.216506351,0.5)(0.433012702,0.375)(0.433012702,0.125)(0.216506351,0)}
\multips(0.866025404,1.5)(0.433012702,0){14}{\pspolygon*[linecolor=Gray](0,0.125)(0,0.375)(0.216506351,0.5)(0.433012702,0.375)(0.433012702,0.125)(0.216506351,0)}
\multips(1.08253175,1.875)(0.433012702,0){13}{\pspolygon*[linecolor=Gray](0,0.125)(0,0.375)(0.216506351,0.5)(0.433012702,0.375)(0.433012702,0.125)(0.216506351,0)}
\multips(1.29903811,2.25)(0.433012702,0){12}{\pspolygon*[linecolor=Gray](0,0.125)(0,0.375)(0.216506351,0.5)(0.433012702,0.375)(0.433012702,0.125)(0.216506351,0)}
\multips(1.51554446,2.625)(0.433012702,0){11}{\pspolygon*[linecolor=Gray](0,0.125)(0,0.375)(0.216506351,0.5)(0.433012702,0.375)(0.433012702,0.125)(0.216506351,0)}
\multips(1.73205081,3)(0.433012702,0){10}{\pspolygon*[linecolor=Gray](0,0.125)(0,0.375)(0.216506351,0.5)(0.433012702,0.375)(0.433012702,0.125)(0.216506351,0)}
\multips(1.94855716,3.375)(0.433012702,0){9}{\pspolygon*[linecolor=Gray](0,0.125)(0,0.375)(0.216506351,0.5)(0.433012702,0.375)(0.433012702,0.125)(0.216506351,0)}
\multips(2.16506351,3.75)(0.433012702,0){8}{\pspolygon*[linecolor=Gray](0,0.125)(0,0.375)(0.216506351,0.5)(0.433012702,0.375)(0.433012702,0.125)(0.216506351,0)}
\multips(2.38156986,4.125)(0.433012702,0){7}{\pspolygon*[linecolor=Gray](0,0.125)(0,0.375)(0.216506351,0.5)(0.433012702,0.375)(0.433012702,0.125)(0.216506351,0)}
\multips(2.59807621,4.5)(0.433012702,0){6}{\pspolygon*[linecolor=Gray](0,0.125)(0,0.375)(0.216506351,0.5)(0.433012702,0.375)(0.433012702,0.125)(0.216506351,0)}
\multips(2.81458256,4.875)(0.433012702,0){5}{\pspolygon*[linecolor=Gray](0,0.125)(0,0.375)(0.216506351,0.5)(0.433012702,0.375)(0.433012702,0.125)(0.216506351,0)}
\multips(3.03108891,5.25)(0.433012702,0){4}{\pspolygon*[linecolor=Gray](0,0.125)(0,0.375)(0.216506351,0.5)(0.433012702,0.375)(0.433012702,0.125)(0.216506351,0)}
\multips(3.24759526,5.625)(0.433012702,0){3}{\pspolygon*[linecolor=Gray](0,0.125)(0,0.375)(0.216506351,0.5)(0.433012702,0.375)(0.433012702,0.125)(0.216506351,0)}
\multips(3.46410162,6)(0.433012702,0){2}{\pspolygon*[linecolor=Gray](0,0.125)(0,0.375)(0.216506351,0.5)(0.433012702,0.375)(0.433012702,0.125)(0.216506351,0)}
\multips(3.68060797,6.375)(0.433012702,0){1}{\pspolygon*[linecolor=Gray](0,0.125)(0,0.375)(0.216506351,0.5)(0.433012702,0.375)(0.433012702,0.125)(0.216506351,0)}}

\multips(0,0)(0,0.75){25}{\multips(0,0)(0.433012702,0){38}{\pspolygon(0,0.125)(0,0.375)(0.216506351,0.5)(0.433012702,0.375)(0.433012702,0.125)(0.216506351,0)}}
\multips(0.216506351,0.375)(0,0.75){24}{\multips(0,0)(0.433012702,0){37}{\pspolygon(0,0.125)(0,0.375)(0.216506351,0.5)(0.433012702,0.375)(0.433012702,0.125)(0.216506351,0)}}

\multirput(0,0)(0,-2.25){9}{
\multirput(0.216506351,18.25)(3.89711432,0){5}{$\mathbf{0}$}
\multirput(0.649519053,18.25)(3.89711432,0){5}{$\mathbf{1}$}
\multirput(1.08253175,18.25)(3.89711432,0){4}{$\mathbf{2}$}
\multirput(1.51554446,18.25)(3.89711432,0){4}{$\mathbf{1}$}
\multirput(1.94855716,18.25)(3.89711432,0){4}{$\mathbf{2}$}
\multirput(2.38156986,18.25)(3.89711432,0){4}{$\mathbf{0}$}
\multirput(2.81458256,18.25)(3.89711432,0){4}{$\mathbf{2}$}
\multirput(3.24759526,18.25)(3.89711432,0){4}{$\mathbf{0}$}
\multirput(3.68060797,18.25)(3.89711432,0){4}{$\mathbf{1}$}}

\multirput(0,0)(0,-2.25){8}{
\multirput(0.433012702,17.875)(3.89711432,0){5}{$\mathbf{1}$}
\multirput(0.866025404,17.875)(3.89711432,0){4}{$\mathbf{0}$}
\multirput(1.29903811,17.875)(3.89711432,0){4}{$\mathbf{0}$}
\multirput(1.73205081,17.875)(3.89711432,0){4}{$\mathbf{0}$}
\multirput(2.16506351,17.875)(3.89711432,0){4}{$\mathbf{2}$}
\multirput(2.59807621,17.875)(3.89711432,0){4}{$\mathbf{2}$}
\multirput(3.03108891,17.875)(3.89711432,0){4}{$\mathbf{2}$}
\multirput(3.46410162,17.875)(3.89711432,0){4}{$\mathbf{1}$}
\multirput(3.89711432,17.875)(3.89711432,0){4}{$\mathbf{1}$}}

\multirput(0,0)(0,-2.25){8}{\multirput(0,-0.75)(0,0){1}{
\multirput(0.216506351,18.25)(3.89711432,0){5}{$\mathbf{2}$}
\multirput(0.649519053,18.25)(3.89711432,0){5}{$\mathbf{1}$}
\multirput(1.08253175,18.25)(3.89711432,0){4}{$\mathbf{0}$}
\multirput(1.51554446,18.25)(3.89711432,0){4}{$\mathbf{0}$}
\multirput(1.94855716,18.25)(3.89711432,0){4}{$\mathbf{2}$}
\multirput(2.38156986,18.25)(3.89711432,0){4}{$\mathbf{1}$}
\multirput(2.81458256,18.25)(3.89711432,0){4}{$\mathbf{1}$}
\multirput(3.24759526,18.25)(3.89711432,0){4}{$\mathbf{0}$}
\multirput(3.68060797,18.25)(3.89711432,0){4}{$\mathbf{2}$}}}

\multirput(0,0)(0,-2.25){8}{\multirput(0,-0.75)(0,0){1}{
\multirput(0.433012702,17.875)(3.89711432,0){5}{$\mathbf{0}$}
\multirput(0.866025404,17.875)(3.89711432,0){4}{$\mathbf{1}$}
\multirput(1.29903811,17.875)(3.89711432,0){4}{$\mathbf{0}$}
\multirput(1.73205081,17.875)(3.89711432,0){4}{$\mathbf{2}$}
\multirput(2.16506351,17.875)(3.89711432,0){4}{$\mathbf{0}$}
\multirput(2.59807621,17.875)(3.89711432,0){4}{$\mathbf{2}$}
\multirput(3.03108891,17.875)(3.89711432,0){4}{$\mathbf{1}$}
\multirput(3.46410162,17.875)(3.89711432,0){4}{$\mathbf{2}$}
\multirput(3.89711432,17.875)(3.89711432,0){4}{$\mathbf{1}$}}}

\multirput(0,0)(0,-2.25){8}{\multirput(0,-1.5)(0,0){1}{
\multirput(0.216506351,18.25)(3.89711432,0){5}{$\mathbf{1}$}
\multirput(0.649519053,18.25)(3.89711432,0){5}{$\mathbf{1}$}
\multirput(1.08253175,18.25)(3.89711432,0){4}{$\mathbf{1}$}
\multirput(1.51554446,18.25)(3.89711432,0){4}{$\mathbf{2}$}
\multirput(1.94855716,18.25)(3.89711432,0){4}{$\mathbf{2}$}
\multirput(2.38156986,18.25)(3.89711432,0){4}{$\mathbf{2}$}
\multirput(2.81458256,18.25)(3.89711432,0){4}{$\mathbf{0}$}
\multirput(3.24759526,18.25)(3.89711432,0){4}{$\mathbf{0}$}
\multirput(3.68060797,18.25)(3.89711432,0){4}{$\mathbf{0}$}}}

\multirput(0,0)(0,-2.25){8}{\multirput(0,-1.5)(0,0){1}{
\multirput(0.433012702,17.875)(3.89711432,0){5}{$\mathbf{2}$}
\multirput(0.866025404,17.875)(3.89711432,0){4}{$\mathbf{2}$}
\multirput(1.29903811,17.875)(3.89711432,0){4}{$\mathbf{0}$}
\multirput(1.73205081,17.875)(3.89711432,0){4}{$\mathbf{1}$}
\multirput(2.16506351,17.875)(3.89711432,0){4}{$\mathbf{1}$}
\multirput(2.59807621,17.875)(3.89711432,0){4}{$\mathbf{2}$}
\multirput(3.03108891,17.875)(3.89711432,0){4}{$\mathbf{0}$}
\multirput(3.46410162,17.875)(3.89711432,0){4}{$\mathbf{0}$}
\multirput(3.89711432,17.875)(3.89711432,0){4}{$\mathbf{1}$}}}

\end{pspicture}
\end{center}
\caption{\label{fig14}Balanced Steinhaus figures in the orbit of $\IAP((0,1,2),(1,1,1))$ in $\Z/3\Z$.}
\end{figure}

The case where $a_0=0$, $a_1=1$, $a_2=2$ and $d=1$ in $\Z/3\Z$, i.e., the orbit associated with the sequence $\IAP((0,1,2),(1,1,1))$, is illustrated in Figure~\ref{fig14}. In this example, balanced Steinhaus figures are depicted in gray: there are a balanced Steinhaus triangle of order $18$, a balanced Pascal triangle of order $35$ and a balanced lozenge of order $35$.

\section{The antisymmetric case}

In this section, we refine Theorem~\ref{thm3} by considering antisymmetric sequences in $\Zn$.

A finite sequence $S=(a_0,\ldots,a_{m-1})$ of length $m\geq1$ in $\Zn$, or in $\Z$, is said to be \textit{antisymmetric} if $a_{m-1-j}=-a_j$ for all $j\in\{0,1,\ldots,m-1\}$.

For examples, the sequences $(1,4,0,3,6)$ and $(2,6,1,5)$ are antisymmetric in $\Z/7\Z$. It is known, see \cite{Chappelon2008}, that the antisymmetry of finite sequences is preserved by the derivation process.

\begin{prop}\label{prop6}
Let $n$ be a positive integer and let $S=(a_0,\ldots,a_{m-1})$ be a finite sequence in $\Zn$, or in $\Z$. Then, the sequence $S$ is antisymmetric if and only if its derived sequence $\partial S$ is also antisymmetric and $a_{\left\lfloor m/2\right\rfloor}+a_{m-\left\lfloor m/2\right\rfloor}=0$, where $\left\lfloor m/2\right\rfloor$ is the floor of $m/2$.
\end{prop}

\begin{proof}
We set $\partial S=(b_0,\ldots,b_{m-2})=(a_0+a_1,\ldots,a_{m-2}+a_{m-1})$. If $S$ is antisymmetric, then $\partial S$ is also antisymmetric since, for all $j\in\{0,1,\ldots,m-2\}$, we have $b_{m-2-j} = a_{m-2-j}+a_{m-1-j} = -a_{j+1}-a_{j} = -b_j$. Conversely, if $\partial S$ is antisymmetric and $a_{\left\lfloor m/2\right\rfloor}+a_{m-\left\lfloor m/2\right\rfloor}=0$, we proceed by decreasing induction on $j$. Since
$$
\sum_{k=j}^{m-2-j}b_k = \sum_{k=j}^{m-2-j}(a_{k}+a_{k+1}) = a_{j} + 2\sum_{k=j+1}^{m-2-j}a_k + a_{m-1-j},
$$
it follows that
$$
a_j+a_{m-1-j}=\sum_{k=j}^{m-2-j}b_{k}-2\sum_{k=j+1}^{m-2-j}a_k=0,
$$
by the decreasing induction hypothesis. This completes the proof.
\end{proof}

The main interest of the antisymmetric sequences in $\Zn$ is that their multiplicity function admits a certain symmetry. Indeed, it is clear that, if $S$ is an antisymmetric sequence in $\Zn$, then its multiplicity function $\m_S$ satisfies $\m_S(x)=\m_S(-x)$, for all $x$ in $\Zn$. The same equality appears for the multiplicity function of Steinhaus or Pascal triangles generated by antisymmetric sequences.

\begin{prop}\label{prop7}
Let $n$ be a positive integer and let $S$ be an antisymmetric sequence of length $m\geq1$ in $\Zn$. Then, we have $\m_{\nabla S}(x)=\m_{\nabla S}(-x)$ for all $x\in\Zn$.
\end{prop}

\begin{proof}
Since each derived sequence $\partial^{i}S$, for $i\in\{0,1,\ldots,m-1\}$, is antisymmetric by Proposition~\ref{prop6}, it follows that $m_{\nabla S}(x)=\sum_{i=0}^{m-1}\m_{\partial^{i}S}(x)=\sum_{i=0}^{m-1}\m_{\partial^{i}S}(-x)=m_{\nabla S}(-x)$ for all $x\in\Zn$.
\end{proof}

\begin{prop}\label{prop11}
Let $n$ be a positive integer and let $S$ be an antisymmetric sequence of length $2m-1\geq1$ in $\Zn$. Then, we have $\m_{\Delta S}(x)=\m_{\Delta S}(-x)$ for all $x\in\Zn$.
\end{prop}

Now, for $n$ odd, we determine all the sequences generating $\IDAO$ in $\Z$ and such that the first $3n$ terms of their projection in $\Zn$ are antisymmetric.

For every doubly infinite sequence $S=(a_j)_{j\in\Z}$ in $\Zn$, or in $\Z$, and for all integers $j_0$ and $j_1$ such that $j_0\leq j_1$, we let $S[j_0,j_1]$ denote the subsequence of $S$ indexed between $j_0$ and $j_1$, that is $S[j_0,j_1]=(a_{j_0},a_{j_0+1},\ldots,a_{j_1})$.

\begin{prop}\label{prop8}
Let $n\in\N$ be odd. Let $a_0,a_1,a_2,d\in\Zn$ and let $\Sigma=a_0+a_1+a_2$. Then, the subsequence $S_m=\IAP((a_0,a_1,a_2),(d,-2d-3\Sigma,d+3\Sigma))[0,m-1]$, of length $m\equiv0\pmod{3n}$ in $\Zn$, is antisymmetric if and only if $\Sigma=0$ and $a_1=-d$, i.e., if we have $S_m=\IAP((a,-d,d-a),(d,-2d,d))[0,m-1]$.
\end{prop}

\begin{proof}
Set $m=3\lambda n$ and $S_m=\IAP((a_0,a_1,a_2),(d,-2d-3\Sigma,d+3\Sigma))[0,m-1]=(a_0,\ldots,a_{m-1})$ in $\Zn$. If $S_m$ is antisymmetric, then its terms $a_j$ must satisfy
$$
\left\{\begin{array}{l}
a_{3j}+a_{3(\lambda n-j-1)+2} = 0 \\
a_{3j+1}+a_{3(\lambda n-j-1)+1} = 0 \\
a_{3j+2}+a_{3(\lambda n-j-1)} = 0
\end{array}\right.\quad\Longleftrightarrow\quad
\left\{\begin{array}{l}
a_0+a_2-d-3(j+1)\Sigma = 0 \\
2a_1+2d+3\Sigma = 0\\
a_0+a_2-d+3j\Sigma = 0
\end{array}\right.,\ \text{for\ all}\ 0\leq j\leq n-1.
$$
This leads to $a_1=-d$, $a_2=d-a_0$ and $\Sigma=0$, since $n$ is odd, and thus $S_m=\IAP((a_0,-d,d-a_0),(d,-2d,d))[0,m-1]$, as announced.
\end{proof}

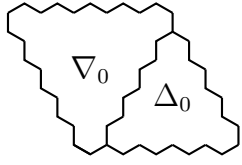
\begin{figure}[!h]
\begin{center}
\begin{pspicture}(3.24759526,2)

\multips(0,1.5)(0.433012702,0){5}{\psline(0,0.375)(0.216506351,0.5)(0.433012702,0.375)}
\multips(0,1.5)(0.216506351,-0.375){5}{\psline(0,0.375)(0,0.125)(0.216506351,0)}
\multips(1.73205081,1.5)(0.216506351,-0.375){4}{\psline(0.433012702,0.375)(0.433012702,0.125)(0.649519053,0)}
\multips(0.866025404,0)(0.433012702,0){5}{\psline(0,0.125)(0.216506351,0)(0.433012702,0.125)}
\psline(3.03108891,0.125)(3.03108891,0.375)
\multips(0.866025404,0)(0.216506351,0.375){4}{\psline(0.433012702,0.125)(0.433012702,0.375)(0.649519053,0.5)}

\rput(1.08253175,1.125){$\nabla_0$}
\rput(2.16506351,0.75){$\Delta_0$}
\end{pspicture}
\caption{\label{fig10}$\nabla_0 = \nabla S[0,3n-1]$ and $\Delta_0=\Delta\partial S[1,6n-3]$.}
\end{center}
\end{figure}

Let $n$ be an odd number and let $a$ and $d$ be two elements in $\Zn$ with $d$ invertible. We refine Theorem~\ref{thm3} by considering the orbit $\mathcal{O}_S$ of the sequence $S=IAP((a,-d,d-a),(d,-2d,d))$. Let $\nabla_0$ be the Steinhaus triangle, of order $3n$, generated by the first $3n$ terms of $S$ and let $\Delta_0$ be the Pascal triangle, of order $6n-3$, adjacent with $\nabla_0$ as depicted in Figure~\ref{fig10}, that are $\nabla_0 = \nabla S[0,3n-1]$ and $\Delta_0=\Delta\partial S[1,6n-3]$.

We begin by showing that these triangles are balanced in $\Zn$.

\begin{prop}\label{prop12}
Let $n\in\N$ be odd and let $a,d\in\Zn$ with $d$ invertible. Consider the $3$-interlaced arithmetic progression $S=\IAP((a,-d,d-a)(d,-2d,d))$. Then, the triangles $\nabla_0=\nabla S[0,3n-1]$ and $\Delta_0=\Delta\partial S[1,6n-3]$ are balanced in $\Zn$.
\end{prop}

\begin{proof}
First, since the derived sequences of $S=\IAP((a,-d,d-a),(d,-2d,d))$ are:
$$
\begin{array}{rcl}
\partial^{3i}S & = & (-1)^i\IAP((a-id,-(i+1)d,(2i+1)d-a),(d,-2d,d)),\\
\partial^{3i+1}S & = & (-1)^i\IAP((a-(2i+1)d,id-a,(i+2)d),(-d,-d,2d)),\\
\partial^{3i+2}S & = & (-1)^i\IAP((-(i+1)d,(2i+2)d-a,a-id),(-2d,d,d)),
\end{array}
$$
for all $i\in\N$, it follows that $\partial^{3n}S=-S$. Moreover, the orbit $\mathcal{O}_S$ is $(6,3)$-interlaced doubly arithmetic and thus each row (resp. each diagonal) of $\mathcal{O}_S$ is periodic of period $3n$ (resp. of period $6n$). This leads to the following periodic decomposition of the orbit $\mathcal{O}_S$ into triangles $\nabla_0$ and $\Delta_0$:
\begin{center}
\begin{pspicture}(9.74278579,2.875)

\multips(0,1.875)(0,-1.875){2}{\multips(0,0)(1.08253175,0){9}{
\multips(0,0.75)(0.216506351,0){5}{\psline(0,0.1875)(0.108253175,0.25)(0.216506351,0.1875)}
\multips(0,0.75)(0.108253175,-0.1875){5}{\psline(0,0.1875)(0,0.0625)(0.108253175,0)}
\multips(0.866025404,0.75)(-0.108253175,-0.1875){5}{\psline(0.216506351,0.1875)(0.216506351,0.0625)(0.108253175,0)}}}

\multips(0.541265877,0.9375)(1.08253175,0){8}{
\multips(0,0.75)(0.216506351,0){5}{\psline(0,0.1875)(0.108253175,0.25)(0.216506351,0.1875)}
\multips(0,0.75)(0.108253175,-0.1875){5}{\psline(0,0.1875)(0,0.0625)(0.108253175,0)}
\multips(0.866025404,0.75)(-0.108253175,-0.1875){5}{\psline(0.216506351,0.1875)(0.216506351,0.0625)(0.108253175,0)}}

\multips(0.324759526,1.6875)(-0.108253175,-0.1875){5}{\psline(0.216506351,0.1875)(0.216506351,0.0625)(0.108253175,0)}
\multips(9.20151992,1.6875)(0.108253175,-0.1875){5}{\psline(0,0.1875)(0,0.0625)(0.108253175,0)}
\multips(0.541265877,0)(0.216506351,0){40}{\psline(0,0)(0.108253175,0.0625)(0.216506351,0)}

\multirput(0.541265877,0.625)(1.08253175,0){9}{\scriptsize$\nabla_0$}
\multirput(1.08253175,1.5625)(1.08253175,0){8}{\scriptsize-$\nabla_0$}
\multirput(0.541265877,2.5)(1.08253175,0){9}{\scriptsize$\nabla_0$}

\multirput(1.08253175,0.375)(1.08253175,0){8}{\scriptsize$\Delta_0$}
\multirput(0.541265877,1.3125)(1.08253175,0){9}{\scriptsize-$\Delta_0$}
\multirput(1.08253175,2.25)(1.08253175,0){8}{\scriptsize$\Delta_0$}

\end{pspicture}
\end{center}
Particularly, the Steinhaus triangle $\nabla S[0,6n-1]$, of order $6n$, and the Pascal triangle $\Delta\partial S[1,12n-3]$, of order $12n-3$, which are balanced by Theorem~\ref{thm3}, admit the decomposition:
\begin{center}
$\nabla S[0,6n-1]=$\!\!\!\begin{tabular}{c}\begin{pspicture}(2.16506351,1.9375)
\multips(0,0.9375)(1.08253175,0){2}{
\multips(0,0.75)(0.216506351,0){5}{\psline(0,0.1875)(0.108253175,0.25)(0.216506351,0.1875)}
\multips(0,0.75)(0.108253175,-0.1875){5}{\psline(0,0.1875)(0,0.0625)(0.108253175,0)}
\multips(0.866025404,0.75)(-0.108253175,-0.1875){5}{\psline(0.216506351,0.1875)(0.216506351,0.0625)(0.108253175,0)}}
\multips(0.541265877,0)(1.08253175,0){1}{
\multips(0,0.75)(0.216506351,0){5}{\psline(0,0.1875)(0.108253175,0.25)(0.216506351,0.1875)}
\multips(0,0.75)(0.108253175,-0.1875){5}{\psline(0,0.1875)(0,0.0625)(0.108253175,0)}
\multips(0.866025404,0.75)(-0.108253175,-0.1875){5}{\psline(0.216506351,0.1875)(0.216506351,0.0625)(0.108253175,0)}}
\rput(0.541265877,1.5625){\scriptsize$\nabla_0$}
\rput(1.62379763,1.5625){\scriptsize$\nabla_0$}
\rput(1.08253175,0.625){\scriptsize-$\nabla_0$}
\rput(1.08253175,1.3125){\scriptsize$\Delta_0$}
\end{pspicture}\end{tabular}\quad and\quad $\Delta\partial S[1,12n-3]=$\!\!\!\begin{tabular}{c}\begin{pspicture}(1.94855716,1.75)
\multips(0.974278579,1.75)(-0.108253175,-0.1875){9}{\psline(0,0)(-0.108253175,-0.0625)(-0.108253175,-0.1875)}
\multips(0.974278579,1.75)(0.108253175,-0.1875){9}{\psline(0,0)(0.108253175,-0.0625)(0.108253175,-0.1875)}
\multips(0,0)(0.216506351,0){9}{\psline(0,0.0625)(0.108253175,0)(0.216506351,0.0625)}
\multips(0.541265877,0.9375)(0.216506351,0){4}{\psline(0,0.0625)(0.108253175,0)(0.216506351,0.0625)}
\multips(0.433012702,0.8125)(0.108253175,-0.1875){4}{\psline(0,0)(0.108253175,-0.0625)(0.108253175,-0.1875)}
\multips(1.51554446,0.8125)(-0.108253175,-0.1875){4}{\psline(0,0)(-0.108253175,-0.0625)(-0.108253175,-0.1875)}
\rput(0.974278579,1.3125){\scriptsize$\Delta_0$}
\rput(0.974278579,0.625){\scriptsize-$\nabla_0$}
\rput(0.433012702,0.3125){\scriptsize-$\Delta_0$}
\rput(1.51554446,0.3125){\scriptsize-$\Delta_0$}
\end{pspicture}\end{tabular}.
\end{center}
The sequences $S[0,3n-1]$ and $\partial S[1,6n-3]$ are antisymmetric in $\Zn$, by Propositions~\ref{prop8} and \ref{prop6}, and thus we deduce, from Propositions~\ref{prop7} and \ref{prop11}, that the multiplicity functions $\m_{-\nabla_0}$ and $\m_{-\Delta_0}$ correspond to $\m_{\nabla_0}$ and $\m_{\Delta_0}$, since $\m_{\nabla_0}(x)=\m_{\nabla_0}(-x)=\m_{-\nabla_0}(x)$ and $\m_{\Delta_0}(x)=\m_{\Delta_0}(-x)=\m_{-\Delta_0}(x)$, for all $x\in\Zn$. Finally, the multiplicity functions $\m_{\nabla_0}$ and $\m_{\Delta_0}$ are constant because they are solutions of the following system of equations
$$
\begin{array}{l}
\displaystyle 3\m_{\nabla_0}+\m_{\Delta_0}=\m_{\nabla S[0,3n-1]}=\frac{1}{n}\binom{3n+1}{2},\\[2ex]
\displaystyle \m_{\nabla_0}+3\m_{\Delta_0}=\m_{\Delta\partial S[1,6n-3]}=\frac{1}{n}\binom{3n}{2}.
\end{array}
$$
Therefore, the elementary triangles $\nabla_0$ and $\Delta_0$ are balanced in $\Zn$.
\end{proof}

Finally, we obtain the refinement of Theorem~\ref{thm3} announced above.

\begin{thm}\label{thm4}
Let $n\in\N$ be odd and let $a,d\in\Zn$ with $d$ invertible. Then, the following Steinhaus figures, contained in the orbit of $S=\IAP((a,-d,d-a),(d,-2d,d))$, are balanced:
\begin{itemize}
\item
the Steinhaus triangles $\nabla S[0,3\lambda n-1]$ of order $3\lambda n$, and $\nabla\partial S[0,3\lambda n-2]$ of order $3\lambda n-1$, for every integer $\lambda\geq1$,
\item
the Steinhaus trapezoid $\ST(S[0,3\lambda n-1],h)$ of order $3\lambda n$ and of height $h$, for every integer $\lambda\geq1$ and for every $h\equiv0$ or $1\pmod{3n}$; the Steinhaus trapezoid $\ST(\partial S[0,3\lambda n-2],h)$ of order $3\lambda n-1$ and of height $h$, for every integer $\lambda\geq1$ and for every $h\equiv -1$ or $0\pmod{3n}$,
\item
the Pascal triangle $\Delta\partial S[-m,m-2]$ of order $2m-1$, for every $m\equiv0$ or $-1\pmod{3n}$,
\item
the Pascal trapezoid $\PT(\partial S[-m,m-2],h)$ of order $2m-1$ and of height $h$, for every $m\equiv0$ or $-1\pmod{3n}$ and for every $h\equiv m$ or $m+1\pmod{3n}$,
\item
the lozenge $\lozenge\partial S[-m,m-2]$ of order $2m-1$, for every $m\equiv0\pmod{3n}$.
\end{itemize}
\end{thm}

\begin{proof}
For every integer $\lambda\geq1$, the Steinhaus triangle $\nabla S[0,3\lambda n-1]$  and the Pascal triangle $\Delta\partial S[-3\lambda n,3\lambda n-2]$ are balanced because they are multiset unions of the elementary triangles $\nabla_0$, $-\nabla_0$, $\Delta_0$ and $-\Delta_0$, which are balanced in $\Zn$ by Proposition~\ref{prop12}. The Steinhaus triangle $\nabla\partial S[0,3\lambda n-2]$ is balanced, since it is obtained from $\nabla S[0,3\lambda n-1]$ by rejecting the first row, which is a $3$-interlaced arithmetic progression with invertible common differences and of length $3\lambda n$ and thus contains $3\lambda$ times each element of $\Zn$. Similarly, the Pascal triangle $\Delta\partial S[-3\lambda n+1,3\lambda n-3]$ is balanced, since it is obtained from $\Delta\partial S[-3\lambda n,3\lambda n-2]$ by rejecting the last row, which is also balanced. The Steinhaus trapezoids (resp. the Pascal trapezoids) listed in this theorem can be seen as multiset differences of Steinhaus triangles (resp. Pascal triangles). Namely, we have
$$
\begin{array}{l}
\ST(S[0,3\lambda n-1],h) = \nabla S[0,3\lambda n-1]\setminus\nabla\partial^{h}S[0,3\lambda n-1-h],\\
\ST(\partial S[0,3\lambda n-2],h) = \nabla\partial S[0,3\lambda n-2]\setminus\nabla\partial^{h+1}S[0,3\lambda n-2-h],\\
\PT(\partial S[-m,m-2],h) = \Delta\partial S[-m,m-2]\setminus\Delta\partial S[-m+h,m-2-h].
\end{array}
$$
We have shown that these triangles are balanced. Therefore the trapezoids of this theorem also are balanced. Finally, the lozenge $\lozenge\partial S[-3\lambda n,3\lambda n-2]$ is the multiset union of the Pascal triangle $\Delta\partial S[-3\lambda n,3\lambda n-2]$ and the Steinhaus triangle $\nabla(-1)^{\lambda}\partial S[-3\lambda n,-2]=\nabla(-1)^{\lambda}\partial S[0,3\lambda n-2]$, which are balanced, for all integers $\lambda\geq1$.
\end{proof}

\section{The universal sequence modulo an odd number}

Let $\univ=\IAP((0,-1,1),(1,-2,1))$ be the universal sequence of integers introduced in Section~1. In this section, we refine Theorem~\ref{thm4} by studying this universal sequence modulo an odd number $n$, namely the sequence
$$
S = d\pi_n(\univ) = \IAP((0,-d,d),(d,-2d,d)),
$$
where $d$ is invertible in $\Zn$. It corresponds to the sequence $S$ of Theorem~\ref{thm4} with $a=0$. First, each element of its orbit $\mathcal{O}_{S}=\left(a_{i,j}\middle|a_{i+1,j}=a_{i,j}+a_{i,j+1},i\in\N,j\in\Z\right)$ can be expressed as a function of $d$.

\begin{prop}\label{prop13}
Let $n\in\N$ be odd and let $d\in\Zn$ be invertible. Consider the orbit $\mathcal{O}_{S}=\left(a_{i,j}\middle|a_{i+1,j}=a_{i,j}+a_{i,j+1},i\in\N,j\in\Z\right)$ of the sequence $S=\IAP((0,-d,d),(d,-2d,d))$ in $\Zn$. Then, for all $i,j\in\N$, we have
$$
a_{i,j} = (-1)^{i}\sum_{k>0}\binom{k}{j+2i-k}(-1)^{k}(k-i)d.
$$
\end{prop}

\begin{proof}
We begin by proving this equality for $i=0$. Let $(u_j)_{j\in\N}$ and $(v_j)_{j\in\N}$ be the sequences, in $\Zn$, defined by $u_j=\sum_{k>0}\binom{k}{j-k}{(-1)}^kkd$ and $v_j=\sum_{k>0}\binom{k}{j-k}{(-1)}^kd$, for all $j\in\N$. Then, for every integer $j\geq2$, we have
$$
u_{j} =\sum_{k>0}\left(\binom{k-1}{j-k-1}+\binom{k-1}{j-k}\right){(-1)}^k(k-1)d+\sum_{k>0}\binom{k}{j-k}{(-1)}^kd = -u_{j-2}-u_{j-1}+v_j.
$$
In the same way, we can prove that the sequence $(v_j)_{j\in\N}$ satisfies the relation $v_j+v_{j-1}+v_{j-2}=0$, for all integers $j\geq2$. It follows that $v_{3j}=d$, $v_{3j+1}=-d$ and $v_{3j+2}=0$, for all $j\in\N$. We complete the proof by induction on $j$. If we suppose that $u_{3j}=jd$, $u_{3j+1}=-(1+2j)d$ and $u_{3j+2}=(1+j)d$, then we obtain that $u_{3j+3}=-u_{3j+2}-u_{3j+1}+v_{3j+3}=(j+1)d$, $u_{3j+4}=-u_{3j+3}-u_{3j+2}+v_{3j+4}=-(3+2j)d$ and $u_{3j+5}=-u_{3j+4}-u_{3j+3}+v_{3j+5}=(2+j)d$. Therefore, we have $a_{0,j}=u_j=\sum_{k>0}\binom{k}{j-k}(-1)^kkd$, for all $j\in\N$, and this completes the proof for $i=0$. Finally, for all integers $i,j\geq1$, we obtain
$$
a_{i,j}
\begin{array}[t]{l}
\displaystyle = \sum_{l=0}^{i}\binom{i}{l}a_{0,j+l} = \sum_{l=0}^{i}\binom{i}{l}\sum_{k>0}\binom{k}{j+l-k}(-1)^{k}kd\\
\displaystyle = \sum_{k>0}\sum_{l=0}^{i}\binom{i}{k}\binom{k}{j+l-k}(-1)^kkd = \sum_{k>0}\binom{i+k}{j+i-k}(-1)^{k}kd\\
\displaystyle = (-1)^i\sum_{k>i}\binom{k}{j+2i-k}(-1)^k(k-i)d = (-1)^i\sum_{k>0}\binom{k}{j+2i-k}(-1)^k(k-i)d.
\end{array}
$$
\end{proof}

\begin{figure}[!h]
\begin{center}
\begin{pspicture}(7.57772228,2)
\multips(0,0)(2.16506351,0){3}{
\multips(0,1.5)(0.433012702,0){5}{\psline(0,0.375)(0.216506351,0.5)(0.433012702,0.375)}
\multips(0,1.5)(0.216506351,-0.375){5}{\psline(0,0.375)(0,0.125)(0.216506351,0)}
\multips(1.73205081,1.5)(-0.216506351,-0.375){5}{\psline(0.433012702,0.375)(0.433012702,0.125)(0.216506351,0)}}
\multips(1.29903811,0)(0.433012702,0){14}{\psline(0,0.125)(0.216506351,0)(0.433012702,0.125)}
\multips(6.49519053,1.625)(0.216506351,-0.375){4}{\psline(0,0)(0.216506351,-0.125)(0.216506351,-0.375)}
\rput(1.08253175,1.125){$\nabla_1$}
\rput(3.24759526,1.125){$\nabla_2$}
\rput(5.41265877,1.125){$\nabla_3$}
\rput(2.16506351,0.75){$\Delta_1$}
\rput(4.33012702,0.75){$\Delta_2$}
\rput(6.49519053,0.75){$\Delta_3$}
\end{pspicture}
\end{center}
\caption{\label{fig7}The elementary triangles $\nabla_1$, $\nabla_2$, $\nabla_3$, $\Delta_1$, $\Delta_2$ and $\Delta_3$.}
\end{figure}
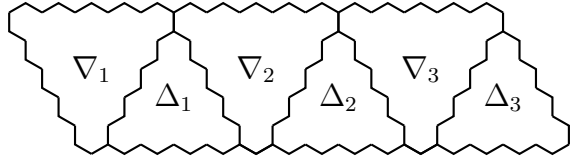

In the sequel of this section, we suppose that $n$ is an odd number and that $S$ is the universal sequence modulo $n$, that is $S=\IAP((0,-d,d),(d,-2d,d))$, where $d$ is an invertible element in $\Zn$. Let $\nabla_1$, $\nabla_2$ and $\nabla_3$ be the Steinhaus triangles of order $n$ associated with the sequences $S[0,n-1]$, $S[n,2n-1]$ and $S[2n,3n-1]$ respectively and let $\Delta_1$, $\Delta_2$ and $\Delta_3$ be their adjacent Pascal triangles of order $2n-3$, as depicted in Figure~\ref{fig7}, that are: $\nabla_1=\nabla S[0,n-1]$, $\nabla_2=\nabla S[n,2n-1]$, $\nabla_3=\nabla S[2n,3n-1]$, $\Delta_1=\Delta\partial S[1,2n-3]$, $\Delta_2=\Delta\partial S[n+1,3n-3]$ and $\Delta_3=\Delta\partial S[2n+1,4n-3]$.

We begin by showing that these triangles, or unions of them, are balanced in $\Zn$.

\begin{prop}\label{prop10}
Let $n\in\N$ be odd and let $d\in\Zn$ be invertible. Consider the universal sequence $S=\IAP((0,-d,d),(d,-2d,d))$ modulo $n$ and the elementary triangles $\nabla_1=\nabla S[0,n-1]$, $\nabla_2=\nabla S[n,2n-1]$, $\nabla_3=\nabla S[2n,3n-1]$, $\Delta_1=\Delta\partial S[1,2n-3]$, $\Delta_2=\Delta\partial S[n+1,3n-3]$ and $\Delta_3=\Delta\partial S[2n+1,4n-3]$. Then, the multisets $\nabla_2$, $\nabla_1\bigcup\nabla_3$, $\Delta_3$ and $\Delta_1\bigcup\Delta_2$ are balanced in $\Zn$.
\end{prop}

The proof of this proposition is based on the following lemma.

A finite sequence $S=(a_0,\ldots,a_{m-1})$ of length $m\geq1$ in $\Zn$ is said to be \textit{symmetric} if $a_j=a_{m-1-j}$ for all $j\in\{0,1,\ldots,m-1\}$.

\begin{lem}\label{lem4}
Let $n\in\N$ be odd and let $\nabla=\left\{a_{i,j}\middle|0\leq i\leq m-1,0\leq j\leq m-1-i\right\}$ be a Steinhaus triangle of order $m\geq1$ in $\Zn$. Then, the anti-diagonals $AD_{2j}$ and $AD_{2j+1}$ of $\nabla$ are respectively antisymmetric and symmetric for all integers $j$ such that $0\leq 2j\leq 2j+1\leq m-1$ if, and only if, we have $a_{i,i}=0$ for all $i\in\{0,1,\ldots,(m-1)/2\}$.
\end{lem}

\begin{proof}
For every $i\in\{0,1,\ldots,(m-1)/2\}$, if the sequence $AD_{2i}=(a_{0,2i},\ldots,a_{2i,0})$ is antisymmetric, then it follows that we have $2a_{i,i}=0$ and thus $a_{i,i}=0$, since $n$ is odd. Conversely, suppose now that $a_{i,i}=0$ for all $i\in\{0,1,\ldots,(m-1)/2\}$. We proceed by induction on $j$. For $j=0$, it is clear that $AD_{0}=(a_{0,0})=(0)$ is antisymmetric and that $AD_{1}=(a_{0,1},a_{1,0})=(a_1,a_1)$ is symmetric. Suppose that the result is true for $j-1$, i.e., that the sequences $AD_{2j-2}$ and $AD_{2j-1}$ are respectively antisymmetric and symmetric, and prove it for $j$. We begin by showing that $a_{j-k,j+k}=-a_{j+k,j-k}$ for all $k\in\{0,1,\ldots,j\}$. For $k=0$, it comes from hypothesis $a_{j,j}=0$. Suppose it is true for all integers in $\{0,\ldots,k-1\}$. Since $a_{j-k,j+k-1}=a_{j+k-1,j-k}$ by symmetry of $AD_{2j-1}$, we obtain that $a_{j-k,j+k}=a_{j-(k-1),j+k-1}-a_{j-k,j+k-1}=-a_{j+k-1,j-(k-1)}-a_{j+k-1,j-k}=-a_{j+k,j-k}$ and thus $AD_{2j}$ is antisymmetric. We now prove that $a_{j-k,j+1+k}=a_{j+1+k,j-k}$ for all $k\in\{0,1,\ldots,j\}$. For $k=0$, it follows from the equality $a_{j+1,j}=a_{j,j}+a_{j,j+1}=a_{j,j+1}$. Suppose it is true for all integers in $\{0,\ldots,k-1\}$. Since $a_{j-k,j+k}=-a_{j+k,j-k}$ by antisymmetry of $AD_{2j}$, we have $a_{j-k,j+k+1}=a_{j-k+1,j+k}-a_{j-k,j+k}=a_{j+k,j-(k-1)}+a_{j+k,j-k}=a_{j+k+1,j-k}$ and thus $AD_{2j+1}$ is symmetric. This concludes the proof.
\end{proof}

\begin{proof}[Proof of Proposition~\ref{prop10}]
First, we consider the Steinhaus triangle $\nabla_0=\nabla S[0,3n-1]$ of order $3n$ and the Pascal triangle $\Delta_0=\Delta\partial S[1,6n-3]$ of order $6n-3$, which are balanced by Proposition~\ref{prop12}. If we denote by $\mathcal{O}_{S}=\left(a_{i,j}\middle|a_{i+1,j}=a_{i,j}+a_{i,j+1},i\in\N,j\in\Z\right)$ the orbit associated with the universal sequence $S=\IAP((0,-d,d),(d,-2d,d))$ in $\Zn$, then Proposition~\ref{prop13} implies that we have $a_{n,j}=-a_{0,2n+j}$, $a_{2n,j}=a_{0,n+j}$ and $a_{3n,j}=-a_{0,j}$ for all $j\in\Z$. Moreover, we have $a_{0,3n+j}=a_{0,j}$ for all $j\in\Z$, since the sequence $S$ is periodic of period $3n$. This leads to the following decomposition of $\nabla_0$ and $\Delta_0$ into elementary triangles $\nabla_1$, $\nabla_2$, $\nabla_3$, $\Delta_1$, $\Delta_2$ and $\Delta_3$:

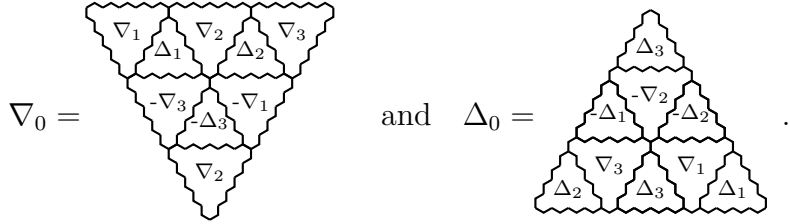
\begin{center}
$\nabla_0=$\!\!\begin{tabular}{c}\begin{pspicture}(3.24759526,2.875)
\multips(0,1.875)(1.08253175,0){3}{
\multips(0,0.75)(0.216506351,0){5}{\psline(0,0.1875)(0.108253175,0.25)(0.216506351,0.1875)}
\multips(0,0.75)(0.108253175,-0.1875){5}{\psline(0,0.1875)(0,0.0625)(0.108253175,0)}
\multips(0.866025404,0.75)(-0.108253175,-0.1875){5}{\psline(0.216506351,0.1875)(0.216506351,0.0625)(0.108253175,0)}}
\multips(0.541265877,0.9375)(1.08253175,0){2}{
\multips(0,0.75)(0.216506351,0){5}{\psline(0,0.1875)(0.108253175,0.25)(0.216506351,0.1875)}
\multips(0,0.75)(0.108253175,-0.1875){5}{\psline(0,0.1875)(0,0.0625)(0.108253175,0)}
\multips(0.866025404,0.75)(-0.108253175,-0.1875){5}{\psline(0.216506351,0.1875)(0.216506351,0.0625)(0.108253175,0)}}
\multips(1.08253175,0)(1.08253175,0){1}{
\multips(0,0.75)(0.216506351,0){5}{\psline(0,0.1875)(0.108253175,0.25)(0.216506351,0.1875)}
\multips(0,0.75)(0.108253175,-0.1875){5}{\psline(0,0.1875)(0,0.0625)(0.108253175,0)}
\multips(0.866025404,0.75)(-0.108253175,-0.1875){5}{\psline(0.216506351,0.1875)(0.216506351,0.0625)(0.108253175,0)}}
\rput(0.541265877,2.5){\scriptsize$\nabla_1$}
\rput(1.62379763,2.5){\scriptsize$\nabla_2$}
\rput(2.70632939,2.5){\scriptsize$\nabla_3$}
\rput(1.08253175,2.25){\scriptsize$\Delta_1$}
\rput(2.16506351,2.25){\scriptsize$\Delta_2$}
\rput(1.08253175,1.5625){\scriptsize-$\nabla_3$}
\rput(2.16506351,1.5625){\scriptsize-$\nabla_1$}
\rput(1.62379763,1.3125){\scriptsize-$\Delta_3$}
\rput(1.62379763,0.625){\scriptsize$\nabla_2$}
\end{pspicture}\end{tabular}\quad and\quad 
$\Delta_0=$\!\!\!\begin{tabular}{c}\begin{pspicture}(3.03108891,2.6875)

\multips(0.541265877,0.9375)(0,0){1}{
\multips(0.974278579,1.75)(-0.108253175,-0.1875){9}{\psline(0,0)(-0.108253175,-0.0625)(-0.108253175,-0.1875)}
\multips(0.974278579,1.75)(0.108253175,-0.1875){9}{\psline(0,0)(0.108253175,-0.0625)(0.108253175,-0.1875)}
\multips(0,0)(0.216506351,0){9}{\psline(0,0.0625)(0.108253175,0)(0.216506351,0.0625)}
\multips(0.541265877,0.9375)(0.216506351,0){4}{\psline(0,0.0625)(0.108253175,0)(0.216506351,0.0625)}
\multips(0.433012702,0.8125)(0.108253175,-0.1875){4}{\psline(0,0)(0.108253175,-0.0625)(0.108253175,-0.1875)}
\multips(1.51554446,0.8125)(-0.108253175,-0.1875){4}{\psline(0,0)(-0.108253175,-0.0625)(-0.108253175,-0.1875)}}

\multips(0,0)(1.08253175,0){2}{
\multips(0.974278579,1.75)(-0.108253175,-0.1875){9}{\psline(0,0)(-0.108253175,-0.0625)(-0.108253175,-0.1875)}
\multips(0.974278579,1.75)(0.108253175,-0.1875){9}{\psline(0,0)(0.108253175,-0.0625)(0.108253175,-0.1875)}
\multips(0,0)(0.216506351,0){9}{\psline(0,0.0625)(0.108253175,0)(0.216506351,0.0625)}
\multips(0.541265877,0.9375)(0.216506351,0){4}{\psline(0,0.0625)(0.108253175,0)(0.216506351,0.0625)}
\multips(0.433012702,0.8125)(0.108253175,-0.1875){4}{\psline(0,0)(0.108253175,-0.0625)(0.108253175,-0.1875)}
\multips(1.51554446,0.8125)(-0.108253175,-0.1875){4}{\psline(0,0)(-0.108253175,-0.0625)(-0.108253175,-0.1875)}}

\rput(1.51554446,2.1875){\scriptsize$\Delta_3$}
\rput(0.974278579,1.3125){\scriptsize-$\Delta_1$}
\rput(2.05681033,1.3125){\scriptsize-$\Delta_2$}
\rput(1.51554446,1.5625){\scriptsize-$\nabla_2$}
\rput(0.433012702,0.3125){\scriptsize$\Delta_2$}
\rput(1.51554446,0.3125){\scriptsize$\Delta_3$}
\rput(2.59807621,0.3125){\scriptsize$\Delta_1$}
\rput(0.974278579,0.625){\scriptsize$\nabla_3$}
\rput(2.05681033,0.625){\scriptsize$\nabla_1$}
\end{pspicture}\end{tabular}.
\end{center}

For every $k\in\{0,1,2,3\}$, we denote by $D_{j}(\nabla_k)$ and $AD_j(\nabla_k)$ the $j$th diagonal and the $j$th anti-diagonal of $\nabla_k$, for every $j\in\{0,1,\ldots,n-1\}$, and by $D_{j}(\Delta_k)$ and $AD_j(\Delta_k)$ the $j$th diagonal and the $j$th anti-diagonal of $\Delta_k$, for every $j\in\{0,1,\ldots,n-2\}$. Since we have $a_{i,i}=0$ for all $i\in\N$, from the general expression of $a_{i,j}$ appearing in the proof of Proposition~\ref{prop4}, it follows, from Lemma~\ref{lem4}, that the sequences $AD_{2j}(\nabla_0)$ and $AD_{2j+1}(\nabla_0)$ are respectively antisymmetric and symmetric, for all integers $j$ such that $0\leq 2j\leq 2j+1\leq 3n-1$.
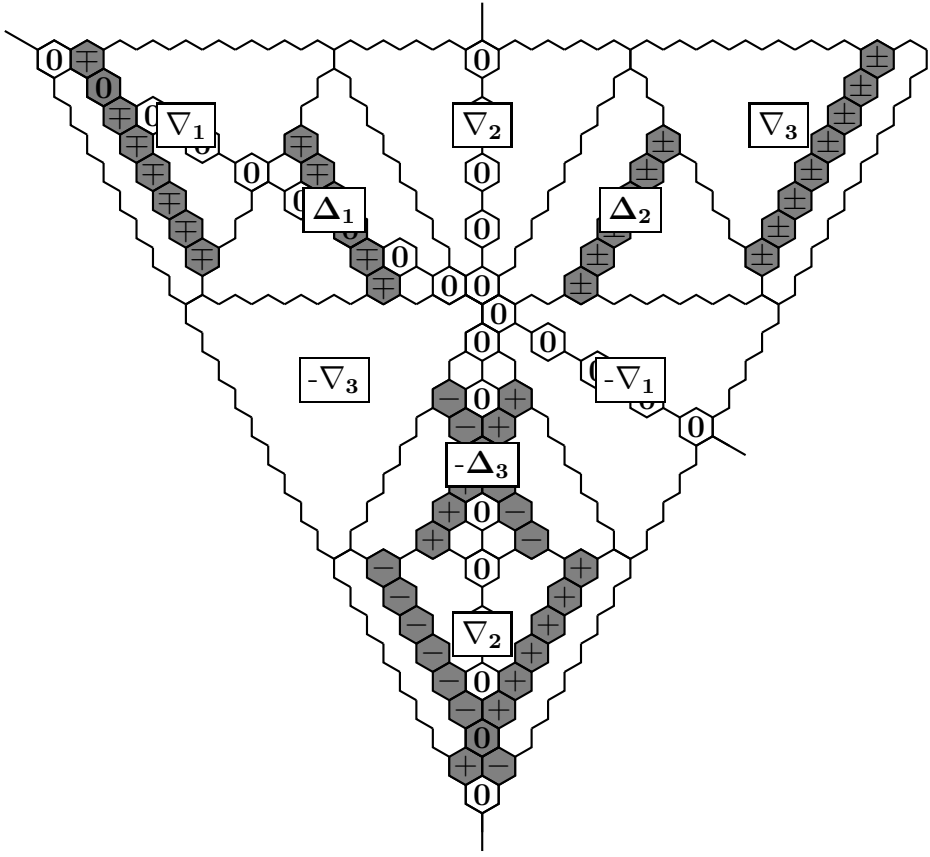
\begin{figure}[!h]
\begin{center}
\begin{pspicture}(11.691343,11.25)

\multips(0,7.25)(3.89711432,0){3}{
\multips(0,3)(0.433012702,0){9}{\psline(0,0.375)(0.216506351,0.5)(0.433012702,0.375)}
\multips(0,3)(0.216506351,-0.375){9}{\psline(0,0.375)(0,0.125)(0.216506351,0)}
\multips(3.46410162,3)(-0.216506351,-0.375){9}{\psline(0.433012702,0.375)(0.433012702,0.125)(0.216506351,0)}}
\multips(1.94855716,3.875)(3.89711432,0){2}{
\multips(0,3)(0.433012702,0){9}{\psline(0,0.375)(0.216506351,0.5)(0.433012702,0.375)}
\multips(0,3)(0.216506351,-0.375){9}{\psline(0,0.375)(0,0.125)(0.216506351,0)}
\multips(3.46410162,3)(-0.216506351,-0.375){9}{\psline(0.433012702,0.375)(0.433012702,0.125)(0.216506351,0)}}
\multips(3.89711432,0.5)(3.89711432,0){1}{
\multips(0,3)(0.433012702,0){9}{\psline(0,0.375)(0.216506351,0.5)(0.433012702,0.375)}
\multips(0,3)(0.216506351,-0.375){9}{\psline(0,0.375)(0,0.125)(0.216506351,0)}
\multips(3.46410162,3)(-0.216506351,-0.375){9}{\psline(0.433012702,0.375)(0.433012702,0.125)(0.216506351,0)}}

\multips(5.41265877,0.875)(0.216506351,0.375){8}{\pspolygon*[linecolor=Gray](0,0.125)(0,0.375)(0.216506351,0.5)(0.433012702,0.375)(0.433012702,0.125)(0.216506351,0)}
\multips(5.41265877,0.875)(0.216506351,0.375){8}{\pspolygon(0,0.125)(0,0.375)(0.216506351,0.5)(0.433012702,0.375)(0.433012702,0.125)(0.216506351,0)}
\rput(5.62916512,1.125){$\mathbf{+}$}
\multirput(6.06217783,1.875)(0.216506351,0.375){6}{$\mathbf{+}$}

\multips(9.30977309,7.625)(0.216506351,0.375){8}{\pspolygon*[linecolor=Gray](0,0.125)(0,0.375)(0.216506351,0.5)(0.433012702,0.375)(0.433012702,0.125)(0.216506351,0)}
\multips(9.30977309,7.625)(0.216506351,0.375){8}{\pspolygon(0,0.125)(0,0.375)(0.216506351,0.5)(0.433012702,0.375)(0.433012702,0.125)(0.216506351,0)}
\multirput(9.52627944,7.875)(0.216506351,0.375){8}{$\mathbf{\pm}$}

\multips(0.433012702,10.25)(0.216506351,-0.375){8}{\pspolygon*[linecolor=Gray](0,0.125)(0,0.375)(0.216506351,0.5)(0.433012702,0.375)(0.433012702,0.125)(0.216506351,0)}
\multips(0.433012702,10.25)(0.216506351,-0.375){8}{\pspolygon(0,0.125)(0,0.375)(0.216506351,0.5)(0.433012702,0.375)(0.433012702,0.125)(0.216506351,0)}
\rput(0.649519053,10.5){$\mathbf{\mp}$}
\multirput(1.08253175,9.75)(0.216506351,-0.375){6}{$\mathbf{\mp}$}

\multips(5.84567148,0.875)(-0.216506351,0.375){8}{\pspolygon*[linecolor=Gray](0,0.125)(0,0.375)(0.216506351,0.5)(0.433012702,0.375)(0.433012702,0.125)(0.216506351,0)}
\multips(5.84567148,0.875)(-0.216506351,0.375){8}{\pspolygon(0,0.125)(0,0.375)(0.216506351,0.5)(0.433012702,0.375)(0.433012702,0.125)(0.216506351,0)}
\rput(6.06217783,1.125){$\mathbf{-}$}
\multirput(5.62916512,1.875)(-0.216506351,0.375){6}{$\mathbf{-}$}

\multips(4.97964607,3.875)(0.216506351,0.375){6}{\pspolygon*[linecolor=Gray](0,0.125)(0,0.375)(0.216506351,0.5)(0.433012702,0.375)(0.433012702,0.125)(0.216506351,0)}
\multips(4.97964607,3.875)(0.216506351,0.375){6}{\pspolygon(0,0.125)(0,0.375)(0.216506351,0.5)(0.433012702,0.375)(0.433012702,0.125)(0.216506351,0)}
\multirput(5.19615242,4.125)(0.216506351,0.375){3}{$\mathbf{+}$}
\multirput(6.06217783,5.625)(0.216506351,0.375){2}{$\mathbf{+}$}

\multips(6.92820323,7.25)(0.216506351,0.375){6}{\pspolygon*[linecolor=Gray](0,0.125)(0,0.375)(0.216506351,0.5)(0.433012702,0.375)(0.433012702,0.125)(0.216506351,0)}
\multips(6.92820323,7.25)(0.216506351,0.375){6}{\pspolygon(0,0.125)(0,0.375)(0.216506351,0.5)(0.433012702,0.375)(0.433012702,0.125)(0.216506351,0)}
\multirput(7.14470958,7.5)(0.216506351,0.375){6}{$\mathbf{\pm}$}

\multips(4.33012702,7.25)(-0.216506351,0.375){6}{\pspolygon*[linecolor=Gray](0,0.125)(0,0.375)(0.216506351,0.5)(0.433012702,0.375)(0.433012702,0.125)(0.216506351,0)}
\multips(4.33012702,7.25)(-0.216506351,0.375){6}{\pspolygon(0,0.125)(0,0.375)(0.216506351,0.5)(0.433012702,0.375)(0.433012702,0.125)(0.216506351,0)}
\multirput(4.54663337,7.5)(-0.216506351,0.375){2}{$\mathbf{\mp}$}
\multirput(3.89711432,8.625)(-0.216506351,0.375){3}{$\mathbf{\mp}$}

\multips(6.27868418,3.875)(-0.216506351,0.375){6}{\pspolygon*[linecolor=Gray](0,0.125)(0,0.375)(0.216506351,0.5)(0.433012702,0.375)(0.433012702,0.125)(0.216506351,0)}
\multips(6.27868418,3.875)(-0.216506351,0.375){6}{\pspolygon(0,0.125)(0,0.375)(0.216506351,0.5)(0.433012702,0.375)(0.433012702,0.125)(0.216506351,0)}
\multirput(6.49519053,4.125)(-0.216506351,0.375){3}{$\mathbf{-}$}
\multirput(5.62916512,5.625)(-0.216506351,0.375){2}{$\mathbf{-}$}

\psline(5.84567148,0)(5.84567148,0.5)
\psline(5.84567148,11.25)(5.84567148,10.75)
\multips(5.62916512,10.25)(0,-0.75){14}{\pspolygon(0,0.125)(0,0.375)(0.216506351,0.5)(0.433012702,0.375)(0.433012702,0.125)(0.216506351,0)\psline(0.216506351,0)(0.216506351,-0.25)}
\multirput(5.84567148,10.5)(0,-0.75){14}{$\mathbf{0}$}

\psline(0,10.625)(-0.433012702,10.875)
\psline(8.87676039,5.5)(9.30977309,5.25)
\multips(0,10.25)(0.649519053,-0.375){14}{\pspolygon(0,0.125)(0,0.375)(0.216506351,0.5)(0.433012702,0.375)(0.433012702,0.125)(0.216506351,0)\psline(0.433012702,0.125)(0.649519053,0)}
\multirput(0.216506351,10.5)(0.649519053,-0.375){14}{$\mathbf{0}$}

\rput(1.94855716,9.625){\fcolorbox{black}{white}{$\mathbf{\nabla_1}$}}
\rput(5.84567148,9.625){\fcolorbox{black}{white}{$\mathbf{\nabla_2}$}}
\rput(9.74278579,9.625){\fcolorbox{black}{white}{$\mathbf{\nabla_3}$}}

\rput(3.89711432,8.5){\fcolorbox{black}{white}{$\mathbf{\Delta_1}$}}
\rput(7.79422863,8.5){\fcolorbox{black}{white}{$\mathbf{\Delta_2}$}}

\rput(3.89711432,6.25){\fcolorbox{black}{white}{-$\mathbf{\nabla_3}$}}
\rput(7.79422863,6.25){\fcolorbox{black}{white}{-$\mathbf{\nabla_1}$}}

\rput(5.84567148,5.125){\fcolorbox{black}{white}{-$\mathbf{\Delta_3}$}}

\rput(5.84567148,2.875){\fcolorbox{black}{white}{$\mathbf{\nabla_2}$}}

\end{pspicture}
\end{center}
\caption{\label{fig6}The Steinhaus triangle $\nabla_0$.}
\end{figure}
This implies the following equalities on the multiplicity functions of the anti-diagonals of $\nabla_2$ and $\nabla_3$:
$$
\m_{AD_{2j}(\nabla_3)}(x)=\m_{AD_{2j}(\nabla_2)}(-x)\quad\text{and}\quad\m_{AD_{2j+1}(\nabla_3)}(x)=\m_{AD_{2j+1}(\nabla_2)}(x),
$$
for all $x\in\Zn$ and for all integers $j$ such that $0\leq 2j\leq2j+1\leq n-1$. Moreover, we know, from Proposition~\ref{prop8}, that the sequence $S[0,3n-1]$ is antisymmetric and, thus, all the rows of $\nabla_0$ are also antisymmetric by Proposition~\ref{prop6}. Therefore, we have
$$
\m_{AD_{j}(\nabla_3)}(x)=\m_{D_{n-1-j}(\nabla_1)}(-x)\quad\text{and}\quad\m_{AD_{j}(\nabla_2)}(x)=\m_{D_{n-1-j}(\nabla_2)}(-x),
$$
for all $x\in\Zn$ and for all $j\in\{0,1,\ldots,n-1\}$. This leads to the equality
$$
\m_{\nabla_1}(x)+\m_{\nabla_3}(x)
\begin{array}[t]{l}
\displaystyle = \sum_{j=0}^{n-1}\m_{D_{j}(\nabla_1)}(x) + \sum_{j=0}^{n-1}\m_{AD_j(\nabla_3)}(x) = \sum_{j=0}^{n-1}\m_{AD_{n-1-j}(\nabla_3)}(-x) + \sum_{j=0}^{n-1}\m_{AD_j(\nabla_3)}(x)\\
\displaystyle = \sum_{j=0}^{n-1}\left(\m_{AD_{j}(\nabla_3)}(-x)+\m_{AD_j(\nabla_3)}(x)\right) = \sum_{j=0}^{n-1}\left(\m_{AD_{j}(\nabla_2)}(-x)+\m_{AD_j(\nabla_2)}(x)\right)\\
\displaystyle = \sum_{j=0}^{n-1}\m_{AD_j(\nabla_2)}(x) + \sum_{j=0}^{n-1}\m_{AD_{n-1-j}(\nabla_2)}(-x) = 2\m_{\nabla_2}(x),\\
\end{array}
$$
for all $x\in\Zn$. Similarly, if we consider the diagonals and the anti-diagonals of the triangles $\Delta_1$, $\Delta_2$ and $-\Delta_3$, as depicted in Figure~\ref{fig6}, then we obtain that $\m_{\Delta_1}+\m_{\Delta_2}=2\m_{-\Delta_3}$. The antisymmetry in $\nabla_0$ also implies the following equalities: $\m_{\nabla_1}=\m_{-\nabla_3}$, $\m_{\nabla_3}=\m_{-\nabla_1}$, $\m_{\nabla_2}=\m_{-\nabla_2}$, $\m_{\Delta_1}=\m_{-\Delta_2}$, $\m_{\Delta_2}=\m_{-\Delta_1}$ and $\m_{\Delta_3}=\m_{-\Delta_3}$. Therefore, the multiplicity functions of these elementary triangles verify the following equations:
$$
\begin{array}{l}
\m_{\nabla_1}+\m_{\nabla_3}=2\m_{\nabla_2},\\
\m_{\Delta_1}+\m_{\Delta_2}=2\m_{\Delta_3}.
\end{array}
$$
Finally, since the triangles $\nabla_0$ and $\Delta_0$ are balanced in $\Zn$, it follows that the multiplicity functions $\m_{\nabla_2}$ and $\m_{\Delta_3}$ are solutions of the following system of equations
$$
\begin{array}{l}
\displaystyle 6\m_{\nabla_2}+3\m_{\Delta_3}=\m_{\nabla_0}=\frac{1}{n}\binom{3n+1}{2},\\[3ex]
\displaystyle 3\m_{\nabla_2}+6\m_{\Delta_3}=\m_{\Delta_0}=\frac{1}{n}\binom{3n}{2}.
\end{array}
$$
We conclude that the triangles $\nabla_2$, $\Delta_3$ and the multisets $\nabla_1\bigcup\nabla_3$ and $\Delta_1\bigcup\Delta_2$ are balanced.
\end{proof}

We are now ready to prove Theorem~\ref{thm0}, the main result of this paper.

\begin{thm}\label{thm5}
Let $n\in\N$ be odd and let $d\in\Zn$ be invertible. Then, the following Steinhaus figures, contained in the orbit associated with the universal sequence $S=\IAP((0,-d,d),(d,-2d,d))$ in $\Zn$, are balanced:
\begin{itemize}
\item
the Steinhaus triangles $\nabla S[m,2m-1]$, for every $m\equiv0\pmod{n}$, and $\nabla\partial S[0,m-1]$, for every $m\equiv-1\pmod{3n}$,
\item
the Steinhaus trapezoids $\ST(S[m,2m-1],h)$, for every $m\equiv0\pmod{n}$ and for every $h\equiv0\pmod{n}$ or $h\equiv m+1\pmod{3n}$, and $\ST(\partial S[0,m-1],h)$, for every $m\equiv-1\pmod{3n}$ and for every $h\equiv-1\pmod{n}$ or $h\equiv 0\pmod{3n}$,
\item
the Pascal triangle $\Delta\partial S[-m,m-2]$, for every $m\equiv-1\pmod{n}$ or $m\equiv0\pmod{3n}$,
\item
the Pascal trapezoid $\PT(\partial S[-m,m-2],h)$, for every $m\equiv-1\pmod{n}$ or $m\equiv0\pmod{3n}$ and for every $h\equiv m\pmod{n}$ or $h\equiv m+1\pmod{3n}$,
\item
the lozenge $\lozenge\partial S[-m,m-2]$, for every $m\equiv0\pmod{n}$.
\end{itemize}
\end{thm}

\begin{proof}
The Steinhaus figures of this theorem are unions of the multisets $\pm\nabla_2$, $\pm(\nabla_1\bigcup\nabla_3)$, $\pm\Delta_3$ and $\pm(\Delta_1\bigcup\Delta_2)$, which are balanced in $\Zn$ by Proposition~\ref{prop10}. More precisely, let $\lambda$ be a positive integer. We know, from Theorem~\ref{thm4}, that the Steinhaus triangles $\nabla S[3\lambda n,6\lambda n-1]$, of order $3\lambda n$, and $\nabla\partial S[0,3\lambda n-2]$, of order $3\lambda n-1$, are balanced. As depicted in Figure~\ref{fig11}, the Steinhaus triangle $\nabla S[(3\lambda+1)n,(6\lambda+2)n-1]$, of order $(3\lambda+1)n$, is the union of $\lambda+1$ triangles $\nabla_2$, $\lambda$ multisets $\nabla_1\bigcup\nabla_3$, $\lambda$ triangles $\Delta_3$, $\lambda$ multisets $\Delta_1\bigcup\Delta_2$ and the Steinhaus triangle $\nabla\partial^nS[(3\lambda+1)n,(6\lambda+2)n-1]=-\nabla S[3\lambda n,6\lambda n-1]$. This leads to the equality
$$
\m_{\nabla S[(3\lambda+1)n,(6\lambda+2)n-1]} = (3\lambda+1)\m_{\nabla_2} + 3\lambda\m_{\Delta_3} + \m_{\nabla S[3\lambda n,6\lambda n-1]}.
$$
Similarly, the Steinhaus triangle $\nabla S[(3\lambda+2)n,(6\lambda+4)n-1]$, of order $(3\lambda+2)n$, is the union of $\lambda$ triangles $\nabla_2$, $\lambda+1$ multisets $\nabla_1\bigcup\nabla_3$, $\lambda+1$ triangles $\Delta_3$, $\lambda$ multisets $\Delta_1\bigcup\Delta_2$ and the Steinhaus triangle $\nabla\partial^nS[(3\lambda+2)n,(6\lambda+4)n-1]=-\nabla S[(3\lambda+1)n,(6\lambda+2)n-1]$. Therefore, we obtain
$$
\m_{\nabla S[(3\lambda+2)n,(6\lambda+4)n-1]} = (3\lambda+2)\m_{\nabla_2} + (3\lambda+1)\m_{\Delta_3} + \m_{\nabla S[(3\lambda+1)n,(6\lambda+2)n-1]}.
$$
This completes the proof that the Steinhaus triangle $\nabla S[m,2m-1]$ is balanced for all $m\equiv0\pmod{n}$. A similar decomposition shows that the Pascal triangle $\Delta\partial S[-m,m-2]$ is balanced for all $m\equiv-1\pmod{n}$. First, we know, from Theorem~\ref{thm4}, that the Pascal triangles $\Delta\partial S[-3\lambda n,3\lambda n-2]$, of order $6\lambda n-1$, and $\Delta\partial S[-3\lambda n+1,3\lambda n-3]$, of order $6\lambda n-3$, are balanced in $\Zn$. The other cases come from the decomposition into elementary triangles, as depicted in Figure~\ref{fig12}, which implies the following equalities:
$$
\begin{array}{l}
\m_{\Delta\partial S[-(3\lambda+1)n+1,(3\lambda+1)n-3]} = 3\lambda\m_{\nabla_2} + (3\lambda+1)\m_{\Delta_3} + \m_{\Delta\partial S[-3\lambda n+1,3\lambda n-3]},\\[2ex]
\m_{\Delta\partial S[-(3\lambda+2)n+1,(3\lambda+2)n-3]} = (3\lambda+1)\m_{\nabla_2} + (3\lambda+2)\m_{\Delta_3} + \m_{\Delta\partial S[-(3\lambda+1)n+1,(3\lambda+1)n-3]}.
\end{array}
$$
The Steinhaus trapezoids (resp. the Pascal trapezoids) listed in this theorem can be seen as multiset differences of Steinhaus triangles (resp. Pascal triangles). Namely, we have
$$
\begin{array}{l}
\ST(S[m,2m-1],h) = \nabla S[m,2m-1]\setminus\nabla\partial^{h}S[m,2m-1-h],\\
\ST(\partial S[0,m-1],h) = \nabla\partial S[0,m-1]\setminus\nabla\partial^{h+1}S[0,m-1-h],\\
\PT(\partial S[-m,m-2],h) = \Delta\partial S[-m,m-2]\setminus\Delta\partial S[-m+h,m-2-h].
\end{array}
$$
We have shown that these triangles are balanced. Therefore the trapezoids of this theorem also are balanced. Finally, the lozenge $\lozenge\partial S[-m+1,m-2]$ is the multiset union of the Pascal triangle $\Delta\partial S[-m+1,m-3]$ and of the Steinhaus triangle $(-1)^{m}\nabla S[m,2m-1]$, which are balanced in $\Zn$ for all $m\equiv0\pmod{n}$.

\begin{figure}[!p]
\begin{center}
\begin{pspicture}(10.8253175,9.4375)

\multips(0,9.1875)(0.216506351,0){50}{\psline(0,0.1875)(0.108253175,0.25)(0.216506351,0.1875)}
\multips(0,9.1875)(0.108253175,-0.1875){50}{\psline(0,0.1875)(0,0.0625)(0.108253175,0)}
\multips(10.6088112,9.1875)(-0.108253175,-0.1875){50}{\psline(0.216506351,0.1875)(0.216506351,0.0625)(0.108253175,0)}

\multips(2.16506351,5.4375)(0.216506351,0){30}{\psline(0,0.1875)(0.108253175,0.25)(0.216506351,0.1875)}

\multips(0,8.4375)(1.08253175,0){10}{
\multips(0,0.75)(0.216506351,0){5}{\psline(0,0.1875)(0.108253175,0.25)(0.216506351,0.1875)}
\multips(0,0.75)(0.108253175,-0.1875){5}{\psline(0,0.1875)(0,0.0625)(0.108253175,0)}
\multips(0.866025404,0.75)(-0.108253175,-0.1875){5}{\psline(0.216506351,0.1875)(0.216506351,0.0625)(0.108253175,0)}}
\multips(0.541265877,7.5)(1.08253175,0){9}{
\multips(0,0.75)(0.216506351,0){5}{\psline(0,0.1875)(0.108253175,0.25)(0.216506351,0.1875)}
\multips(0,0.75)(0.108253175,-0.1875){5}{\psline(0,0.1875)(0,0.0625)(0.108253175,0)}
\multips(0.866025404,0.75)(-0.108253175,-0.1875){5}{\psline(0.216506351,0.1875)(0.216506351,0.0625)(0.108253175,0)}}
\multips(1.08253175,6.5625)(1.08253175,0){8}{
\multips(0,0.75)(0.216506351,0){5}{\psline(0,0.1875)(0.108253175,0.25)(0.216506351,0.1875)}
\multips(0,0.75)(0.108253175,-0.1875){5}{\psline(0,0.1875)(0,0.0625)(0.108253175,0)}
\multips(0.866025404,0.75)(-0.108253175,-0.1875){5}{\psline(0.216506351,0.1875)(0.216506351,0.0625)(0.108253175,0)}}
\multips(1.62379763,5.625)(1.08253175,0){7}{
\multips(0,0.75)(0.216506351,0){5}{\psline(0,0.1875)(0.108253175,0.25)(0.216506351,0.1875)}
\multips(0,0.75)(0.108253175,-0.1875){5}{\psline(0,0.1875)(0,0.0625)(0.108253175,0)}
\multips(0.866025404,0.75)(-0.108253175,-0.1875){5}{\psline(0.216506351,0.1875)(0.216506351,0.0625)(0.108253175,0)}}

\multirput(0.541265877,9)(3.24759526,0){4}{\scriptsize$\nabla_2$}
\multirput(1.62379763,9)(3.24759526,0){3}{\scriptsize$\nabla_3$}
\multirput(2.70632939,9)(3.24759526,0){3}{\scriptsize$\nabla_1$}

\multirput(1.08253175,8.8125)(3.24759526,0){3}{\scriptsize$\Delta_2$}
\multirput(2.16506351,8.8125)(3.24759526,0){3}{\scriptsize$\Delta_3$}
\multirput(3.24759526,8.8125)(3.24759526,0){3}{\scriptsize$\Delta_1$}

\multirput(1.08253175,8.0625)(3.24759526,0){3}{\scriptsize-$\nabla_1$}
\multirput(2.16506351,8.0625)(3.24759526,0){3}{\scriptsize-$\nabla_2$}
\multirput(3.24759526,8.0625)(3.24759526,0){3}{\scriptsize-$\nabla_3$}

\multirput(1.62379763,7.875)(3.24759526,0){3}{\scriptsize-$\Delta_1$}
\multirput(2.70632939,7.875)(3.24759526,0){3}{\scriptsize-$\Delta_2$}
\multirput(3.78886114,7.875)(3.24759526,0){2}{\scriptsize-$\Delta_3$}

\multirput(1.62379763,7.125)(3.24759526,0){3}{\scriptsize$\nabla_3$}
\multirput(2.70632939,7.125)(3.24759526,0){3}{\scriptsize$\nabla_1$}
\multirput(3.78886114,7.125)(3.24759526,0){2}{\scriptsize$\nabla_2$}

\multirput(2.16506351,6.9375)(3.24759526,0){3}{\scriptsize$\Delta_3$}
\multirput(3.24759526,6.9375)(3.24759526,0){2}{\scriptsize$\Delta_1$}
\multirput(4.33012702,6.9375)(3.24759526,0){2}{\scriptsize$\Delta_2$}

\multirput(2.16506351,6.1875)(3.24759526,0){3}{\scriptsize-$\nabla_2$}
\multirput(3.24759526,6.1875)(3.24759526,0){2}{\scriptsize-$\nabla_3$}
\multirput(4.33012702,6.1875)(3.24759526,0){2}{\scriptsize-$\nabla_1$}

\multirput(2.70632939,6)(3.24759526,0){2}{\scriptsize-$\Delta_2$}
\multirput(3.78886114,6)(3.24759526,0){2}{\scriptsize-$\Delta_3$}
\multirput(4.8713929,6)(3.24759526,0){2}{\scriptsize-$\Delta_1$}

\rput(5.41265877,3.375){\scriptsize$\nabla S[3\lambda n,6\lambda n-1]$}

\end{pspicture}
\end{center}
\caption{\label{fig11}The Steinhaus triangle $\nabla S[m,2m-1]$ for $m\equiv0\pmod{n}$.}
\end{figure}
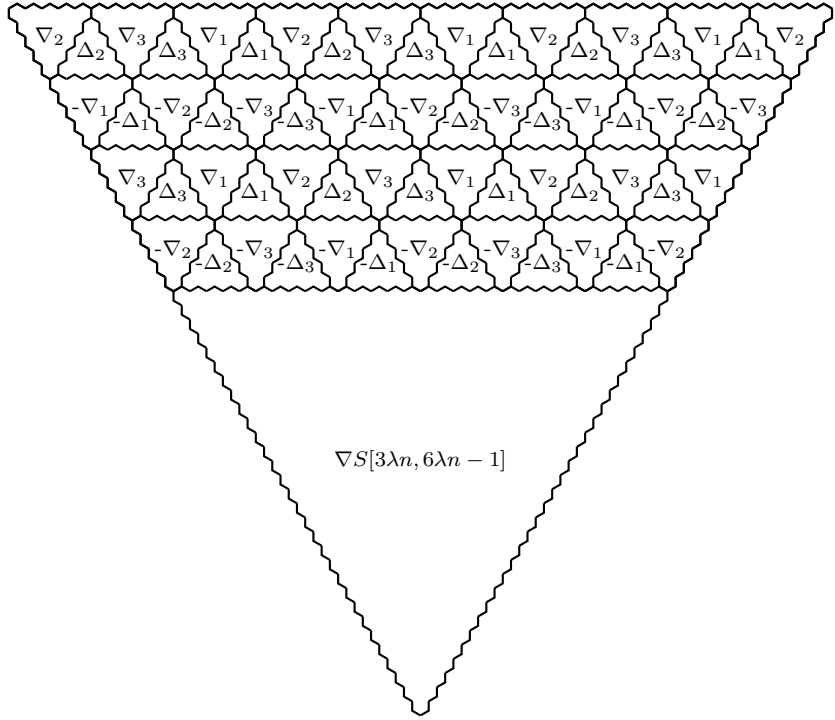

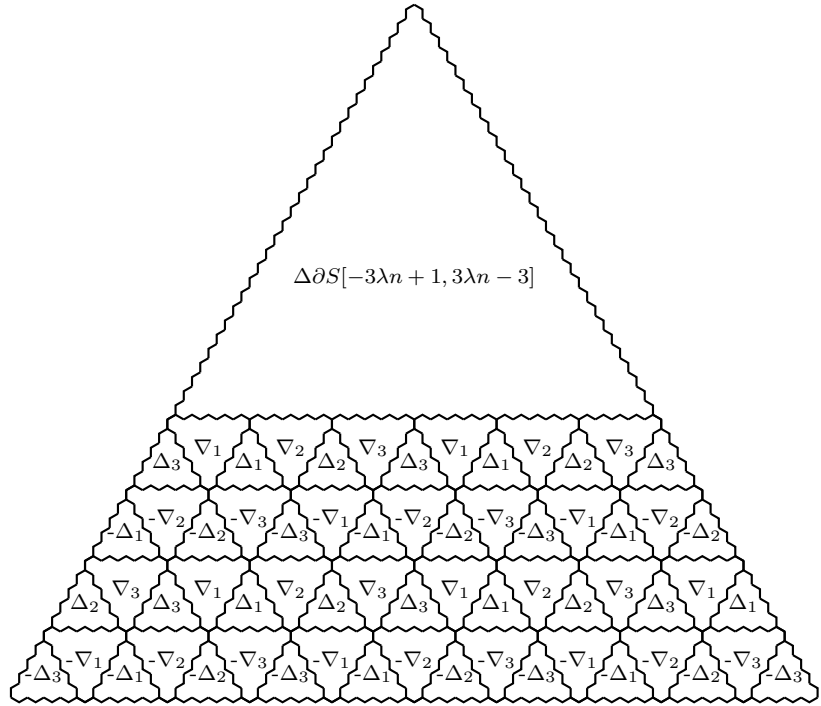
\begin{figure}[!p]
\begin{center}
\begin{pspicture}(10.6088112,9.25)

\multips(0,0)(0.216506351,0){49}{\psline(0,0.0625)(0.108253175,0)(0.216506351,0.0625)}
\multips(0,0)(0.108253175,0.1875){49}{\psline(0,0.0625)(0,0.1875)(0.108253175,0.25)}
\multips(10.3923048,0)(-0.108253175,0.1875){49}{\psline(0.216506351,0.0625)(0.216506351,0.1875)(0.108253175,0.25)}

\multips(2.05681033,2.8125)(1.08253175,0){6}{
\multips(0,0.75)(0.216506351,0){5}{\psline(0,0.1875)(0.108253175,0.25)(0.216506351,0.1875)}
\multips(0,0.75)(0.108253175,-0.1875){5}{\psline(0,0.1875)(0,0.0625)(0.108253175,0)}
\multips(0.866025404,0.75)(-0.108253175,-0.1875){5}{\psline(0.216506351,0.1875)(0.216506351,0.0625)(0.108253175,0)}}

\multips(1.51554446,1.875)(1.08253175,0){7}{
\multips(0,0.75)(0.216506351,0){5}{\psline(0,0.1875)(0.108253175,0.25)(0.216506351,0.1875)}
\multips(0,0.75)(0.108253175,-0.1875){5}{\psline(0,0.1875)(0,0.0625)(0.108253175,0)}
\multips(0.866025404,0.75)(-0.108253175,-0.1875){5}{\psline(0.216506351,0.1875)(0.216506351,0.0625)(0.108253175,0)}}

\multips(0.974278579,0.9375)(1.08253175,0){8}{
\multips(0,0.75)(0.216506351,0){5}{\psline(0,0.1875)(0.108253175,0.25)(0.216506351,0.1875)}
\multips(0,0.75)(0.108253175,-0.1875){5}{\psline(0,0.1875)(0,0.0625)(0.108253175,0)}
\multips(0.866025404,0.75)(-0.108253175,-0.1875){5}{\psline(0.216506351,0.1875)(0.216506351,0.0625)(0.108253175,0)}}

\multips(0.433012702,0)(1.08253175,0){9}{
\multips(0,0.75)(0.216506351,0){5}{\psline(0,0.1875)(0.108253175,0.25)(0.216506351,0.1875)}
\multips(0,0.75)(0.108253175,-0.1875){5}{\psline(0,0.1875)(0,0.0625)(0.108253175,0)}
\multips(0.866025404,0.75)(-0.108253175,-0.1875){5}{\psline(0.216506351,0.1875)(0.216506351,0.0625)(0.108253175,0)}}

\multirput(2.59807621,3.375)(3.24759526,0){2}{\scriptsize$\nabla_1$}
\multirput(3.68060797,3.375)(3.24759526,0){2}{\scriptsize$\nabla_2$}
\multirput(4.76313972,3.375)(3.24759526,0){2}{\scriptsize$\nabla_3$}

\multirput(2.05681033,3.1875)(3.24759526,0){3}{\scriptsize$\Delta_3$}
\multirput(3.13934209,3.1875)(3.24759526,0){2}{\scriptsize$\Delta_1$}
\multirput(4.22187384,3.1875)(3.24759526,0){2}{\scriptsize$\Delta_2$}

\multirput(2.05681033,2.4375)(3.24759526,0){3}{\scriptsize-$\nabla_2$}
\multirput(3.13934209,2.4375)(3.24759526,0){2}{\scriptsize-$\nabla_3$}
\multirput(4.22187384,2.4375)(3.24759526,0){2}{\scriptsize-$\nabla_1$}

\multirput(1.51554446,2.25)(3.24759526,0){3}{\scriptsize-$\Delta_1$}
\multirput(2.59807621,2.25)(3.24759526,0){3}{\scriptsize-$\Delta_2$}
\multirput(3.68060797,2.25)(3.24759526,0){2}{\scriptsize-$\Delta_3$}

\multirput(1.51554446,1.5)(3.24759526,0){3}{\scriptsize$\nabla_3$}
\multirput(2.59807621,1.5)(3.24759526,0){3}{\scriptsize$\nabla_1$}
\multirput(3.68060797,1.5)(3.24759526,0){2}{\scriptsize$\nabla_2$}

\multirput(0.974278579,1.3125)(3.24759526,0){3}{\scriptsize$\Delta_2$}
\multirput(2.05681033,1.3125)(3.24759526,0){3}{\scriptsize$\Delta_3$}
\multirput(3.13934209,1.3125)(3.24759526,0){3}{\scriptsize$\Delta_1$}

\multirput(0.974278579,0.5625)(3.24759526,0){3}{\scriptsize-$\nabla_1$}
\multirput(2.05681033,0.5625)(3.24759526,0){3}{\scriptsize-$\nabla_2$}
\multirput(3.13934209,0.5625)(3.24759526,0){3}{\scriptsize-$\nabla_3$}

\multirput(0.433012702,0.375)(3.24759526,0){4}{\scriptsize-$\Delta_3$}
\multirput(1.51554446,0.375)(3.24759526,0){3}{\scriptsize-$\Delta_1$}
\multirput(2.59807621,0.375)(3.24759526,0){3}{\scriptsize-$\Delta_2$}

\rput(5.3044056,5.625){\scriptsize$\Delta\partial S[-3\lambda n+1,3\lambda n-3]$}

\end{pspicture}
\end{center}
\caption{\label{fig12}The Pascal triangle $\Delta\partial S[-m,m-2]$ for $m\equiv-1\pmod{n}$.}
\end{figure}
\end{proof}

\section{Conclusions and open problems}

In this section, we analyse the results about the generalized Molluzzo problem obtained in this paper and two possible extensions of this work are proposed.

\subsection{Conclusions on the generalized Molluzzo problem}

As listed in Theorem~\ref{thm0} and detailed in Theorem~\ref{thm5}, there exist, for every odd number $n$, infinitely many balanced figures in $\Zn$, for each kind of figure. These results partially solve Problem~\ref{prob2}, the generalized Molluzzo problem. For Steinhaus triangles, since a Steinhaus triangle of order $m$ has cardinality $\binom{m+1}{2}$ and since the set of all the integers $m$ such that the binomial coefficient $\binom{m+1}{2}$ is divisible by $n$ is an union of $2^{\omega(n)}$ classes of integers modulo $n$, where $\omega(n)$ is the number of distinct prime factors of $n$, including the classes of $0$ and $-1$, we have proved, in this paper, that there exist balanced Steinhaus triangles for at least $2/(3.2^{\omega(n)-1})$ of the admissible orders. Particularly, in the case where $n$ is an odd prime power, this proportion becomes $2/3$. In \cite{Chappelon2008}, the author proved that arithmetic progressions with invertible common difference generate balanced Steinhaus triangles for $1/(2^{\omega(n)-1}\beta(n))$ of the admissible orders, where $\beta(n)$ is the order of $2^n$ in the multiplicative quotient group $\left(\Zn\right)^{*}/\{-1,1\}$. This completely solved the Molluzzo problem in $\Z/3^k\Z$ for all $k\geq1$. A new proof of this result, shorter and based on doubly arithmetic triangles, will appear in a forthcoming paper. For Pascal triangles, the proportion of balanced Pascal triangles that we have highlighted is the same: $2/(3.2^{\omega(n)-1})$ for every odd number $n$ and, thus, $2/3$ if $n$ is an odd prime power. Finally, for lozenges, since a lozenge of order $2m-1$ has cardinality $m^2$, the orbit of the universal sequence contains balanced lozenges for all admissible orders in $\Zn$, in the case where $n$ is a square-free odd number. This completely solves Problem~\ref{prob2} for lozenges in the square-free odd case.

\subsection{Additive cellular automata}

Other derivation maps can be considered. For all positive integers $n$ and $r$ and for every $(2r+1)$-tuple of integers $W=(\omega_{-r},\ldots,\omega_0,\ldots,\omega_{r})$, we define the derivation map $\partial_{W}$ by
$$
\partial_{W}\left(a_j\right)_{j\in\Z} = \left(\sum_{k=-r}^{r}\omega_{k}a_{j+k}\right)_{j\in\Z},
$$
for every doubly infinite sequence $\left(a_j\right)_{j\in\Z}$ in $\Zn$. Then, the derivation map $\partial$ of previous sections corresponds to $\partial_{(0,1,1)}$. Now, we naturally wonder, for every $(2r+1)$-tuple of integers $W$, if there exist balanced Steinhaus figures in the \textit{additive cellular automaton} associated with the derivation map $\partial_W$ in $\Zn$.

\begin{prob}
Let $n$ and $r$ be two positive integers and let $W$ be a $(2r+1)$-tuple of integers. Do balanced Steinhaus figures exist in the additive cellular automaton associated with the derivation map $\partial_W$ in $\Zn$?
\end{prob}

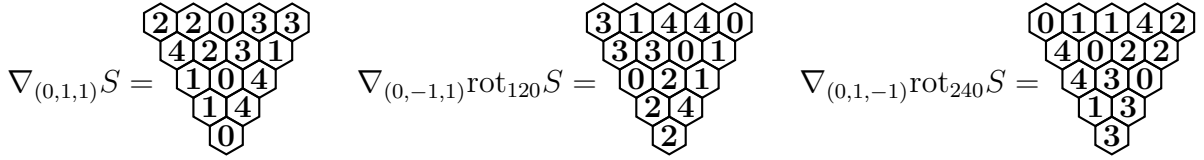
\begin{figure}[!h]
\begin{center}
\begin{tabular}{ccc}
$\nabla_{(0,1,1)}S=$\!\!\!\!\!\begin{tabular}{c}\begin{pspicture}(2.16506351,2)
\multips(0,1.875)(0.433012702,0){5}{\pspolygon(0,0)(0.216506351,0.125)(0.433012702,0)(0.433012702,-0.25)(0.216506351,-0.375)(0,-0.25)}
\multips(0.216506351,1.5)(0.433012702,0){4}{\pspolygon(0,0)(0.216506351,0.125)(0.433012702,0)(0.433012702,-0.25)(0.216506351,-0.375)(0,-0.25)}
\multips(0.433012702,1.125)(0.433012702,0){3}{\pspolygon(0,0)(0.216506351,0.125)(0.433012702,0)(0.433012702,-0.25)(0.216506351,-0.375)(0,-0.25)}
\multips(0.649519053,0.75)(0.433012702,0){2}{\pspolygon(0,0)(0.216506351,0.125)(0.433012702,0)(0.433012702,-0.25)(0.216506351,-0.375)(0,-0.25)}
\multips(0.866025404,0.375)(0.433012702,0){1}{\pspolygon(0,0)(0.216506351,0.125)(0.433012702,0)(0.433012702,-0.25)(0.216506351,-0.375)(0,-0.25)}
\rput(0.216506351,1.75){$\mathbf{2}$}
\rput(0.649519053,1.75){$\mathbf{2}$}
\rput(1.08253175,1.75){$\mathbf{0}$}
\rput(1.51554446,1.75){$\mathbf{3}$}
\rput(1.94855716,1.75){$\mathbf{3}$}
\rput(0.433012702,1.375){$\mathbf{4}$}
\rput(0.866025404,1.375){$\mathbf{2}$}
\rput(1.29903811,1.375){$\mathbf{3}$}
\rput(1.73205081,1.375){$\mathbf{1}$}
\rput(0.649519053,1){$\mathbf{1}$}
\rput(1.08253175,1){$\mathbf{0}$}
\rput(1.51554446,1){$\mathbf{4}$}
\rput(0.866025404,0.625){$\mathbf{1}$}
\rput(1.29903811,0.625){$\mathbf{4}$}
\rput(1.08253175,0.25){$\mathbf{0}$}
\end{pspicture}\end{tabular}
&
$\nabla_{(0,-1,1)}\mathrm{rot}_{120}S=$\!\!\!\!\!\begin{tabular}{c}\begin{pspicture}(2.16506351,2)
\multips(0,1.875)(0.433012702,0){5}{\pspolygon(0,0)(0.216506351,0.125)(0.433012702,0)(0.433012702,-0.25)(0.216506351,-0.375)(0,-0.25)}
\multips(0.216506351,1.5)(0.433012702,0){4}{\pspolygon(0,0)(0.216506351,0.125)(0.433012702,0)(0.433012702,-0.25)(0.216506351,-0.375)(0,-0.25)}
\multips(0.433012702,1.125)(0.433012702,0){3}{\pspolygon(0,0)(0.216506351,0.125)(0.433012702,0)(0.433012702,-0.25)(0.216506351,-0.375)(0,-0.25)}
\multips(0.649519053,0.75)(0.433012702,0){2}{\pspolygon(0,0)(0.216506351,0.125)(0.433012702,0)(0.433012702,-0.25)(0.216506351,-0.375)(0,-0.25)}
\multips(0.866025404,0.375)(0.433012702,0){1}{\pspolygon(0,0)(0.216506351,0.125)(0.433012702,0)(0.433012702,-0.25)(0.216506351,-0.375)(0,-0.25)}
\rput(0.216506351,1.75){$\mathbf{3}$}
\rput(0.649519053,1.75){$\mathbf{1}$}
\rput(1.08253175,1.75){$\mathbf{4}$}
\rput(1.51554446,1.75){$\mathbf{4}$}
\rput(1.94855716,1.75){$\mathbf{0}$}
\rput(0.433012702,1.375){$\mathbf{3}$}
\rput(0.866025404,1.375){$\mathbf{3}$}
\rput(1.29903811,1.375){$\mathbf{0}$}
\rput(1.73205081,1.375){$\mathbf{1}$}
\rput(0.649519053,1){$\mathbf{0}$}
\rput(1.08253175,1){$\mathbf{2}$}
\rput(1.51554446,1){$\mathbf{1}$}
\rput(0.866025404,0.625){$\mathbf{2}$}
\rput(1.29903811,0.625){$\mathbf{4}$}
\rput(1.08253175,0.25){$\mathbf{2}$}
\end{pspicture}\end{tabular}
&
$\nabla_{(0,1,-1)}\mathrm{rot}_{240}S=$\!\!\!\!\!\begin{tabular}{c}\begin{pspicture}(2.16506351,2)
\multips(0,1.875)(0.433012702,0){5}{\pspolygon(0,0)(0.216506351,0.125)(0.433012702,0)(0.433012702,-0.25)(0.216506351,-0.375)(0,-0.25)}
\multips(0.216506351,1.5)(0.433012702,0){4}{\pspolygon(0,0)(0.216506351,0.125)(0.433012702,0)(0.433012702,-0.25)(0.216506351,-0.375)(0,-0.25)}
\multips(0.433012702,1.125)(0.433012702,0){3}{\pspolygon(0,0)(0.216506351,0.125)(0.433012702,0)(0.433012702,-0.25)(0.216506351,-0.375)(0,-0.25)}
\multips(0.649519053,0.75)(0.433012702,0){2}{\pspolygon(0,0)(0.216506351,0.125)(0.433012702,0)(0.433012702,-0.25)(0.216506351,-0.375)(0,-0.25)}
\multips(0.866025404,0.375)(0.433012702,0){1}{\pspolygon(0,0)(0.216506351,0.125)(0.433012702,0)(0.433012702,-0.25)(0.216506351,-0.375)(0,-0.25)}
\rput(0.216506351,1.75){$\mathbf{0}$}
\rput(0.649519053,1.75){$\mathbf{1}$}
\rput(1.08253175,1.75){$\mathbf{1}$}
\rput(1.51554446,1.75){$\mathbf{4}$}
\rput(1.94855716,1.75){$\mathbf{2}$}
\rput(0.433012702,1.375){$\mathbf{4}$}
\rput(0.866025404,1.375){$\mathbf{0}$}
\rput(1.29903811,1.375){$\mathbf{2}$}
\rput(1.73205081,1.375){$\mathbf{2}$}
\rput(0.649519053,1){$\mathbf{4}$}
\rput(1.08253175,1){$\mathbf{3}$}
\rput(1.51554446,1){$\mathbf{0}$}
\rput(0.866025404,0.625){$\mathbf{1}$}
\rput(1.29903811,0.625){$\mathbf{3}$}
\rput(1.08253175,0.25){$\mathbf{3}$}
\end{pspicture}\end{tabular}
\end{tabular}
\end{center}
\caption{\label{fig13}The Steinhaus triangles of $S=(2,2,0,3,3)$, $\mathrm{rot}_{120}S$ and $\mathrm{rot}_{240}S$ in $\Z/5\Z$.}
\end{figure}

Consider the simpler case $W=(0,\omega_0,\omega_1)$ in the sequel and denote by $\nabla_{W}S$ the $W$-Steinhaus triangle and by $\Delta_{W}S$ the $W$-Pascal triangle associated with a finite sequence $S$ in $\Zn$. Then, for every odd number $n$ and for every invertible $d\in\Zn$, the universal sequence $S=\IAP((0,-d,d),(d,-2d,d))$, in $\Zn$, has a $(0,1,1)$-orbit which contains infinitely many balanced $(0,1,-1)$-Steinhaus and Pascal triangles and infinitely many balanced $(0,-1,1)$-Steinhaus and Pascal triangles. Indeed, as illustrated in Figure~\ref{fig13}, the rotation of $120$ degrees defined on the set of finite sequences of length $m\geq1$ in $\Zn$ by
$$
\mathrm{rot}_{120}\left((a_j\right)_{0\leq j\leq m-1})=\left(\sum_{k=0}^{j}\binom{j}{k}a_{m-1-k}\right)_{0\leq j\leq m-1},
$$
induces an isomorphism between $(0,1,1)$-Steinhaus triangles (resp. $(0,1,1)$-Pascal triangles) and $(0,-1,1)$-Steinhaus triangles (resp. $(0,-1,1)$-Pascal triangles), which conserves multiplicity. Similarly, the rotation of $240$ degrees, which assigns to a sequence $(a_j)_{0\leq j\leq m-1}$ of length $m$ in $\Zn$ the sequence
$$
\mathrm{rot}_{240}\left((a_j)_{0\leq j\leq m-1}\right)=\left(\sum_{k=0}^{m-1-j}\binom{m-1-j}{k}a_k\right)_{0\leq j\leq m-1},
$$
induces an isomorphism between $(0,1,1)$-Steinhaus triangles (resp. $(0,1,1)$-Pascal triangles) and $(0,1,-1)$-Steinhaus triangles (resp. $(0,1,-1)$-Pascal triangles), which conserves multiplicity. These sequences can be seen as the right side, for $\mathrm{rot}_{120}(S)$, and the left side, for $\mathrm{rot}_{240}(S)$, of the $(0,1,1)$-Steinhaus triangle $\nabla_{(0,1,1)}S$ associated with $S$.

Finally, since there exist balanced $(0,1,1)$-Steinhaus triangles of order $m$ for every $m\equiv0\pmod{n}$ or $m\equiv-1\pmod{3n}$, in $\Zn$ with $n$ odd, then there exist balanced $(0,-1,1)$ and $(0,1,-1)$-Steinhaus triangles of the same orders in $\Zn$. For an odd prime power $n$, this corresponds to $2/3$ of the admissible orders. Similarly, there exist balanced $(0,-1,1)$ and $(0,1,-1)$-Pascal triangles of order $2m-1$ for every $m\equiv-1\pmod{n}$ or $m\equiv0\pmod{3n}$, in $\Zn$ with $n$ odd. This also corresponds to $2/3$ of the admissible orders, in the case where $n$ is an odd prime power.

\subsection{Steinhaus and Pascal tetrahedra}

In this paper, we have studied balanced Steinhaus figures appearing in the cellular automaton of dimension $1$ that generates the standard Pascal triangle. We may also consider similar figures in higher dimension, in the cellular automaton of dimension $2$ generating the standard Pascal tetrahedron, for instance. Let $n$ be a positive integer and let $S=(a_{i,j})_{i,j\in\Z}$ be a doubly infinite double sequence of terms in $\Zn$. The derived sequence $\partial S$ of $S$ is the sequence defined by $\partial S=\left(a_{i,j}+a_{i,j+1}+a_{i+1,j}\right)_{i,j\in\Z}$ and the orbit of $S$ is the sequence of iterated derived sequences $\mathcal{O}_S=\left(\partial^{k}S\right)_{k\in\N}$. This orbit can also be seen as the $(\N\times{\Z}^2)$-indexed sequence of elements in $\Zn$, defined by
$$
\mathcal{O}_S = \left( \sum_{i'=0}^{k}\sum_{j'=0}^{k-i'}\binom{k}{i',j'}a_{i+i',j+j'}\ \middle|\ i\in\Z,j\in\Z,k\in\N  \right),
$$
where $\binom{k}{i',j'}$ is the trinomial coefficient $\binom{k}{i',j'}=\frac{k!}{i'!j'!(k-i'-j')!}$. The finite orbit of a triangle $T=\left\{a_{i',j'}\middle|0\leq i'\leq m-1,0\leq j'\leq m-1-i'\right\}$, of size $\binom{m+1}{2}$ in $S$, is called the Steinhaus tetrahedron associated with $T$ and of order $\binom{m+1}{2}$. A Steinhaus tetrahedron of order $\binom{m+1}{2}$ has cardinality $\binom{m+2}{3}$. The Molluzzo problem on Steinhaus triangles can then be generalized as follows:

\begin{prob}
Let $n$ be a positive integer. For every $m\geq1$ such that $\binom{m+2}{3}$ is divisible by $n$, does there exist a balanced Steinhaus tetrahedron of order $\binom{m+1}{2}$ in $\Zn$?
\end{prob}

As for Pascal triangles of order $2m-1$ defined from Steinhaus triangles of order $2m-1$, a Pascal tetrahedron of order $\binom{3m-1}{2}$ is a tetrahedron of height $m$, built from the top to the base, that appears in a Steinhaus tetrahedron of order $\binom{3m-1}{2}$. A tetrahedron of order $\binom{3m-1}{2}$ has cardinality $\binom{m+2}{3}$. The Pascal tetrahedron of order $\binom{3m-1}{2}$ associated with the triangle with a $1$ in the middle and $0$ elsewhere corresponds to the first $m$ floors of the standard Pascal tetrahedron modulo $n$. The problem of determining the existence of balanced Pascal tetrahedra in $\Zn$ can be posed.

\begin{prob}
Let $n$ be a positive integer. For every $m\geq1$ such that $\binom{m+2}{3}$ is divisible by $n$, does there exist a balanced Pascal tetrahedron of order $\binom{3m-1}{2}$ in $\Zn$?
\end{prob}

\begin{figure}[!h]
\begin{center}
\begin{tabular}{ccccccc}

\begin{pspicture}(3.03108891,2.75)
\multips(1.29903811,1.5)(0.433012702,0){1}{\pspolygon*[linecolor=Gray](0,0.125)(0,0.375)(0.216506351,0.5)(0.433012702,0.375)(0.433012702,0.125)(0.216506351,0)}

\multips(0,2.25)(0.433012702,0){7}{\pspolygon(0,0.125)(0,0.375)(0.216506351,0.5)(0.433012702,0.375)(0.433012702,0.125)(0.216506351,0)}
\multips(0.216506351,1.875)(0.433012702,0){6}{\pspolygon(0,0.125)(0,0.375)(0.216506351,0.5)(0.433012702,0.375)(0.433012702,0.125)(0.216506351,0)}
\multips(0.433012702,1.5)(0.433012702,0){5}{\pspolygon(0,0.125)(0,0.375)(0.216506351,0.5)(0.433012702,0.375)(0.433012702,0.125)(0.216506351,0)}
\multips(0.649519053,1.125)(0.433012702,0){4}{\pspolygon(0,0.125)(0,0.375)(0.216506351,0.5)(0.433012702,0.375)(0.433012702,0.125)(0.216506351,0)}
\multips(0.866025404,0.75)(0.433012702,0){3}{\pspolygon(0,0.125)(0,0.375)(0.216506351,0.5)(0.433012702,0.375)(0.433012702,0.125)(0.216506351,0)}
\multips(1.08253175,0.375)(0.433012702,0){2}{\pspolygon(0,0.125)(0,0.375)(0.216506351,0.5)(0.433012702,0.375)(0.433012702,0.125)(0.216506351,0)}
\multips(1.29903811,0)(0.433012702,0){1}{\pspolygon(0,0.125)(0,0.375)(0.216506351,0.5)(0.433012702,0.375)(0.433012702,0.125)(0.216506351,0)}

\rput(0.216506351,2.5){$\mathbf{0}$}
\rput(0.649519053,2.5){$\mathbf{4}$}
\rput(1.08253175,2.5){$\mathbf{4}$}
\rput(1.51554446,2.5){$\mathbf{3}$}
\rput(1.94855716,2.5){$\mathbf{1}$}
\rput(2.38156986,2.5){$\mathbf{0}$}
\rput(2.81458256,2.5){$\mathbf{0}$}

\rput(0.433012702,2.125){$\mathbf{2}$}
\rput(0.866025404,2.125){$\mathbf{1}$}
\rput(1.29903811,2.125){$\mathbf{2}$}
\rput(1.73205081,2.125){$\mathbf{0}$}
\rput(2.16506351,2.125){$\mathbf{1}$}
\rput(2.59807621,2.125){$\mathbf{3}$}

\rput(0.649519053,1.75){$\mathbf{4}$}
\rput(1.08253175,1.75){$\mathbf{3}$}
\rput(1.51554446,1.75){$\mathbf{2}$}
\rput(1.94855716,1.75){$\mathbf{0}$}
\rput(2.38156986,1.75){$\mathbf{1}$}

\rput(0.866025404,1.375){$\mathbf{4}$}
\rput(1.29903811,1.375){$\mathbf{3}$}
\rput(1.73205081,1.375){$\mathbf{0}$}
\rput(2.16506351,1.375){$\mathbf{2}$}

\rput(1.08253175,1){$\mathbf{2}$}
\rput(1.51554446,1){$\mathbf{3}$}
\rput(1.94855716,1){$\mathbf{3}$}

\rput(1.29903811,0.625){$\mathbf{0}$}
\rput(1.73205081,0.625){$\mathbf{4}$}

\rput(1.51554446,0.25){$\mathbf{4}$}
\end{pspicture}

&

\begin{pspicture}(2.59807621,2.375)

\multips(1.08253175,1.5)(0.433012702,0){1}{\pspolygon*[linecolor=Gray](0,0.125)(0,0.375)(0.216506351,0.5)(0.433012702,0.375)(0.433012702,0.125)(0.216506351,0)}
\multips(0.866025404,1.125)(0.433012702,0){2}{\pspolygon*[linecolor=Gray](0,0.125)(0,0.375)(0.216506351,0.5)(0.433012702,0.375)(0.433012702,0.125)(0.216506351,0)}

\multips(0,1.875)(0.433012702,0){6}{\pspolygon(0,0.125)(0,0.375)(0.216506351,0.5)(0.433012702,0.375)(0.433012702,0.125)(0.216506351,0)}
\multips(0.216506351,1.5)(0.433012702,0){5}{\pspolygon(0,0.125)(0,0.375)(0.216506351,0.5)(0.433012702,0.375)(0.433012702,0.125)(0.216506351,0)}
\multips(0.433012702,1.125)(0.433012702,0){4}{\pspolygon(0,0.125)(0,0.375)(0.216506351,0.5)(0.433012702,0.375)(0.433012702,0.125)(0.216506351,0)}
\multips(0.649519053,0.75)(0.433012702,0){3}{\pspolygon(0,0.125)(0,0.375)(0.216506351,0.5)(0.433012702,0.375)(0.433012702,0.125)(0.216506351,0)}
\multips(0.866025404,0.375)(0.433012702,0){2}{\pspolygon(0,0.125)(0,0.375)(0.216506351,0.5)(0.433012702,0.375)(0.433012702,0.125)(0.216506351,0)}
\multips(1.08253175,0)(0.433012702,0){1}{\pspolygon(0,0.125)(0,0.375)(0.216506351,0.5)(0.433012702,0.375)(0.433012702,0.125)(0.216506351,0)}

\rput(0.216506351,2.125){$\mathbf{1}$}
\rput(0.649519053,2.125){$\mathbf{4}$}
\rput(1.08253175,2.125){$\mathbf{4}$}
\rput(1.51554446,2.125){$\mathbf{4}$}
\rput(1.94855716,2.125){$\mathbf{2}$}
\rput(2.38156986,2.125){$\mathbf{3}$}

\rput(0.433012702,1.75){$\mathbf{2}$}
\rput(0.866025404,1.75){$\mathbf{1}$}
\rput(1.29903811,1.75){$\mathbf{4}$}
\rput(1.73205081,1.75){$\mathbf{1}$}
\rput(2.16506351,1.75){$\mathbf{0}$}

\rput(0.649519053,1.375){$\mathbf{1}$}
\rput(1.08253175,1.375){$\mathbf{3}$}
\rput(1.51554446,1.375){$\mathbf{2}$}
\rput(1.94855716,1.375){$\mathbf{3}$}

\rput(0.866025404,1){$\mathbf{4}$}
\rput(1.29903811,1){$\mathbf{1}$}
\rput(1.73205081,1){$\mathbf{0}$}

\rput(1.08253175,0.625){$\mathbf{0}$}
\rput(1.51554446,0.625){$\mathbf{0}$}

\rput(1.29903811,0.25){$\mathbf{3}$}

\end{pspicture}

&

\begin{pspicture}(2.16506351,2)

\multips(0.866025404,1.5)(0.433012702,0){1}{\pspolygon*[linecolor=Gray](0,0.125)(0,0.375)(0.216506351,0.5)(0.433012702,0.375)(0.433012702,0.125)(0.216506351,0)}
\multips(0.649519053,1.125)(0.433012702,0){2}{\pspolygon*[linecolor=Gray](0,0.125)(0,0.375)(0.216506351,0.5)(0.433012702,0.375)(0.433012702,0.125)(0.216506351,0)}
\multips(0.433012702,0.75)(0.433012702,0){3}{\pspolygon*[linecolor=Gray](0,0.125)(0,0.375)(0.216506351,0.5)(0.433012702,0.375)(0.433012702,0.125)(0.216506351,0)}

\multips(0,1.5)(0.433012702,0){5}{\pspolygon(0,0.125)(0,0.375)(0.216506351,0.5)(0.433012702,0.375)(0.433012702,0.125)(0.216506351,0)}
\multips(0.216506351,1.125)(0.433012702,0){4}{\pspolygon(0,0.125)(0,0.375)(0.216506351,0.5)(0.433012702,0.375)(0.433012702,0.125)(0.216506351,0)}
\multips(0.433012702,0.75)(0.433012702,0){3}{\pspolygon(0,0.125)(0,0.375)(0.216506351,0.5)(0.433012702,0.375)(0.433012702,0.125)(0.216506351,0)}
\multips(0.649519053,0.375)(0.433012702,0){2}{\pspolygon(0,0.125)(0,0.375)(0.216506351,0.5)(0.433012702,0.375)(0.433012702,0.125)(0.216506351,0)}
\multips(0.866025404,0)(0.433012702,0){1}{\pspolygon(0,0.125)(0,0.375)(0.216506351,0.5)(0.433012702,0.375)(0.433012702,0.125)(0.216506351,0)}

\rput(0.216506351,1.75){$\mathbf{2}$}
\rput(0.649519053,1.75){$\mathbf{4}$}
\rput(1.08253175,1.75){$\mathbf{2}$}
\rput(1.51554446,1.75){$\mathbf{2}$}
\rput(1.94855716,1.75){$\mathbf{0}$}

\rput(0.433012702,1.375){$\mathbf{4}$}
\rput(0.866025404,1.375){$\mathbf{3}$}
\rput(1.29903811,1.375){$\mathbf{2}$}
\rput(1.73205081,1.375){$\mathbf{4}$}

\rput(0.649519053,1){$\mathbf{3}$}
\rput(1.08253175,1){$\mathbf{1}$}
\rput(1.51554446,1){$\mathbf{0}$}

\rput(0.866025404,0.625){$\mathbf{0}$}
\rput(1.29903811,0.625){$\mathbf{1}$}

\rput(1.08253175,0.25){$\mathbf{3}$}

\end{pspicture}

&

\begin{pspicture}(1.73205081,1.625)

\multips(0,1.125)(0.433012702,0){4}{\pspolygon(0,0.125)(0,0.375)(0.216506351,0.5)(0.433012702,0.375)(0.433012702,0.125)(0.216506351,0)}
\multips(0.216506351,0.75)(0.433012702,0){3}{\pspolygon(0,0.125)(0,0.375)(0.216506351,0.5)(0.433012702,0.375)(0.433012702,0.125)(0.216506351,0)}
\multips(0.433012702,0.375)(0.433012702,0){2}{\pspolygon(0,0.125)(0,0.375)(0.216506351,0.5)(0.433012702,0.375)(0.433012702,0.125)(0.216506351,0)}
\multips(0.649519053,0)(0.433012702,0){1}{\pspolygon(0,0.125)(0,0.375)(0.216506351,0.5)(0.433012702,0.375)(0.433012702,0.125)(0.216506351,0)}

\rput(0.216506351,1.375){$\mathbf{0}$}
\rput(0.649519053,1.375){$\mathbf{4}$}
\rput(1.08253175,1.375){$\mathbf{1}$}
\rput(1.51554446,1.375){$\mathbf{1}$}

\rput(0.433012702,1){$\mathbf{0}$}
\rput(0.866025404,1){$\mathbf{1}$}
\rput(1.29903811,1){$\mathbf{1}$}

\rput(0.649519053,0.625){$\mathbf{4}$}
\rput(1.08253175,0.625){$\mathbf{2}$}

\rput(0.866025404,0.25){$\mathbf{4}$}

\end{pspicture}

&

\begin{pspicture}(1.29903811,1.25)

\multips(0,0.75)(0.433012702,0){3}{\pspolygon(0,0.125)(0,0.375)(0.216506351,0.5)(0.433012702,0.375)(0.433012702,0.125)(0.216506351,0)}
\multips(0.216506351,0.375)(0.433012702,0){2}{\pspolygon(0,0.125)(0,0.375)(0.216506351,0.5)(0.433012702,0.375)(0.433012702,0.125)(0.216506351,0)}
\multips(0.433012702,0)(0.433012702,0){1}{\pspolygon(0,0.125)(0,0.375)(0.216506351,0.5)(0.433012702,0.375)(0.433012702,0.125)(0.216506351,0)}

\rput(0.216506351,1){$\mathbf{4}$}
\rput(0.649519053,1){$\mathbf{1}$}
\rput(1.08253175,1){$\mathbf{3}$}

\rput(0.433012702,0.625){$\mathbf{0}$}
\rput(0.866025404,0.625){$\mathbf{4}$}

\rput(0.649519053,0.25){$\mathbf{0}$}

\end{pspicture}

&

\begin{pspicture}(0.866025404,0.875)

\multips(0,0.375)(0.433012702,0){2}{\pspolygon(0,0.125)(0,0.375)(0.216506351,0.5)(0.433012702,0.375)(0.433012702,0.125)(0.216506351,0)}
\multips(0.216506351,0)(0.433012702,0){1}{\pspolygon(0,0.125)(0,0.375)(0.216506351,0.5)(0.433012702,0.375)(0.433012702,0.125)(0.216506351,0)}

\rput(0.216506351,0.625){$\mathbf{0}$}
\rput(0.649519053,0.625){$\mathbf{3}$}

\rput(0.433012702,0.25){$\mathbf{4}$}

\end{pspicture}

&

\begin{pspicture}(0.433012702,0.5)

\multips(0,0)(0.433012702,0){1}{\pspolygon(0,0.125)(0,0.375)(0.216506351,0.5)(0.433012702,0.375)(0.433012702,0.125)(0.216506351,0)}

\rput(0.216506351,0.25){$\mathbf{2}$}

\end{pspicture}

\end{tabular}
\end{center}
\caption{A Steinhaus tetrahedron in $\Z/5\Z$, with a Pascal tetrahedron in gray.}
\end{figure}
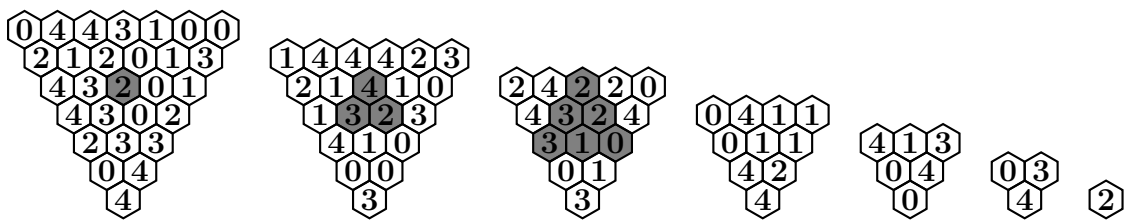

\section*{Acknowledgments}

I would like to thank Shalom Eliahou for his unfailing support and his enthusiasm for this work. I also thank the anonymous referees for their helpful comments.

\end{document}